\documentclass[11pt]{article} 

\usepackage{tgpagella}
\usepackage[T1]{fontenc}

\usepackage[disable]{todonotes}
\usepackage{graphicx}
\usepackage[square,sort,comma,numbers]{natbib}
\usepackage{etex}
\usepackage{enumitem}
\setlist{parsep = -0em, itemsep = 0.25em}
\usepackage[utf8]{inputenc}
\usepackage[dvips,letterpaper,margin=1in]{geometry}
\usepackage{amsmath}
\usepackage{amssymb}
\usepackage{amsthm}
\usepackage{bm}
\usepackage{colonequals}
\usepackage{comment}
\usepackage{graphicx}
\usepackage{subcaption}
\usepackage{xcolor}
\usepackage{tikz}
\usetikzlibrary{shadows}
\usepackage{authblk}
\usepackage[utf8]{inputenc} 
\usepackage[T1]{fontenc}
\usepackage{url}            
\usepackage{booktabs}       
\usepackage{amsfonts}       
\usepackage{nicefrac}      
\usepackage{microtype}      
\usepackage{mathtools}
\usepackage{dsfont}

\usepackage{hyperref}
\hypersetup{
    colorlinks = true,
    linkcolor = blue,
    urlcolor  = black,
    citecolor = magenta,
    anchorcolor = blue
}

\usepackage{amsmath,amsthm,amssymb}
\usepackage{algorithm,algpseudocode}
\usepackage{wrapfig}
\usepackage[toc,page,header]{appendix}
\usepackage{minitoc}

\usepackage{etoolbox}

\makeatletter
\patchcmd{\@maketitle}{\LARGE \@title}{\fontsize{16}{19.2}\selectfont\@title}{}{}
\makeatother

\usepackage{esint}

\newcommand\blfootnote[1]{
  \begingroup
  \renewcommand\thefootnote{}\footnote{#1}
  \addtocounter{footnote}{-1}
  \endgroup
}

\newtheorem{theorem}{Theorem}[section]
\newtheorem{proposition}[theorem]{Proposition}
\newtheorem{lemma}[theorem]{Lemma}
\newtheorem{remark}[theorem]{Remark}
\newtheorem{assumption}[theorem]{Assumption}
\newtheorem{corollary}[theorem]{Corollary}

\DeclareMathOperator{\supp}{supp}

\DeclareMathOperator{\argmin}{argmin}

\def\RR{\mathbb{R}} \def\NN{\mathbb{N}} 
\def\EE{\mathbb{E}}\def\PP{\mathbb{P}}

\title{A Dynamical Central Limit Theorem for Shallow Neural Networks}

\author[ a]{Zhengdao Chen}
\author[ a, b]{Grant M. Rotskoff}
\author[ a, c]{Joan Bruna}
\author[ a]{Eric Vanden-Eijnden}

\affil[a]{Courant Institute of Mathematical Sciences, New York
  University, New York}
\affil[b]{Department of Chemistry, Stanford University, California}
\affil[c]{Center for Data Science, New York University, New York}

\def\thetab{\boldsymbol{\theta}}
\def\Thetab{\boldsymbol{\Theta}}

\def\xib{\boldsymbol{\xi}}
\def\ab{\boldsymbol{a}}
\def\bb{\boldsymbol{b}}

\def\ub{\boldsymbol{u}}
\def\xb{\boldsymbol{x}}

\def\zb{\boldsymbol{z}}

\def\Tb{\boldsymbol{T}}

\def\E{\mathbb{E}}
\def\R{\mathbb{R}}

\newcommand\restr[2]{{
  \left.\kern-\nulldelimiterspace
  #1
  \vphantom{\big|}
  \right|_{#2}
  }}

\def\RR{\mathbb{R}} \def\NN{\mathbb{N}} 
\def\EE{\mathbb{E}}\def\PP{\mathbb{P}}

\begin{document}
\maketitle
\blfootnote{Correspondence to: \href{mailto:zc1216@nyu.edu}{zc1216@nyu.edu}, \href{mailto:rotskoff@stanford.edu}{rotskoff@stanford.edu}, \href{mailto:bruna@cims.nyu.edu}{bruna@cims.nyu.edu} and \href{mailto:eve2@cims.nyu.edu}{eve2@cims.nyu.edu}.}

\begin{abstract}
Recent theoretical works have characterized the dynamics of wide shallow neural networks trained via gradient descent in an asymptotic mean-field limit when the width tends towards infinity. At initialization, the random sampling of the parameters leads to deviations from the mean-field limit dictated by the classical Central Limit Theorem (CLT). However, since gradient descent induces correlations among the parameters,  it is of interest to analyze how these fluctuations evolve. Here, we use a dynamical CLT to prove that the asymptotic fluctuations around the mean limit remain bounded in mean square throughout training. The upper bound is given by a Monte-Carlo resampling error, with a variance that that depends on the $2$-norm of the underlying measure, which also controls the generalization error. This motivates the use of this $2$-norm as a regularization term during training. Furthermore, if the mean-field dynamics converges to a measure that interpolates the training data, we prove that the asymptotic deviation eventually vanishes in the CLT scaling. 
We also complement these results with numerical experiments.
\end{abstract}

\doparttoc
\faketableofcontents

\part{}
\section{Introduction}
Theoretical analyses of neural networks aim to understand their computational and statistical advantages seen in practice. On the computation side, the training of neural networks often succeed despite being a non-convex  optimization problem known to be hard in certain settings \cite{livni2014computational,goel2020superpolynomial,diakonikolas2020algorithms}. On the statistics side, neural networks often generalize well despite having large numbers of parameters \cite{zhang2016understanding, belkin_reconciling_2019}. In this context, the notion of \emph{over-parametrization} has been useful, by providing insights into the optimization and generalization properties as the network widths tend to infinity \cite{jacot2018neural, du2018gradient, allen2019convergence, arora2019fine, soltanolkotabi2018theoretical, JMLR:v20:18-674, lee2017deep}. In particular, under appropriate scaling, one can view shallow (a.k.a. single-hidden-layer or two-layer) networks as interacting particle systems that admit a mean-field limit. Their training dynamics can then be studied as Wasserstein Gradient Flows \cite{mei2018mean, rotskoff2018parameters, chizat2018global, sirignano2018dgm}, leading to global convergence guarantees in the mean-field limit under certain assumptions. 
On the statistics side, such an approach lead to powerful generalization guarantees for learning high-dimensional functions with hidden low-dimensional structures, as compared to learning in Reproducing Kernel Hilbert Spaces (RKHS) \cite{bach2017breaking, ghorbani2020neural}. 
However, since ultimately we are concerned with neural networks of finite width, it is key to study the deviation of finite-width networks from their infinite-width limits, and how it scales with the width $m$. 
At the random initial state, neurons do not interact and therefore a standard Monte-Carlo (MC) argument shows that the fluctuations in the underlying measure scale as $m^{-1/2}$, which we refer to as the Central Limit Theorem (CLT) scaling. 
As optimization introduces complex dependencies among the parameters, the key question is to understand how the fluctuation evolves during training. To make this investigation tractable, we aim to obtain insight on an asymptotic scale as the width grows, and focus on the evolution in time. 
An application of Gr\"onwall's inequality  shows that this asymptotic deviation remains bounded at all finite time \cite{mei2019mean}, but the dependence on time is exponential, making it difficult to assess the long-time behavior.

The main focus of this paper is to investigate this question in-depth, by analyzing the interplay between the deviations from the mean-field limit and the gradient flow dynamics. 
First, we prove a dynamical CLT to capture how the fluctuations away from the mean-field limit evolve as a function of training time to show that the fluctuations remain on the initial $m^{-1/2}$-scale for all finite times.
Next, we examine the long-time behavior of the fluctuations, proving that, in several scenarios, the long-time fluctuations are controlled by the error of Monte-Carlo resampling from the limiting measure.
We focus on two main setups relevant for supervised learning and scientific computing: the unregularized case with global convergence of mean-field gradient flows to minimizers that interpolate the data, and the regularized case where the limiting measure has atomic support and is nondegenerate. In the former setup, we prove particularly that the fluctuations eventually vanish in the CLT scaling.
These asymptotic predictions are complemented by empirical results in a teacher-student model.

\paragraph{Related Works:}
This paper continues the line of work initiated in \cite{mei2018mean, chizat2018global, rotskoff2018parameters, sirignano2018dgm} that studies optimization of over-parameterized shallow neural networks under the mean-field scaling.
Global convergence for the unregularized setting is discussed in \cite{mei2018mean,mei2019mean,sirignano2018dgm,rotskoff2018parameters}.
In the regularized setting, \cite{chizat2018global} establishes global convergence in the mean-field limit under specific homogeneity conditions on the neuron activation.
Other works that study asymptotic properties of wide neural networks include \cite{ghorbani2019limitations, geiger2019disentangling, bai2019beyond, dyer2019asymptotics, hanin2019finite, huang2019dynamics, woodworth2020kernel, luo2020phase, adlamneural}, notably investigating the transition between the so-called \emph{lazy} and \emph{active} regimes \cite{chizat2019lazy}, corresponding respectively to linear versus nonlinear learning. Our focus is on the dynamics under the mean-field scaling, which encompasses the active, nonlinear regime. 

A relevant work concerning the sparse optimization of measures is \cite{chizat2019sparse}, where under a different metric for gradient flow and additional assumptions on the nature of the minimizer, it can be established that fluctuations vanish for sufficiently large $m$. 
Our results are only asymptotic in $m$ but apply to broader settings in the context of shallow neural networks.
Concerning the next-order deviations of finite neural networks from their mean-field limit, \cite{rotskoff2018parameters} show that the scale of fluctuations is below that of MC resampling for unregularized problems using non-rigorous arguments. \cite{sirignano2020mean} provides a CLT for the fluctuations at finite time under stochastic gradient descent (SGD) and proves that the fluctuations decay in time in the case where there is a single critical point in the parameter space. Our focus is on the long-time behavior of the fluctuations in more general settings. Another relevant topic is the propagation of chaos in McKean-Vlasov systems, which study the deviations of randomly-forced interacting particle systems from their infinite-particle limits \citep{braun1977vlasov, sznitman1991topics, spohn2012large, baladron2011mean}. In particular, a line of work provides uniform-in-time bounds to the fluctuations in various settings \citep{del2000branching, cortez2016uniform, salem2018gradient, salhi2018uniform, durmus2018elementary}, but the conditions are not applicable to shallow neural networks. Concurrently to our work, \cite{de2020quantitative} studies quantitative propagation of chaos of shallow neural networks trained by SGD, but the bound grows exponentially in time, and therefore cannot address the long-time behavior of the fluctuations.

Learning with neural networks exhibits the phenomenon that generalization error can decrease with the level of overparameterization \citep{belkin_reconciling_2019, spigler2018jamming}. \cite{neal2018modern} proposes a bias-variance decomposition that contains a variance term initialization in optimization. They show in experiments that this term decreases as the width of the network increases, and justifies this theoretically under the strong assumption that model parameters remain Gaussian-distributed in the components that are irrelevant for the task, which does not hold in the scenario we consider, for example. \cite{geiger2020scaling} provides scaling arguments for the dependence of this term on the width of the network. Our work provides a more rigorous analysis of the dependence of this term on the width of the network and training time.

\section{Background}
\label{sec:meanfield}
\subsection{Shallow Neural Networks and the Integral Representation}
On a data space $\Omega \subseteq \R^d$, we consider parameterized models of the following form
\begin{equation}
\label{eq:basic}
    f^{(m)}(\xb) = \frac{1}{m} \sum_{i=1}^m \varphi(\thetab_i, \xb),
\end{equation}
where $\xb \in \Omega$,  $\{ \thetab_i \}_{i=1}^m \subseteq D$ is the set of model parameters, and $\varphi:D\times \Omega \to\R$ is the activation function. Of particular interest are shallow neural network models, which admit a more specific form:
\begin{assumption}[Shallow neural networks setting]
\label{ass:shallownn}
$D = \mathbb{R} \times \hat{D}$, $\thetab = (c, \zb) \in D$, and $\varphi(\thetab, \xb) = c \hat{\varphi}(\zb, \xb)$ with $\hat \varphi: \hat D\times \Omega \to\R$. Thus, \eqref{eq:basic} can be rewritten as $f^{(m)}(\xb) = \frac{1}{m} \sum_{i=1}^m c_i \hat{\varphi}(\zb_i, \xb)$.
\end{assumption}
\noindent
As many of our results hold for general models of the form \eqref{eq:basic}, we will invoke Assumption \ref{ass:shallownn} only when needed. We shall also assume the following:
\begin{assumption}
\label{ass:unit_1}
$\Omega$ is compact; $D$ is a Euclidean space (or a subset thereof);  
$\varphi(\thetab, \xb)$ is twice differentiable in $\thetab$; $\nabla_{\thetab} \nabla_{\thetab} \varphi(\thetab, \xb)$ is Lipschitz in $\thetab$, uniformly in $\xb$.
\end{assumption}
\noindent
The regularity assumptions are standard in the literature \citep{chizat2019sparse, lancellotti2009fluctuations, braun1977vlasov}.
We note that they are not satisfied by ReLU units (i.e., $\hat{\varphi}(\zb, \xb) = \max\{0, \langle \ab, \xb \rangle + b \}$, where $\zb = ( \ab, b )^\intercal$, with $\ab \in \mathbb{R}^d$ and $b \in \mathbb{R}$), though prior work \citep{chizat2018global, chizat2020implicit} has considered differentiable approximations of these models.

As observed in ~\cite{mei2018mean,chizat2018global, rotskoff2018parameters, sirignano2018dgm, e2019machine}, a model of the form \eqref{eq:basic} can be expressed in integral form in terms of a probability measure over $D$ as $f^{(m)} = f[\mu^{(m)}]$, where we define
\begin{equation}
\label{eq:integralrep}
     f[\mu](\xb) = \int_D \varphi(\thetab, \xb) \mu(d\thetab)~,
\end{equation}
and $\mu^{(m)}$ is the empirical measure of the parameters $\{ \thetab_i \}_{i=1}^m$:
\begin{equation}
\label{eq:empmeasure}
    \mu^{(m)}(d \thetab) = \frac{1}{m} \sum_{i=1}^m \delta_{\thetab_i}(d \thetab)~.
\end{equation}  

Suppose we are given a dataset $\{(\xb_l, y_l)\}_{l = 1}^n$, which can be represented by an empirical data measure $\hat \nu = \tfrac{1}{n} \sum_{l=1}^n \delta_{\xb_l}$, and $y_l = f_*(\xb_l)$ are generated by an target function $f_*$ that we wish to estimate using least-squares regression. A canonical approach to this regression task is to consider an Empirical Risk Minimization (ERM) problem of the form 
\begin{equation}
\label{eq:ermmeasure}
    \min_{\mu \in \mathcal{P}(D)} \mathcal{L}(\mu) \quad \text{with} \quad  \mathcal{L}(\mu)  := 
    \tfrac{1}{2} \left\| f[\mu] - f_*\right\|_{\hat\nu} ^2 + \lambda \int_D r(\thetab) \mu(d \thetab) ~.
\end{equation}
where $\mathcal{P}(D)$ is the space of probability measures on $D$, $\| f - f_*\|_{\hat \nu}^2 := \int_\Omega |f(\xb) - f_*(\xb)|^2\hat \nu(d\xb)$ denotes the $L_2$ function reconstruction error averaged over the data, and $\lambda \int_D r(\thetab) \mu(d \thetab)$ is some optional regularization term. While we can allow $r$ to be a general convex function, in Section~\ref{sec:var_norm} we will motivate a choice of $r$ in the shallow neural networks setting that is related to the variation norm \cite{bach2017breaking} or Barron norm \cite{ma2019barron} of functions.


\subsection{Approximation and Optimization with a Finite Number of Neurons} 
Integral representations with a probability measure such as those defined in \eqref{eq:integralrep} are amenable to efficient approximation in high dimensions via Monte-Carlo sampling. Namely, if the parameters~$\thetab_i$ in $f^{(m)}$ are drawn i.i.d. from an underlying measure $\mu$ on $D$, then by the Law of Large Numbers (LLN), the resulting empirical measure $\mu^{(m)}$ converges $\mu$ almost surely, and moreover, 
\begin{equation}
    \label{eq:MCbounda}
    \E_{\mu^{(m)}} \|f[\mu^{(m)}] - f[\mu] \|_{\hat\nu}^2 = \frac{1}{m} \left( \int_D \| \varphi(\thetab, \cdot)\|_{\hat{\nu}}^2 \mu (d \thetab) - \|f[\mu]\|_{\hat\nu}^2\right)~,
\end{equation}
Such a Monte-Carlo estimator showcases the benefit of normalized integral representations for high-dimensional approximation, as the ambient dimension appears in the rate of approximation only through the term $\int_D \| \varphi(\thetab, \cdot)\|_{\hat{\nu}}^2 \mu (d \thetab)$. In the case of shallow neural networks, this is connected to the variation norm or Barron norm of the function we wish to approximate \citep{bach2017breaking, ma2019barron} (see Section \ref{sec:var_norm} for details).  

While the Monte-Carlo sampling strategy above can be seen as a `static' approximation of a function representable as \eqref{eq:integralrep}, it also gives rise to an efficient algorithm to optimize (\ref{eq:ermmeasure}). Indeed, in terms of the empirical distribution $\mu^{(m)}$, the loss $\mathcal{L}(\mu^{(m)})$ becomes a function of the parameters $\{ \thetab_i \}_{i=1}^{m}$, which we can seek to minimize by adjusting the parameters:
\begin{equation}
\label{eq:pop_loss_n}
{L}(\thetab_1, \dots, \thetab_m) =\frac{1}{2}\|f^{(m)}- f_*\|_{\hat\nu}^2 + \frac{\lambda}{m}\sum_{i=1}^m r(\thetab_i)~.
\end{equation}
In the shallow neural network setting, with suitable choices of the function $r$, the regularization term corresponds to \emph{weight decay} over the parameters. 
\subsection{From Particle to Wasserstein Gradient Flows}
Expanding \eqref{eq:pop_loss_n}, we get
\begin{equation}
\label{eq:loss}
    {L}(\thetab_1, \dots, \thetab_m) = C_{f_*} - \frac{1}{m} \sum_{i=1}^m F(\thetab_i) + \frac{1}{2 m^2} \sum_{i, j = 1}^m K(\thetab_i, \thetab_j),
\end{equation}
where we have defined $C_f = \frac{1}{2} \|f\|_{\hat{\nu}}^2$, and 
\begin{equation}
     F(\thetab) = \int_\Omega f_*(\xb) \varphi(\thetab, \xb) \hat{\nu}(d \xb) -\lambda r(\thetab), \qquad K(\thetab, \thetab') = \int_\Omega \varphi(\thetab, \xb)\varphi(\thetab', \xb) \hat{\nu}(d \xb)~.
\end{equation}
Performing GD on  ${L}$ amounts to discretizing in time the following ODE system that governs the evolution of 
for $\{ \thetab_i \}_{i=1}^m$: 
\begin{equation}
\label{eq:gd_particle}
\begin{split}
    \dot{\thetab}_i = -m \partial_{\thetab_i} L(\thetab_1\dots\thetab_m)=& \nabla F(\thetab_i) - \frac{1}{m} \sum_{j=1}^m \nabla K(\thetab_i, \thetab_j) =: - \nabla V(\thetab_i,\mu_t^{(m)}).
\end{split}
\end{equation}
where we defined the potential
\begin{equation}
\label{eq:Vdef}
    V(\thetab,\mu) = -F(\thetab) + \int_D K(\thetab, \thetab') \mu(d \thetab')~.
\end{equation}
Heuristically, the `particles' $\thetab_i$ perform GD according to the potential $V(\thetab,\mu_t^{(m)})$
which itself evolves, depending on the particles positions through their empirical measure. Such dynamics can also be expressed in terms of the empirical measure via the \emph{continuity equation}:
\begin{equation}
\label{eq:mv_n}
    \partial_t \mu_t^{(m)} = \nabla \cdot (\nabla V (\thetab, \mu^{(m)}_t) \mu^{(m)}_t)
\end{equation}
This equation should be understood in the weak sense by testing it against continuous functions $\chi:D\to\RR$, and it can be interpreted as the gradient flow on the loss defined in (\ref{eq:ermmeasure}) under the 2-Wasserstein metric
~\cite{chizat2018global,rotskoff2018parameters,mei2018mean,sirignano2018dgm}. This insight provides powerful analytical tools to understand convergence properties, by considering the mean-field limit when $m \to \infty$. 

\subsection{Law of Large Numbers and Mean-Field Gradient Flow}
From now on, we assume that the particle gradient flow is initialized in the following way:
\begin{assumption}
\label{ass:iid_init}  The ODE \eqref{eq:gd_particle} is solved for the initial condition $\thetab_i(0)=\thetab_i^0$, with
$\thetab_i^0$ drawn i.i.d. from a compactly supported measure $\mu_0\in\mathcal{P}(D)$ for each $i=1,\ldots,m$.  Hence, $\mu^{(m)}_0(d\thetab) = \tfrac{1}{m} \sum_{i=1}^m  \delta_{\thetab_i^0}(d \thetab)$. 
\end{assumption}
\noindent We use $\PP_0$ to denote the probability measure associated with the set $\{\thetab_i^0\}_{i\in\NN}$ with each $\thetab_i^0$ drawn i.i.d. from $\mu_0$, and use $\EE_0$ to denote the expectation under $\PP_0$. The Law of Large Numbers (LLN) indicates
that $\PP_0$-almost surely, $\mu_t^{(m)} \rightharpoonup \mu_t$ as $m \to \infty$, where $\mu_t$ satisfies the mean-field gradient flow~\cite{rotskoff2018neural, chizat2018global, mei2018mean, sirignano2020mean_lln}:
\begin{equation}
\label{eq:mf}
    \partial_t \mu_t = \nabla \cdot (\nabla V(\thetab, \mu_t) \mu_t)~, \qquad \mu_{t=0}=\mu_0~.
\end{equation}
The solution to this equation can be understood via the representation formula 
\begin{equation}
    \label{eq:weakn}
    \int_D \chi(\thetab) \mu_t(d \thetab) = \int_D \chi(\Thetab_t(\thetab)) \mu_0(d \thetab)~,
\end{equation}
where $\chi$ is a continuous test function $\chi: D \to \mathbb{R}$ and $\Thetab_t: D \to D$ is the \emph{characteristic flow} associated with \eqref{eq:mv_n}, which in direct analogy with \eqref{eq:gd_particle} solves
\begin{equation}
    \label{eq:flow}
    \dot{\Thetab}_t(\thetab) = - \nabla V(\Thetab_t(\thetab), \mu_t), \qquad \Thetab_0(\thetab) = \thetab~.
\end{equation}
Using expression~\eqref{eq:Vdef} for $V$ as well as \eqref{eq:weakn}, this equation can be written in closed form explicitly as
\begin{equation}
    \label{eq:flowexplicit}
    \dot{\Thetab}_t(\thetab) = \nabla F(\Thetab_t(\thetab)) - \int_D \nabla K(\Thetab_t(\thetab),\Thetab_t(\thetab'))\mu_0(d\thetab'), \qquad \Thetab_0(\thetab) = \thetab~.
\end{equation}
It is easy to see that this equation is itself a gradient flow since it is the continuous-time limit of a proximal scheme (mirror descent), which we state as:
\begin{proposition}
  \label{th:prop}
  Given $\bar\Thetab_0(\thetab) = \thetab$ and $\tau>0$, for $p\in \NN$ let $\Thetab_{p\tau}$ be specified via
  \begin{equation}
      \label{eq:proxi1}
      \bar\Thetab_{p\tau} \in \argmin \left( \frac1{2\tau} \|\Thetab-\bar\Thetab_{(p-1)\tau}\|_0^2 + \mathcal{E}(\Thetab)~,  \right)
  \end{equation}
  where we defined
  \begin{equation}
      \label{eq:norm0}
      \|\Thetab\|_0^2 = \int_D |\Thetab(\thetab)|^2 \mu_0(d\thetab)
  \end{equation}
  and 
  \begin{equation}
      \label{eq:E0}
      \mathcal{E}(\Thetab) = - \int_D F(\Thetab(\thetab)) \mu_0(d\thetab) + \frac12 \int_D K(\Thetab(\thetab),\Thetab(\thetab')) \mu_0(d\thetab)\mu_0(d\thetab')~.
  \end{equation}
  Then 
  \begin{equation}
      \lim_{\tau\to0}\bar\Thetab_{\lfloor t/\tau\rfloor\tau} = \Thetab_t \qquad \text{$\mu_0$-almost surely}~,
  \end{equation}
  where $\Thetab_t$ solves~\eqref{eq:flowexplicit}.
\end{proposition}

\subsection{Long-Time Properties of the Mean-Field Gradient Flow}

In the shallow neural networks setting, a series of earlier works \cite{chizat2018global,rotskoff2018parameters,mei2018mean,sirignano2018dgm} has established that under certain assumptions $\mu_t$ will converge to a global minimizer of the loss functional $\mathcal{L}$. In particular, \cite{chizat2018global} studies global convergence for the regularized loss $\mathcal{L}$ under homogeneity assumptions on $\hat{\varphi}$, and \cite{rotskoff2019global} considers modified dynamics using \emph{double-lifting}. Here, to study the long time behavior of the fluctuations, we will often work with the following weaker assumptions:
\begin{assumption}
\label{ass:bounded_supp}
The solution to \eqref{eq:flowexplicit} exists for all time, and has a limit:
\begin{equation}
    \label{eq:limflow}
    \Thetab_t \to \Thetab_\infty \quad \text{$\mu_0$-almost surely as $t\to\infty$}.
\end{equation}
\end{assumption}
\begin{assumption}
\label{ass:limflow_local_min}
The limiting $\Thetab_\infty$ is a local minimizer of~\eqref{eq:E0}.
\end{assumption}

\noindent With these assumptions, we have
\begin{proposition}
  \label{th:ltflow}
  Under Assumptions~\ref{ass:iid_init} and \ref{ass:bounded_supp}, we have 
\begin{equation}
\label{eq:bounded_supp}
    \cup_{t \geq 0} \supp {\mu_t} = \cup_{t \geq 0} \{\Thetab_t(\thetab)\,:\, \thetab\in\supp \mu_0\} \ \ \text{is compact,}
\end{equation}  
and $\mu_t \rightharpoonup \mu_\infty$ weakly as $t\to\infty$, with $\mu_\infty$ satisfying
  \begin{equation}
    \label{eq:weakninfty}
    \int_D \chi(\thetab) \mu_\infty(d \thetab) = \int_D \chi(\Thetab_\infty(\thetab)) \mu_0(d \thetab),
\end{equation}
for all continuous test function $\chi: D \to \mathbb{R}$. Additionally, if Assumption~\ref{ass:limflow_local_min} also holds, then
\begin{equation}
    \label{eq:Vlim}
    \nabla \nabla V(\Thetab_\infty(\thetab),\mu_\infty) \ \ \text{is positive semidefinite for $\mu_0$-almost all $\thetab$}
\end{equation}
\end{proposition}
\noindent
We prove this proposition in Appendix \ref{app_local_min}. Here, $\nabla \nabla V(\Thetab_\infty(\thetab),\mu_\infty)$ denotes
\begin{equation}
    \label{eq:Vlimdef}
    \nabla \nabla V(\Thetab_\infty(\thetab),\mu_\infty) = -\nabla \nabla F(\Thetab_\infty(\thetab))+ \int_D \nabla \nabla K(\Thetab_\infty(\thetab),\Thetab_\infty(\thetab')) \mu_0(d\thetab')~,
\end{equation} 
which will become useful in Section~\ref{sec:g_infty} when we analyze the long time properties of the fluctuations around the mean-field limit. 

\begin{remark}
Assumptions \ref{ass:bounded_supp} and \ref{ass:limflow_local_min} impose conditions on the initial measure $\mu_0$ \cite{rotskoff2018parameters, mei2018mean, chizat2018global}. While the convergence of gradient flows in finite-dimensional Euclidean space to local minimizers is guaranteed under mild assumptions \cite{smale1963stable, lee2017first}, its infinite-dimensional counterpart, Assumption~\ref{ass:limflow_local_min}, may require further technical assumptions, left for future study. Also, while Assumption~\ref{ass:bounded_supp} implies that $\mu_\infty$ is a stationary point of~\eqref{eq:mf},  Assumption~\ref{ass:limflow_local_min} does not imply that $\mu_\infty$ minimizes $\mathcal{L}$.
\end{remark} 

\section{Fluctuations from Mean-Field Gradient Flow}
The main goal of this section is to characterize the deviations of finite-particle shallow networks from their mean-field evolution, by first deriving an estimate for $f^{(m)}_t-f_t$ for $t \ge
0$ (Section \ref{sec:formalCLT}), and then analyzing its long-time properties (Section \ref{sec:g_infty}).
In Section \ref{sec:var_norm}, we then motivate a choice of the regularization term in \eqref{eq:ermmeasure} that controls the bound on the long-time fluctuations derived in Section \ref{sec:g_infty}, and which
is also connected to generalization via the variation norm \cite{bach2017breaking} or Barron norm \cite{ma2019barron} of functions.

\subsection{A Dynamical Central Limit Theorem}
\label{sec:formalCLT}
Let us start by defining
\begin{equation}
  g^{(m)}_t := m^{1/2} \big(f^{(m)}_t - f_t\big)~.
\end{equation}
By the static Central Limit Theorem (CLT) we know that, if we draw the
initial values of the parameters $\thetab_i$ independently from
$\mu_0$ as specified in Assumption~\ref{ass:iid_init}, $g^{(m)}_{0}$ has a limit as $m\to\infty$, leading
to estimates similar to~\eqref{eq:MCbounda} with $\mu^{(m)}$ and $\mu$
replaced by the initial $\mu_0^{(m)}$ and $\mu_0$, respectively. For
$t > 0$, however, this estimate is not preserved by the gradient flow:
the static CLT no longer applies and needs to be replaced by a
dynamical
variant~\cite{braun1977vlasov,sznitman1991topics,spohn2012large, sirignano2020mean}. Next,
we derive this dynamical CLT in the context of neural network
optimization.

To this end
let us define the discrepancy measure $\omega^{(m)}_t$ such that
\begin{equation}
    \int_D \chi(\thetab) \omega^{(m)}_t(d\thetab) 
    := m^{1/2}  \int_D \chi(\thetab) \left(\mu^{(m)}_t (d\thetab)- \mu_t(d\thetab) \right)~,
\end{equation}
for any continuous test function $\chi: D\to \RR$. We can then represent $g^{(m)}_t$ in terms of $\omega^{(m)}_t$ as
\begin{equation}
  \label{eq:g_vs_omega}
    g^{(m)}_t = \int_D \varphi(\thetab,\cdot) \omega^{(m)}_t (d\thetab)~.
\end{equation}
Hence, we will first establish how the limit of $\omega^{(m)}_t$ as $m
\to \infty$ evolves over time. This can be done by noting  that the representation formula~\eqref{eq:weakn} implies that
\begin{equation}
    \label{eq:omegat_t_n}
    \int_D \chi(\thetab) \omega^{(m)}_t(d\thetab)=
    m^{1/2}  \int_D \left(\chi(\Thetab^{(m)}_t(\thetab))\mu^{(m)}_0(d\thetab) -  \chi(\Thetab_t(\thetab))\mu_0(d\thetab) \right)~,
\end{equation}
where $\Thetab^{(m)}_t$ solves \eqref{eq:flowexplicit} with $\mu_0$ replaced by $\mu_0^{(m)}$. Defining
\begin{equation}
    \label{eq:identA2}
    \Tb_t^{(m)} (\thetab) = m^{1/2} \big ( \Thetab_t^{(m)}(\thetab) - \Thetab_t(\thetab) \big )~,\end{equation}
we can write~\eqref{eq:omegat_t_n} as
\begin{equation}
    \label{eq:omegat_t_n2}
    \begin{split}
        \int_D \chi(\thetab) \omega_t^{(m)} (d \thetab) 
        = & \int_D \chi(\Thetab_t(\thetab)) \omega_0^{(m)} (d \thetab)\\
        + & \int_0^1 \int_D \nabla \chi\big(\Thetab_t (\thetab) + m^{-1/2} \eta \, \Tb_t^{(m)}(\thetab) \big) \cdot \Tb_t^{(m)}(\thetab)  \mu_0^{(m)}(d \thetab) d\eta~. 
    \end{split}
\end{equation}
As shown in Appendix~\ref{app-dclt1}, we can take the limit $m\to\infty$ of this formula to obtain: 
\begin{proposition}[Dynamical CLT - I]
\label{prop:dclt}
Under Assumptions~\ref{ass:unit_1} and \ref{ass:iid_init}, 
$\forall t \geq 0$,  as $m\to\infty$ we have $\omega_t^{(m)} \rightharpoonup \omega_t$ weakly in law
with respect to $\PP_0$, where $\omega_t$ is such that given a test function $\chi: D\to \RR$,
\begin{equation}
    \label{eq:omega_t}
    \int_D \chi(\thetab) \omega_t(d \thetab) =  \int_D \chi(\Thetab_t(\thetab)) \omega_0 (d \thetab)+ \int_D \nabla \chi(\Thetab_t(\thetab)) \cdot \Tb_t(\thetab) \mu_0(d \thetab)~.  
\end{equation}
Here $\omega_0$ is the Gaussian measure with mean zero and covariance
\begin{equation}
\label{eq:covomega0}
    \EE_{0} \left[\omega_0(d\thetab)\omega_0(d\thetab')\right] =
    \mu_0(d\thetab) \delta_{\thetab}(d\thetab') - \mu_0(d\thetab)\mu_0(d\thetab')~,
  \end{equation}
  where $\EE_0$ denotes expectation over $\PP_0$,
and $\Tb_t= \lim_{m\to\infty} m^{1/2}(\Thetab^{(m)}_t - \Thetab_t) $ is the flow solution to
\begin{equation}
\label{eq:T_t}
    \begin{aligned}
        \dot \Tb_t (\thetab) = & - \nabla \nabla V(\Thetab_t(\thetab),
    \mu_t)  \Tb_t(\thetab)
    - \int_D \nabla \nabla'
    K(\Thetab_t(\thetab), \Thetab_t(\thetab'))  \Tb_t(\thetab')
    \mu_0(d\thetab') \\
    &  - \int_D \nabla K(\Thetab_t(\thetab),\Thetab_t(\thetab'))
    \omega_0(d\thetab')
    \end{aligned}
\end{equation}
with initial condition $\Tb_0=0$ and where $\Thetab_t$ solves~\eqref{eq:flow} and $\nabla \nabla V(\Thetab_t(\thetab), \mu_t)$ is a shorthand for
\begin{equation}
    \label{eq:gradgradV}
    \nabla \nabla V(\Thetab_t(\thetab), \mu_t)
    = -\nabla \nabla F(\Thetab_t(\thetab)) + \int_D \nabla \nabla K(\Thetab_t(\thetab),\Thetab_t(\thetab')) \mu_0(d\thetab')~.
\end{equation}
\end{proposition}
\noindent A direct consequence of this proposition and
formula~\eqref{eq:g_vs_omega} is:
\begin{corollary}
  \label{prop:dcltg}
   Under Assumptions~\ref{ass:unit_1} and \ref{ass:iid_init}, 
$\forall t \geq 0$, as $m\to\infty$ we have 
$g_t^{(m)} \to g_t$ pointwise in law
with respect to $\PP_0$, where $g_t$ is given in terms of the limiting measure $\omega_t$ or the flow $\Tb_t$ as
\begin{equation}
\label{eq:g_t_omega_t0}
    g_t = \int_D \varphi(\thetab, \cdot) \omega_t(d \thetab) = \int_D \varphi(\Thetab_t(\thetab), \cdot) \omega_0 (d \thetab)  + \int_D \nabla \varphi(\Thetab_t(\thetab), \cdot) \cdot \Tb_t(\thetab) \mu_0(d \thetab) ~.
\end{equation}
  \end{corollary}
\noindent
It is interesting to comment on the origin of both terms at the right
hand side of~\eqref{eq:omega_t} and, consequently, \eqref{eq:g_t_omega_t0}.
The first term captures the deviations induced by fluctuations of
$\mu_0^{(m)}$ around $\mu_0$ assuming that the flow $\Thetab^{(m)}_t$
is unaffected by these fluctuations, and remains equal to
$\Thetab_t$. In particular, this term is the one we would obtain if we
were to resample $\mu_t^{(m)} $ from $\mu_t$ at every $t\ge0$,
i.e. use $\bar{\mu}_t^{(m)} = m^{-1} \sum_{i=1}^m
\delta_{\bar{\thetab}_t^i}$ with $\{ \bar{\thetab}_t^i
\}_{i=1}^m$ sampled i.i.d. from $\mu_t$, so that  $\Thetab^{(m)}_t$
is identical to $\Thetab_t$ in \eqref{eq:omegat_t_n}. In this case, the limiting discrepancy measure $\bar{\omega}_t$ would simply be given by

\begin{equation}
    \int_D \chi(\thetab) \bar{\omega}_t(d \thetab) = \int_D \chi(\Thetab_t(\thetab)) \omega_0(d \thetab)~,
\end{equation}
while the associated deviation in the represented function would read 
\begin{equation}
\label{eq:gtbar0}
\bar{g}_t = \int_D \varphi(\thetab, \cdot) \bar{\omega}_t(d \thetab) = \int_D \varphi(\Thetab_t(\thetab),\cdot) \omega_0(d\thetab)~.
\end{equation}
The second term at  right
hand side of~\eqref{eq:omega_t} and~\eqref{eq:g_t_omega_t0} captures
the deviations to the flow~$\Thetab_t$ in~\eqref{eq:flowexplicit}
induced by the perturbation of~$\mu_0$, i.e. how much
$\Thetab^{(m)}_t$ differs from $\Thetab_t$
in~\eqref{eq:omegat_t_n}. In the limit as $m\to\infty$, these
deviations are captured by the solution $\Tb_t$ to~\eqref{eq:T_t}, as is apparent from~\eqref{eq:omegat_t_n2}.

The difference between $g_t$ and $\bar g_t$ can also be quantified via the following  Volterra equation, which can be derived from Proposition~\ref{prop:dclt} and relates the evolution of $g_t$ to that of $\bar g_t$.
\begin{corollary}[Dynamical CLT - II]
\label{prop:dclt2}
Under Assumptions~\ref{ass:unit_1} and \ref{ass:iid_init}, $\forall t \geq 0$, pointwise on $\Omega$, we have $g_t^{(m)} \to g_t$ in law with respect to $\PP_0$ as $m\to\infty$, where $g_t$ solves the Volterra equation
\begin{equation}
    \label{eq:gtvolt}
    g_t(\xb) +\int_0^t \int_\Omega \Gamma_{t,s}(\xb,\xb') g_s(\xb') \hat\nu(d\xb') ds= \bar g_t(\xb)~.
\end{equation}
Here $\bar g_t$ is given in~\eqref{eq:gtbar0} and we defined
\begin{equation}
    \label{eq:Gammats}
    \Gamma_{t,s} (\xb,\xb') = \int_D \langle \nabla_{\thetab} \varphi(\Thetab_t(\thetab)), J_{t,s} (\thetab) \nabla_{\thetab} \varphi(\Thetab_s(\thetab))\rangle \mu_0(d\thetab)~,
\end{equation}
where $J_{t,s}$ is the solution to
\begin{equation}
  \label{eq:Jeq}
    \frac{d}{dt} J_{t,s}(\thetab) =  -\nabla \nabla V(\Thetab_t(\thetab),\mu_t)  J_{t,s}(\thetab), \qquad J_{s,s} (\thetab) = \text{Id}~.
\end{equation}
\end{corollary}
\noindent This corollary is proven in Appendix~\ref{app-dclt2}. In a nutshell, \eqref{eq:gtvolt} can be established using Duhamel's principle on~\eqref{eq:T_t} by considering all terms at the right hand side except the first as the source term (hence the role of $J_{t,s}$) and inserting the result in~\eqref{eq:g_t_omega_t0}.

\subsection{Long-Time Behavior of the Fluctuations}
\label{sec:g_infty}
Next, we study the long-time behavior of $g_t$ and, in particular, evaluate
\begin{equation}
    \label{eq:limlossg}
    \lim_{t\to\infty} \EE_0 \|g_t\|_{\hat \nu}^2=
    \lim_{t\to\infty} \lim_{m\to\infty} m \EE_{0}  \|f^{(m)}_t - f_t \|_{\hat \nu}^2.
\end{equation}
This limit  quantifies the asymptotic  approximation error of $f^{(m)}_t$ around its mean field limit $f_t$ after gradient flow, i.e. if we take $m\to\infty$ first, then $t\to\infty$ -- taking these limits in the opposite order is of interest too but is beyond the scope of the present paper. Our main result is to show that, under certain assumptions to be specified below, the limit in~\eqref{eq:limlossg} is not only finite but necessarily upper-bounded by $\lim_{t\to\infty} \EE_0 \|\bar g_t\|_{\hat\nu}^2$ with $\bar g_t$ given in~\eqref{eq:gtbar0}. That is, the approximation error at the end of training is always no higher than than that obtained by resampling the mean-field measure~$\mu_\infty$ defined in Proposition~\ref{th:ltflow}.

It is useful to start by considering an idealized case, namely when the initial conditions are sampled as in Assumption~\ref{ass:iid_init} with $\mu_0=\mu_\infty$. In that case, there is no evolution at mean field level, i.e. $\Thetab_t(\thetab)=\Thetab_\infty(\thetab)=\thetab$, $\mu_t=\mu_\infty$, and $f_t=f_\infty = \int_D \varphi_\infty (\thetab,\cdot) \mu_\infty(d\thetab)$, but the CLT fluctuations still evolve. In particular, it is easy to see that the Volterra equation in~\eqref{eq:gtvolt} for $g_t$ becomes
\begin{equation}
    \label{eq:gtvoltinf}
    g_t(\xb) +\int_0^t \int_\Omega \Gamma^\infty_{t-s}(\xb,\xb')
    g_s(\xb') \hat\nu(d\xb') ds= \bar g_\infty(\xb)~.
\end{equation}
Here $\Gamma^\infty_{t-s}(\xb,\xb')$ is the Volterra kernel obtained by solving~\eqref{eq:Jeq} with $\nabla \nabla V(\Thetab_t(\thetab),\mu_t)$ replaced by $\nabla \nabla V(\thetab,\mu_\infty)$ and inserting the result in~\eqref{eq:Gammats} with $\Thetab_t(\thetab) = \thetab$ and $\mu_0=\mu_\infty$,
\begin{equation}
    \label{eq:limGamm}
    \Gamma^\infty_{t - s} (\xb,\xb') = \int_D \langle \nabla_{\thetab} \varphi(\thetab,\xb), e^{-(t-s) \nabla \nabla V(\thetab,\mu_\infty)} \nabla_{\thetab} \varphi(\thetab,\xb')\rangle \mu_\infty(d\thetab)~,
 \end{equation}
 and $\bar g_\infty$ is the Gaussian field with variance
\begin{equation}
    \label{eq:varginfty}
    \EE_0 \|\bar g_\infty\|^2_{\hat \nu} 
    = \int_D \|\varphi(\thetab,\cdot)\|_{\hat \nu}^2 \mu_\infty(d\thetab)- \| f_\infty\|_{\hat \nu}^2~.
\end{equation}
From~\eqref{eq:Vlim} in Proposition~\ref{th:ltflow} we know that $\nabla \nabla V(\thetab,\mu_\infty)$ is positive semidefinite for $\mu_\infty$-almost all $\thetab$. As a result, we prove in \ref{app:long_time_reg_asymp} that the Volterra kernel~\eqref{eq:limGamm} viewed as an operator on functions defined on $\Omega \times[0,T]$ is positive semidefinite. Therefore, we have
\begin{equation}
    \label{eq:gtvoltinfnorm}
    \begin{aligned}
    \int_0^T \|g_t\|^2_{\hat\nu} dt  &\le \int_0^T \|g_t\|^2_{\hat\nu} dt +\int_0^T \int_0^t \int_{\Omega\times\Omega} g_t(\xb)\Gamma^\infty_{t-s}(\xb,\xb')
    g_s(\xb') \hat\nu(d\xb)\hat\nu(d\xb') ds dt\\
    & = \int_0^T \EE_{\hat\nu} (g_t \bar g_\infty) dt \le T^{1/2}\|\bar g_\infty\|_{\hat \nu} \left(\int_0^T \|g_t\|^2_{\hat \nu} dt \right)^{1/2}~.
    \end{aligned}
\end{equation}
Together with~\eqref{eq:varginfty}, this implies that
\begin{theorem}
\label{prop:ginfty}
Under Assumptions~\ref{ass:unit_1}, \ref{ass:iid_init}, \ref{ass:bounded_supp} and \ref{ass:limflow_local_min}, with $\mu_0=\mu_\infty$ and $\mu_\infty$ as specified in Proposition~\ref{th:ltflow}, we have
\begin{equation}
    \label{eq:mcbound_lem33}
    \lim_{T\to\infty} \frac1T \int_0^T \EE_0 \|g_t\|^2_{\hat{\nu}} dt  \le \int_D \|\varphi(\thetab,\cdot)\|_{\hat \nu}^2 \mu_\infty(d\thetab)- \| f_\infty\|_{\hat \nu}^2~.
\end{equation}
\end{theorem}
\noindent 
This theorem indicates that, if we knew $\mu_\infty$ and could sample initial conditions for the parameters from it, it would still be favorable to train these parameters  as this would reduce the approximation error. Of course, in practice we have no \textit{a~priori} access to  $\mu_\infty$, and so the relevant question is whether \eqref{eq:mcbound_lem33} also holds if we sample initial conditions from any $\mu_0$ such that Proposition~\ref{th:ltflow} holds.

In light of \eqref{eq:g_t_omega_t0}, one way to address this question is to study the long-time behavior of $\Tb_t$. In the setup without regularization ($\lambda = 0$), we can do so by leveraging existing results that, under certain assumptions, the mean-field gradient flow converges to a global minimizer which interpolates the training data points exactly \cite{rotskoff2018neural, chizat2018global, mei2018mean, sirignano2020mean}. In this case, the following theorem shows that we can actually obtain stronger controls on the fluctuations than \eqref{eq:mcbound_lem33}, which we prove in Appendix~\ref{app:long_time_unreg}.

\begin{theorem}[Long-time fluctuations in the unregularized case]
  \label{thm:long_time_unreg}
Consider the ERM setting with $\lambda=0$ and under Assumptions \ref{ass:unit_1}, \ref{ass:iid_init} and \ref{ass:bounded_supp}. Suppose that as $t \to \infty$, $\mu_t$ converges to a global minimizer $\mu_\infty$ that interpolates the data, i.e. the function $f_\infty = \int_D \varphi(\thetab,\cdot) \mu_\infty(d\theta)$ satisfies
\begin{equation}
\label{eq:interpolate}
     \forall \xb \in \supp \hat{\nu} \ : \  f_\infty(\xb) = f_*(\xb) ~,
\end{equation}
and, furthermore, the convergence satisfies
\begin{equation}
    \label{eq:assump_unreg_maintext_loss}
    \int_0^\infty t \left ( \mathcal{L}(\mu_t) \right )^{1/2} dt < \infty
\end{equation}
Then \eqref{eq:mcbound_lem33} holds. Additionally,
\begin{enumerate}
    \item if Assumption \ref{ass:shallownn} also holds, i.e., in the shallow neural network setting, we further have 
    \begin{equation}
    \label{eq:lim_fluc_0}
    \lim_{T\to\infty} \frac1T \int_0^T \EE_0 \|g_t\|_{\hat{\nu}}^2 dt  = 0~;
\end{equation}
\item if $\mu_0 = \mu_\infty$, then $\| g_t \|_{\hat{\nu}}$ decreases monotonically in $t$.
\end{enumerate}
\end{theorem}
\noindent
Hence, in the shallow neural networks setting and under these assumptions, the fluctuations will eventually vanish in the $O(m^{-1/2})$ scale of CLT.
Note that for \eqref{eq:assump_unreg_maintext_loss} to hold, it is sufficient that $\mathcal{L}(\mu_t)$ decays at an asymptotic rate of $O(t^{-\alpha})$ with $\alpha > 4$. 
For instance, \cite{chen2022on} proves that in an ERM setting where the size of the training dataset is no larger than the input dimension (i.e. $n \leq d$), the loss converges to zero at a linear rate, which will satisfy the condition \eqref{eq:assump_unreg_maintext_loss}. We leave the search for weaker sufficient conditions for future work.

When the limiting measure $\mu_\infty$ does not necessarily interpolate the training data, such as when regularization is added,  we can proceed with the analysis of the long-time behavior of $\Tb_t$ under the following assumption on the long-time behavior of the curvature:

\begin{theorem}[Long-time fluctuations under assumptions on the curvature]
\label{thm:min_eig_ddV}
Let $\Lambda_t(\thetab)$ denote the smallest eigenvalue of the tensor $\nabla \nabla V(\Thetab_t(\thetab), \mu_t)$ defined in~\eqref{eq:gradgradV}
and assume that for a constant $C$ (to be specified in
Appendix~\ref{app:curvature_2}) such that
  \begin{equation}
  \label{eq:curv}
      -\int_D \min\{\Lambda_t(\thetab), 0 \} \mu_0(d \thetab) = O(e^{-Ct}) \qquad \text{as\ \  $t\to\infty$.}
  \end{equation}
Then \eqref{eq:mcbound_lem33} holds.
\end{theorem}
\noindent
This theorem is proven in
Appendix~\ref{app:curvature_2}.
To intuitively  understand~\eqref{eq:curv}, note that we know from \eqref{eq:Vlim} in Proposition~\ref{th:ltflow} that $\Lambda_t(\thetab)\to0$ $\mu_0$-almost surely as $t\to\infty$. Condition~\eqref{eq:curv} can therefore be satisfied  by having $\Lambda_t(\thetab)$ converge to zero sufficiently fast in the regions of $D$ where it is negative, or having the measure of these regions with respect to $\mu_0$ converge to zero sufficiently fast, or both.

Alternatively, in the regularized ($\lambda > 0$) ERM setting, we can obtain the following result when the support of $\mu_\infty$ is atomic, as expected on general grounds \cite{zuho1948, fisher1975spline, bach2017breaking, boyer2019representer, de2020sparsity}:
\begin{theorem}[Long-time fluctuations in the regularized case]
  \label{thm:long_time_reg_maintext}
Consider the ERM setting under Assumptions \ref{ass:unit_1}, \ref{ass:iid_init} and \ref{ass:bounded_supp}. Suppose further that as $t \to \infty$, $\mu_t$ converges to  $\mu_\infty$ satisfying
\begin{flalign}
      &\exists \sigma > 0 \ \text{ s.t. } \ \forall \thetab \in \supp\mu_\infty \ : \  \nabla \nabla V(\thetab,\mu_\infty) \succ \sigma  \text{Id} ~, \textit{ and} \label{eq:pd} \\
     &\text{$\Thetab_t$ admits an asymptotic uniform convergence rate of $O(t^{-\alpha})$ with $\alpha > 3/2$}. \label{eq:unif_rate_char}
\end{flalign}
Then \eqref{eq:mcbound_lem33} holds with the ``$\lim$" replaced by ``$\limsup$" on its LHS.
\end{theorem}
\noindent Theorem~\ref{thm:long_time_reg_maintext} is proven in Appendix~\ref{app:long_time_reg} by analyzing directly the Volterra equation~\eqref{eq:gtvolt} and establishing that its solution coincides with that of~\eqref{eq:gtvoltinf} in the limit as $t\to\infty$, a property that we also expect to hold more generally than under the assumptions of Theorem~\ref{thm:long_time_reg_maintext}. In fact, we prove in Appendix~\ref{app:long_time_reg} that \eqref{eq:unif_rate_char} can be replaced by a weaker condition, \eqref{eq:assump_reg_maintext}. We also discuss the relation between Theorem~\ref{thm:long_time_reg_maintext} and the work of \cite{chizat2019sparse} in Appendix~\ref{app:chizat}.

\subsection{The Monte-Carlo Bound and Regularization}
\label{sec:var_norm}
The bound \eqref{eq:mcbound_lem33} on the long-time fluctuations motivates us to control the term $\int_D \|\varphi(\thetab,\cdot)\|_{\hat \nu}^2 \mu_\infty(d\thetab)$ using a suitable choice of regularization in \eqref{eq:ermmeasure}. In the following, we restrict our attention to the shallow neural networks setting, and further assume that
\begin{assumption}
\label{ass:D_hat_compact}
$\hat{D}$ is compact.
\end{assumption}
\noindent Under this assumption, there is
\begin{equation}
\int_D \| \varphi(\thetab, \cdot)\|_{\hat{\nu}}^2 \mu (d \thetab) = \int_D \int_\Omega |\varphi(\thetab, \xb)|^2 \hat\nu(d\xb)\mu(d\thetab)  \leq  \hat{K}_M \int_D c^2 \mu(d\thetab)~,    \end{equation}
where $\hat K_M = \max_{\zb\in\hat D} \|\hat \varphi(\zb,\cdot)\|_{\hat\nu}^2$. Thus, we consider regularization with $r(\thetab) = \frac{1}{2} c^2$, in which case \eqref{eq:ermmeasure} becomes
\begin{equation}
\label{eq:ermmeasure_2norm}
    \min_{\mu \in \mathcal{P}(D)} \mathcal{L}(\mu) \quad \text{with} \quad  \mathcal{L}(\mu)  := \tfrac{1}{2} \left\| f[\mu] - f_*\right\|_{\hat\nu} ^2 + \tfrac{1}{2} \lambda\int_D c^2 \mu(d \thetab)~.
\end{equation}
Interestingly, this choice of regularization leads to learning in the function space $\mathcal{F}_1$ \cite{bach2017breaking} (or alternatively, the Barron space \citep{ma2019barron}) associated with $\hat{\varphi}$, which is equipped with the \emph{variation norm} (or the \emph{Barron norm}) defined as
\begin{equation}
 \label{eq:bli}
    |\gamma_q( f)| := \inf_{\mu \in \mathcal{P}(D)} \left\{ \textstyle{\int_{D}}  |c|^q \mu(d\thetab);~ f(\xb) = \textstyle{\int_{D}} c \hat{\varphi}(\zb, \xb) \mu(d\thetab)\, \right\} = |\gamma_1( f)|^q ~,\qquad q\ge 1~.
\end{equation}
We call $\int_D |c|^q \mu(d \thetab)$ the \textit{$q$-norm} of $\mu$. One can verify \citep[Proposition 1]{ma2019barron} that indeed, using any $q \geq 1$ above yields the same norm because $\mu$, the object defining the integral representation (\ref{eq:integralrep}), is in fact a \emph{lifted} version of a more `fundamental' object $\gamma=\int_\RR c\mu(dc,\cdot) \in \mathcal{M}(\hat{D})$, the space of signed Radon measures over $\hat{D}$. They are related via the projection 
 \begin{equation}
     \int_{\hat{D}} \chi(\zb) \gamma(d\zb) = \int_{D} c  \chi(\zb) \mu(d\thetab)
 \end{equation}
for all continuous test functions $\chi: \hat{D} \to \R$. One can also verify  \cite{chizat2019sparse} that
$\gamma_1( f) = \inf \{ \| \gamma\|_{\mathrm{TV}};\, f(\xb) = \int_{\hat{D}} \hat{\varphi}(\zb, {\xb}) \gamma(d\zb) \}$, where $\| \gamma\|_{\mathrm{TV}}$ is the \emph{total variation} of $\gamma$ \cite{bach2017breaking}.

The space $\mathcal{F}_1$ contains any RKHS whose kernel is generated as an expectation over features  $k(\xb, \xb') = \int_{\hat D} \hat\varphi(\zb,\xb) \hat\varphi(\zb,\xb') \hat \mu_0(d\zb)$  with a base measure $\hat \mu_0 \in \mathcal{P}(\hat D)$, but it provides crucial approximation advantages over such RKHS at approximating certain non-smooth, high-dimensional functions with hidden low-dimensional structure, giving rise to powerful generalization guarantees \cite{bach2017breaking}. This also motivates the study of overparametrized shallow networks with the scaling as in (\ref{eq:basic}), as opposed to the NTK scaling of $m^{-1/2}$ \cite{jacot2018neural}.

To learn in $\mathcal{F}_1$, a canonical approach is to consider the ERM problem
\begin{equation}
\label{eq:kiki}
    \min_{f \in \mathcal{F}_1}  \tfrac{1}{2}\| f - f_*\|_{\hat \nu}^2 + \tfrac{1}{2} \lambda\gamma_1( f),
\end{equation}
By \eqref{eq:bli}, this is indeed equivalent to \eqref{eq:ermmeasure_2norm}.
In Appendix~\ref{app.prop_min}, we prove the following proposition, which shows that the measure obtained from \eqref{eq:ermmeasure_2norm} indeed has its $2$-norm controlled:
\begin{proposition}
  \label{th:min}
  Under Assumptions \ref{ass:shallownn}, \ref{ass:unit_1}, and \ref{ass:D_hat_compact}, $\mathcal{L}$ has no local minima and its global minimum  value can only be attained at measures $\mu_\lambda \in \mathcal{P}(D)$ such that both $f_\lambda = \int_D \varphi(\thetab,\cdot)\mu_\lambda(d\thetab)$ and $c_\lambda = \int_D |c| \mu_\lambda(d\thetab) = \left ( \int_D |c|^2 \mu_{\lambda}(d \thetab) \right )^{1/2} \le \gamma_1( f_*)$ are unique, and
  \begin{equation}
    \label{eq:4}
    \lambda^2 |c_\lambda|^{2} \hat K_M^{-1} \le 
     \|f_\lambda-f_*\|_{\hat\nu}^2, \qquad 
    \|f_\lambda-f_*\|_{\hat\nu}^2 + \lambda |c_\lambda|^2
    \le \lambda |\gamma_1( f_*)|^2.
  \end{equation}
  where $\hat K_M = \max_{\zb\in\hat D} \|\hat \varphi(\zb,\cdot)\|_{\hat\nu}^2$.
\end{proposition}

\section{Numerical Experiments}
\label{sec.experiments}
\subsection{Student-Teacher Setting}
\label{sec.planted}
We first perform numerical experiments in a  student-teacher setting, using a shallow teacher network as the target function to be learned by shallow student networks with different widths~$m$ of the hidden layer. Both $\hat{D}$ and $\Omega$ are taken to be the unit sphere of $d=16$ dimensions, and we take $\hat{\varphi}(\zb, \xb) = \max(0, \langle \zb, \xb \rangle)$. The teacher network has two neurons, $(c_1, \zb_1)$ and $(c_2, \zb_2)$, in the hidden layer, with $c_1 = c_2 = 1$ and $\zb_1$ and $\zb_2$ sampled i.i.d. from the uniform distribution on $\hat{D}$ and then fixed across the experiments. We vary the width of the student network in the range of $m = 128, 256, 512, 1024$ and $2048$, with their initial $\zb_i$'s sampled i.i.d. from the uniform distribution on $\hat{D}$. We consider two ways for initializing the $c_i$'s of the student networks: 1) \emph{Gaussian-initialization}, where the $c_i$'s are sampled i.i.d. from $\mathcal{N}(0, 1)$; and 2) \emph{zero-initialization}, where each $c_i$ is set to be $0$. 

We train the student networks in two ways: using the \emph{population} loss or the \emph{empirical} loss. For the former scenario, the data distribution $\nu$ is chosen to be uniform on $\Omega$, which allows an analytical formula for the loss as well as its gradient.
The student networks are trained by gradient descent under $L_2$ loss. Moreover, we rescale both the squared loss and the gradient by $d$ in order to adjust to the $\frac{1}{d}$ factor resulting from spherical integrals, and set the learning rate (which is the step size for discretizing \eqref{eq:gd_particle}) to be $1$. The models are trained for $20000$ epochs.
For each choice of $m$, we run the experiment $\kappa=20$ times with different random initializations of the student network. The average fluctuation of the population loss is defined as $\frac{1}{\kappa} \sum_{k=1}^\kappa \|f^{(m)}_k - \bar f^{(m)}\|_{\nu}^2$ for the population loss, with $\bar f^{(m)} = \frac{1}{\kappa} \sum_{k=1}^\kappa f^{(m)}_k$ being the averaged model, similar to the approach in \cite{geiger2019disentangling}. The other plotted quantities -- loss, TV-norm and $2$-norm -- are averaged across the $\kappa$ number of runs. The TV-norm (i.e., $1$-norm) and $2$-norm are defined as in Appendix~\ref{sec:var_norm}. 

\begin{figure}[!htb]
    \centering
    \includegraphics[scale=0.41,trim=15pt 0 0 0, clip]{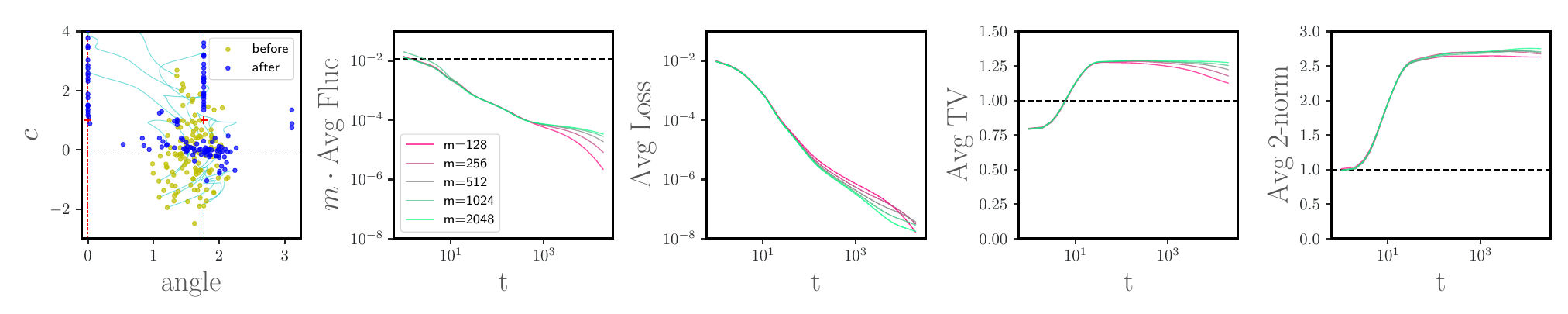} \\
    \vspace{-14pt}
    \includegraphics[scale=0.41,trim=15pt 0 0 0, clip]{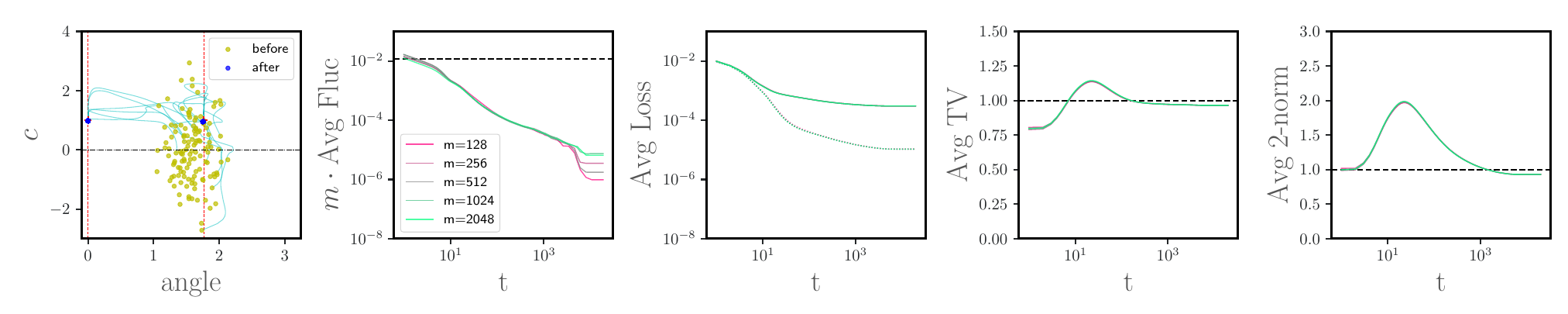} \\
    \vspace{-14pt}
    \includegraphics[scale=0.41,,trim=15pt 0 0 0, clip]{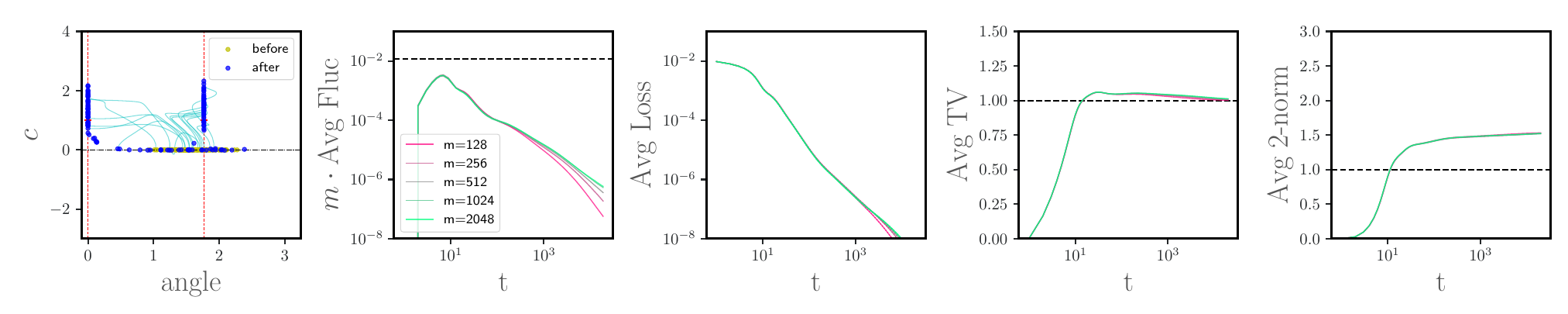}
    \caption{\small Results of the experiments in the student-teacher setting and where the student networks are trained by gradient descent on the \emph{population} loss. Each row corresponds to one setup. \emph{Row 1}: Using unregularized loss and non-zero-initialization; \emph{Row 2}: Using regularized loss with $\lambda=0.01$ and non-zero-initialization; \emph{Row 3}: Using unregularized loss and zero-initialization. In each row, \emph{Column 1} plots the trajectory of the neurons, $\thetab_i = (c_i, \zb_i)$, of a student network of width 128 during its training, with $x$-coordinate being the angle between $\zb_i$ and that of a chosen teacher's neuron and $y$-coordinate being $c_i$. The yellow dots, blue dots and cyan curves mark their initial values, terminal values, and trajectory during training. \emph{Columns 2-5} plot the average fluctuations (scaled by $m$), average loss, average TV norm, and average $2$-norm during training, respectively, computed across $\kappa=20$ runs with different random initializations of the student network for each choice of $m$. In \emph{Column 2}, the \emph{solid} curves give the average fluctuation of the population loss and the black horizontal \emph{dashed} line gives an approximate value of the asymptotic Monte-Carlo bound in \eqref{eq:mcbound_lem33} for this setting computed in Appendix~\ref{app.mc_value}.
    In \emph{Column 3}, the \emph{solid} curves indicate the total population loss, and the \emph{dotted} curves indicate the unregularized population loss (for the regularized case only).
    In \emph{Columns 4} and \emph{5}, the horizontal \emph{dashed} line gives the relevant norm of the teacher network.}
    \label{fig:results_pop}
\end{figure}

The results for the scenario of training under the \emph{population} loss are presented in Figures~ \ref{fig:results_pop}. As seen from \emph{Column 3} the average loss values remain similar over time for different choices of $m$, justifying the approximation by a mean-field dynamics.
In the unregularized case with non-zero initialization, the fluctuation of the population loss (shown in \emph{Column 2}) remains close to a $1/m$ scaling in roughly the first $10^3$ epochs, after which it decays faster for smaller $m$. Interestingly, this coincides with the tendency for the student neurons with $\zb$ not aligned with the teacher neurons to slowly have their $|c|$ decrease to zero due to a finite-$m$ effect, which is also reflected in the decrease in TV-norm. 
Aside from this phenomenon, the fluctuations decay at similar rates for different choices of $m$, which is consistent with our theory, since their dynamics are governed by the same dynamical CLT. 
Also, when regularization is added, each student neuron becomes aligned with one of the teacher neurons in both $\zb$ and $c$ after training; without regularization but using zero-initialization, after training, each student neuron either becomes aligned with one of the teacher neurons in $\zb$ or has $c$ close to zero. Both of these choices result in lower TV-norms and $2$-norms compared to using non-zero initialization and without regularization.

Next, we consider the \emph{empirical} loss scenario (ERM setting), using $n=32$ random vectors sampled i.i.d. from the uniform distribution $\nu$ on $\Omega$ as the training dataset, which then define the empirical data measure $\hat{\nu}(d \xb) = \frac{1}{n} \sum_{l=1}^n \delta_{\xb_l}(d \xb)$. We use the full training dataset for computing the gradient at every iteration. The other training setups are the same as when the population loss is used. We additionally plot the average fluctuation of the training loss, defined as $\frac{1}{\kappa} \sum_{k=1}^\kappa \|f^{(m)}_k - \bar f^{(m)}\|_{\hat{\nu}}^2$.

\begin{figure}[!htb]
    \centering
    \includegraphics[scale=0.41,trim=15pt 0 0 0, clip]{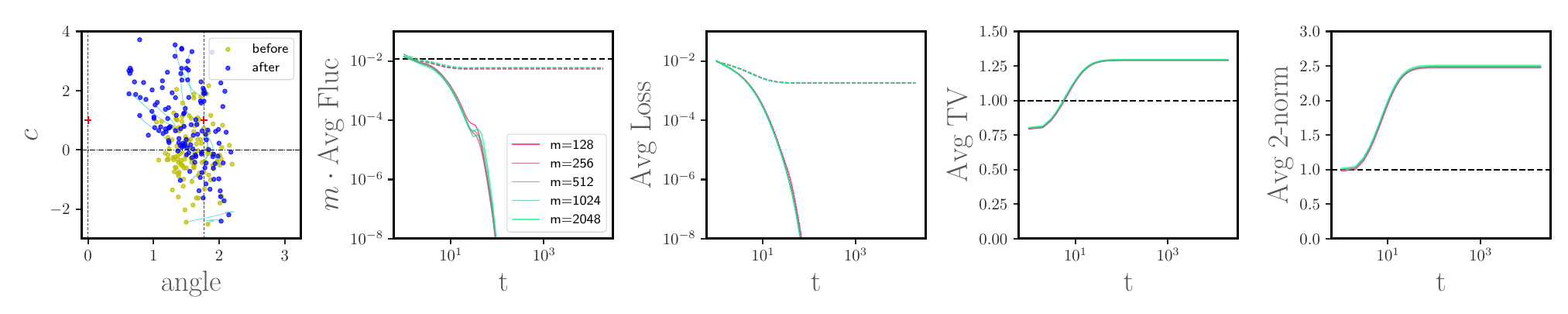} \\
    \vspace{-14pt}
    \includegraphics[scale=0.41,trim=15pt 0 0 0, clip]{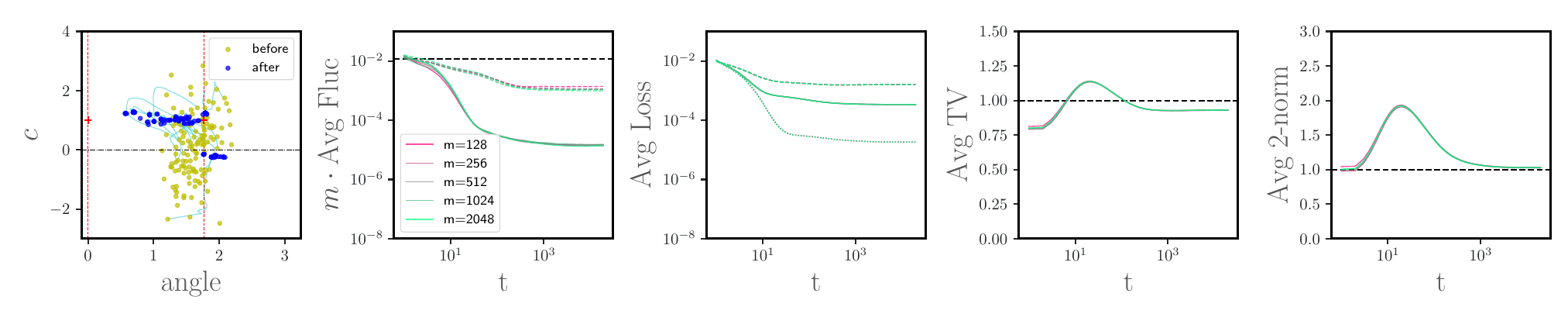} \\
    \vspace{-14pt}
    \includegraphics[scale=0.41,trim=15pt 0 0 0, clip]{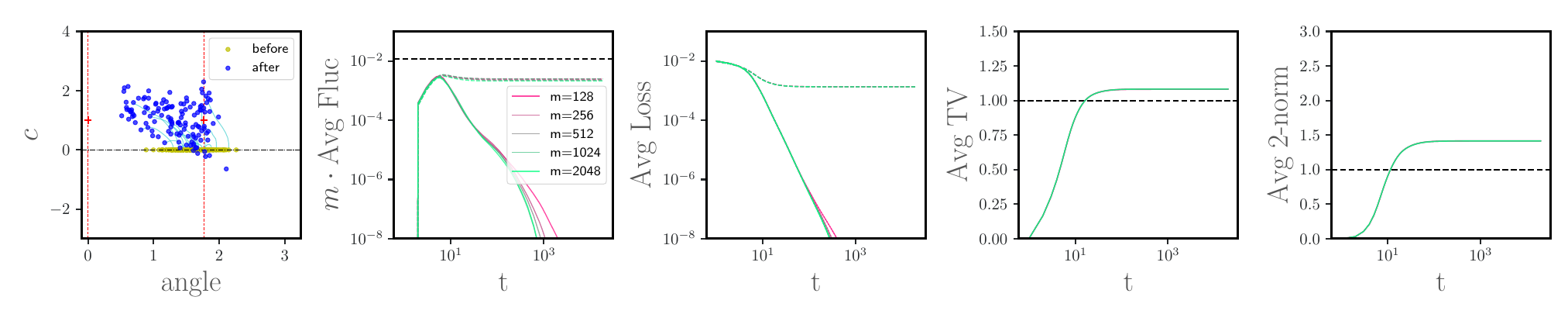}
     \vspace{-5pt}
    \caption{\small Results of the experiments in the student-teacher setting and where the student networks are trained by gradient descent on the \emph{empirical} loss.
    In \emph{Column 2}, the \emph{solid} curves indicate the average fluctuation in the \emph{training} loss, the \emph{dashed} curves indicate the average fluctuation in the \emph{population} loss computed analytically via spherical integrals, and the black horizontal \emph{dashed} line indicates an approximate value of the asymptotic Monte-Carlo bound in \eqref{eq:mcbound_lem33} for this setting computed in Appendix~\ref{app.mc_value}.
    In \emph{Column 3}, the \emph{solid} curves indicate the total \emph{training} loss, the \emph{dotted} curves the unregularized \emph{training} loss (for the regularized case only), and the \emph{dashed} curves the unregularized \emph{population} loss. All the other plot settings are identical to Figure~\ref{fig:results_pop}.}
    \label{fig:results_main}
\end{figure}

The results for the scenario of training under the \emph{empirical} loss is presented in Figures~\ref{fig:results_main}. Compared to the scenario of training under the \emph{population} loss, we wee that in the unregularized cases, both the average training loss and the average fluctuation of the training loss decay to below $10^{-8}$ within $10^3$ iterations, and the latter observation is consistent with \eqref{eq:lim_fluc_0}. In the regularized case, neither of them vanishes, but the average fluctuation of the training loss indeed remains below the asymptotic Monte-Carlo bound given in \eqref{eq:mcbound_lem33}, whose analytical expression and numerical value in this setup (under the approximation of replacing $\mu_\infty$, $f_\infty$ and $\hat{\nu}$ by the target measure, the target function and $\nu$, respectively) are given in Appendix \ref{app.mc_value}. Regularization and zero-initialization have a weaker effect in aligning the student neurons with the teacher neurons after training compared to the scenario of training under the \emph{population} loss, but they still result in lower TV-norm and $2$-norm, and moreover, lower average fluctuation and (slightly) lower average value of the \emph{population} loss. This demonstrates their positive effects on both approximation and generalization.

\subsection{Non-planted Case}
\label{sec:nonplanted}
We also conducted experiments in which the target function is not given by a teacher network but rather by  $f_*(\xb) = \int_{\hat{D}} \hat{\varphi}(\zb, \xb) \hat{\mu}_*(d \zb)$, where $\hat{\mu}_*$ is the uniform measure on the $1$-dimensional great circle in the first $2$ dimensions, i.e., $\{(\cos \theta, \sin \theta), 0, ..., 0: \theta \in [0, 2 \pi) \} \subseteq \mathbb{S}^d$, and where $\hat{D}$, $\Omega$, $\hat{\varphi}$ as well as the widths of the student networks remain the same as in the previous experiments.  The student networks are trained using gradient descent under the population loss where the data distribution $\nu$ is uniform on $\Omega$, which  
allows an analytical formula for the gradient using spherical integrals.

\begin{figure}[h]
    \centering
    \includegraphics[scale=0.45,,trim=13pt 28pt 0 0, clip]{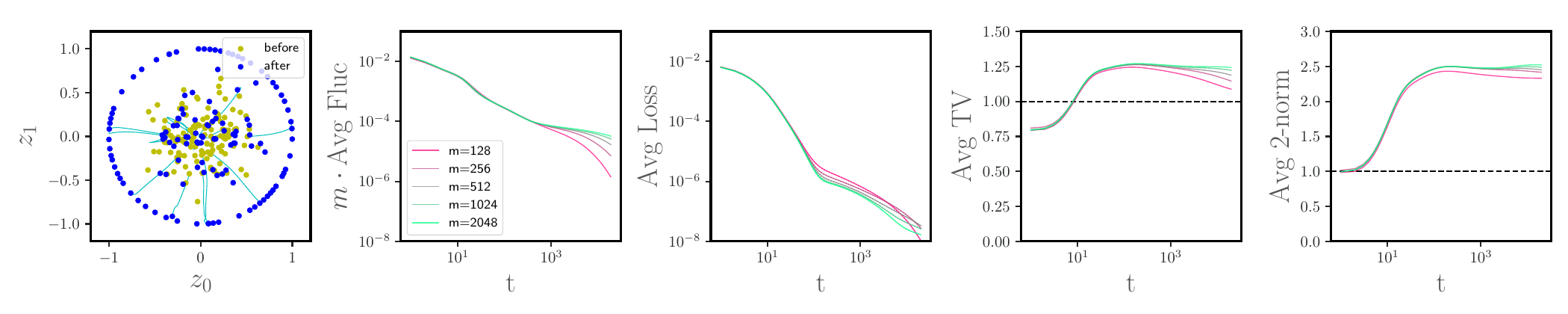}
    \vspace{-14pt}
    \includegraphics[scale=0.45,,trim=13pt 0 0 0, clip]{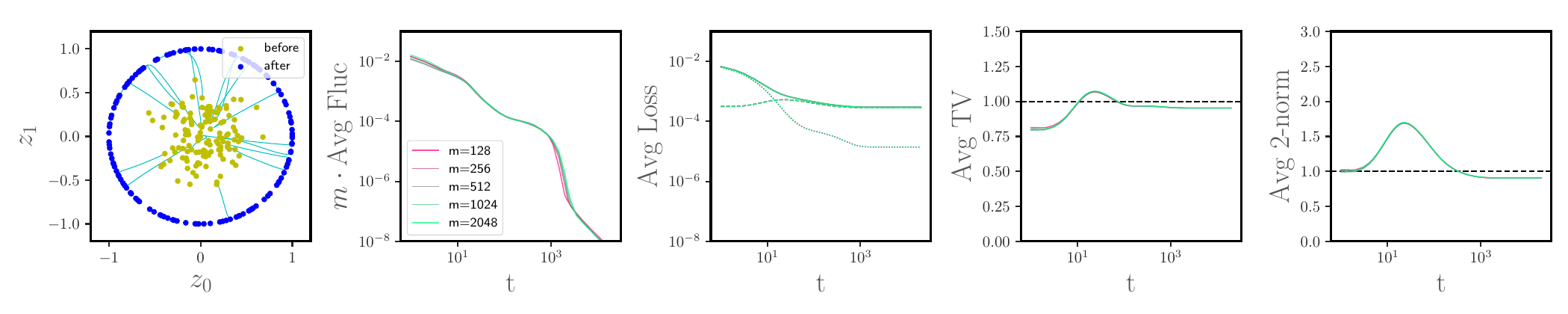}
    \vspace{-14pt}
    \includegraphics[scale=0.45,,trim=13pt 0 0 0, clip]{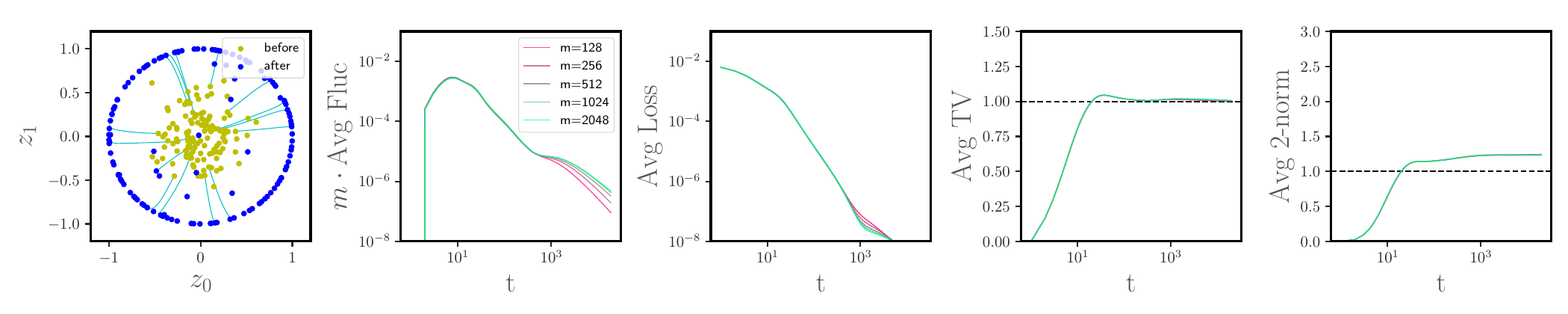}
    \caption{\small Results of the experiments with a non-planted target using the exact population loss, as described in Section~\ref{sec:nonplanted}. \emph{Row 1}: Using unregularized loss and non-zero-initialization; \emph{Row 2}: Using regularized loss with $\lambda=0.01$ and non-zero-initialization; \emph{Row 3}: Using unregularized loss and zero-initialization. In each row, \emph{Column 1} plots the projection in the first two dimensions of the $\zb_i$'s in the student network. The other columns show the same quantities as in Figure~\ref{fig:results_pop}.}
     \vspace{-5pt}
     \label{fig:results_nonp}
\end{figure}

The results are shown in Figure~\ref{fig:results_nonp}. We observe that the behaviors of the fluctuation are qualitatively similar to those found in Figure~\ref{fig:results_pop}.

\section{Conclusions}
Here we studied the deviations of shallow neural networks from their infinite-width limit, and how these deviations evolve during training by gradient flow. 
In the ERM setting, we established that under different sets of conditions, the long-term deviation under the Central Limit Theorem (CLT) scaling is controlled by a Monte Carlo (MC) resampling error,
giving width-asymptotic guarantees that do not depend on the data dimension explicitly. The MC resampling bound motivates a choice of regularization that is also connected to generalization via the variation-norm function spaces.

Our results thus seem to paint a favorable picture for high-dimensional learning, in which the optimization and generalization guarantees for the idealized mean-field limit could be transferred to their finite-width counterparts. However, we stress that these results are asymptotic, in that we take limits both in the width and time. In the face of negative results for the computational efficiency of training shallow networks \cite{manurangsi2018computational, livni2014computational,safran2018spurious,diakonikolas2020algorithms,goel2020superpolynomial}, an important challenge is to leverage additional structure in the problem (such as the empirical data distribution \cite{goldt2020gaussian}, or the structure of the minimizers \cite{de2020sparsity}) to provide nonasymptotic versions of our results, along the lines of \cite{chizat2019sparse} or \cite{pmlr-v125-li20a}. Finally, another clear direction for future research is to extend our techniques to deep neural architectures, in light of recent works that consider deep or residual models \cite{araujo2019mean, sirignano2019mean, nguyen2020rigorous,lu2020mean,wojtowytsch2020banach, fang2020modeling}. 

\section*{Acknowledgements}
This work benefited from  discussions with Lenaic Chizat and Carles Domingo-Enrich, and the authors sincerely thank Jiaheng Chen for pointing out an error in Theorem \ref{thm:long_time_unreg} in the previous version of this manuscript. Z.C. acknowledges support from the Henry MacCraken Fellowship. G.M.R. acknowledges support from the James S. McDonnell Foundation. J.B. acknowledges support from the Alfred P. Sloan Foundation, NSF RI-1816753, NSF CAREER CIF 1845360, and the Institute for Advanced Study. E. V.-E. acknowledges support from the National Science Foundation (NSF) Materials Research Science and Engineering Center Program Grant No. DMR-1420073, and from NSF Grant No. DMS-1522767. 

\bibliography{ref}
\bibliographystyle{plain}
\clearpage
\newpage
\appendix

{\hypersetup{linkcolor=black}
\addcontentsline{toc}{section}{Appendix} 
\part{Appendix} 
\parttoc
}

\section{Notations}
\label{app:notations}
We will use $\nabla \varphi(\thetab, \xb)$ and $\nabla \nabla \varphi(\thetab, \xb)$ to denote $\nabla_{\thetab} \varphi(\thetab, \xb)$ and $\nabla_{\thetab} \nabla_{\thetab} \varphi(\thetab, \xb)$, respectively. We will use $\nabla K(\thetab, \thetab')$ to denote $\nabla_{\thetab} K(\thetab, \thetab')$, $\nabla \nabla K(\thetab, \thetab')$ to denote $\nabla_{\thetab} \nabla_{\thetab} K(\thetab, \thetab')$, $\nabla' \nabla K(\thetab, \thetab')$ to denote $\nabla_{\thetab'} \nabla_{\thetab} K(\thetab, \thetab')$, and $\nabla' \nabla' K(\thetab, \thetab')$ to denote $\nabla_{\thetab'} \nabla_{\thetab'} K(\thetab, \thetab')$. We will write $V_t(\cdot)$ for $V(\cdot, \mu_t)$ and $V_\infty(\cdot)$ for $V(\cdot, \mu_\infty)$.\\

\noindent Let $D' = \cup_{t > 0} \supp \mu_t$. Under Assumption~\ref{ass:bounded_supp} and Proposition~\ref{th:ltflow}, $D'$ is bounded, and we denote its diameter by $|D'|$.
We will use $C_\varphi$, $C_{\nabla \varphi}$ and $C_{\nabla \nabla \varphi}$ to denote the supremum of $|\varphi(\thetab, \xb)|$, $|\nabla \varphi(\thetab, \xb)|$ and $|\nabla \nabla \varphi(\thetab, \xb)|$ over $\thetab \in D'$ and $\xb \in \supp \hat{\nu}$, which are all finite under Assumptions~\ref{ass:unit_1} and the boundedness of $D'$. We will use $L_{\nabla \nabla \varphi}$ to denote the (uniform-in-$\xb$) Lipschitz constant of $\nabla \nabla \varphi(\thetab, \xb)$ in $\thetab$, which is also finite under Assumption~\ref{ass:unit_1}.\\

\noindent The following notations will be used in Appendix~\ref{app:long_time_unreg}: 
Assuming that $D$ is Euclidean (under Assumption \ref{ass:unit_1}), let $\mathcal{V}(D)$ denote the space of random vector fields on $D$. It becomes a Hilbert space once equipped with the inner product 
\begin{equation}
    \begin{split}
    \big \langle \xib_1, \xib_2 \big \rangle_{0} :=
    \EE_{0} \int_D \xib_1(\thetab)\cdot \xib_2(\thetab) \mu_0(d \thetab),
\end{split}
\end{equation}
where $\xib_1$, $\xib_2$ denotes two random vector fields in $ \mathcal{V}(D)$. This inner product gives rise to the norm
\begin{equation}
    \| \xib \|_0^2 := \EE_{0} \int_D |\xib(\thetab)|^2 \mu_0(d \thetab)~.
\end{equation}
For each $t$, we define $\bb_t \in \mathcal{V}(D)$ as
\begin{equation}
\label{eq:b_t}
    \begin{split}
        \bb_t(\thetab) =& \int_D \nabla K (\Thetab_t(\thetab), \Thetab_t(\thetab')) \omega_0(d \thetab')
    \end{split}
\end{equation}
which depends on the random measure $\omega_0$. We define
two linear operators, $\mathcal{A}_t^{(K)}$ and $\mathcal{A}_t^{(V)}$ on $\mathcal{V}(D)$, as
\begin{align}
    (\mathcal{A}_t^{(K)} \xib)(\thetab) =& \int_D \nabla' \nabla K (\Thetab_t(\thetab), \Thetab_t(\thetab'))  \xib(\thetab') \mu_0(d \thetab') \\
    =& \int_\Omega \nabla \varphi(\Thetab_t(\thetab), \xb)  \Big ( \int_D \nabla \varphi(\Thetab_t(\thetab'), \xb)^\intercal \xib(\thetab') \mu_0(d \thetab') \Big ) \hat{\nu}(d \xb)~, \label{eq:A_t_K} \\
    (\mathcal{A}_t^{(V)} \xib)(\thetab) =&
    \nabla \nabla V(\Thetab_t(\thetab), \mu_t) \xib(\thetab)~, \label{eq:A_t_V}
\end{align}
for $\xib \in \mathcal{V}(D)$. Under Assumption \ref{ass:bounded_supp}, we also define $\bb_\infty$, $\mathcal{A}_\infty^{(K)}$, and $\mathcal{A}_\infty^{(V)}$ similarly by replacing $\Thetab_t(\cdot)$ with $\Thetab_\infty(\cdot)$.

Let $\mathcal{W}_n(\Omega)$ denote the space of random functions on $\Omega$. It becomes a Hilbert space once equipped with the inner product
\begin{equation}
    \langle \eta_1, \eta_2 \rangle_{\hat{\nu}, 0} := \EE_0 \int_\Omega \eta_1(\xb) \eta_2(\xb) \hat{\nu}(d \xb) = \frac{1}{n} \EE_0 \sum_{l=1}^n \eta_1(\xb_l) \eta_2(\xb_l)~,
\end{equation}
which gives rise to the norm
\begin{equation}
    \| \eta \|_{\hat{\nu}, 0}^2 :=  \langle \eta, \eta \rangle_{\hat{\nu}, 0} = \EE_0 \| \eta \|_{\hat{\nu}}^2~.
\end{equation}
With an abuse of notation, we will consider elements in $\mathcal{W}_n(\Omega)$ equivalently as random vectors on $\mathbb{R}^L$. Next, we can define $\mathcal{B}_t$ to be the operator that maps $\eta \in \mathcal{W}_n(\Omega)$ into the vector field
\begin{equation}
\label{eq:B_t}
    (\mathcal{B}_t \eta)(\thetab) = \int_\Omega \nabla \varphi(\Thetab_t(\thetab), \xb) \eta(\xb) \hat{\nu}(d \xb)
\end{equation}
in $\mathcal{V}(D)$.
Its transpose is
\begin{equation}
\label{eq:B_t_transpose}
    (\mathcal{B}_t^\intercal \xib)(\xb) = \int_D \nabla \varphi(\Thetab_t(\thetab), \xb) \xib(\thetab) \mu_0(d \thetab),
\end{equation}
which maps a vector field $\xib \in \mathcal{V}(D)$ back into $\mathcal{W}_n(\Omega)$.

\section{Long-Time Properties of the Mean-Field Gradient Flow}
\label{app_local_min}
\textit{Proof of Proposition~\ref{th:ltflow}:}
The compactness of $\cup_{t \geq 0} \supp {\mu_t}$ follows from \eqref{eq:limflow} and the compactness of $\supp {\mu_0}$ assumed in Assumption \ref{ass:iid_init}. $\mu_t \rightharpoonup \mu_\infty$ follows from \eqref{eq:weakn} and \eqref{eq:limflow}.

Under Assumption \ref{ass:bounded_supp}, $\Thetab_\infty$ is a local minimizer of the energy $\mathcal{E}$ defined in \eqref{eq:E0}. Consider a local perturbation $\epsilon \Thetab_{\Delta}$ to $\Thetab$. The energy value after the perturbation is
\begin{equation}
    \begin{split}
        \mathcal{E}(\Thetab_\infty + \epsilon \Thetab_{\Delta}) =& - \int_D F(\Thetab_\infty(\thetab) + \epsilon \Thetab_{\Delta}(\thetab)) \mu_0(d \thetab) \\
        & + \frac{1}{2} \int_D \int_D K(\Thetab_\infty(\thetab) + \epsilon \Thetab_{\Delta}(\thetab), \Thetab_\infty(\thetab') + \epsilon \Thetab_{\Delta}(\thetab')) \mu_0(d \thetab') \mu_0(d \thetab')~.
    \end{split}
\end{equation}
Under Assumptions \ref{ass:unit_1}, using Taylor expansion, we have
\begin{equation}
    \begin{split}
        F(\Thetab_\infty(\thetab) + \epsilon \Thetab_{\Delta}(\thetab)) =& F(\Thetab_\infty(\thetab)) + \epsilon \nabla F(\Thetab_\infty(\thetab)) \cdot \Thetab_{\Delta}(\thetab) \\
        & + \frac{1}{2} \epsilon^2 \langle \Thetab_{\Delta} (\thetab), \nabla \nabla F(\Thetab_\infty(\thetab)) \Thetab_{\Delta} (\thetab) \rangle + O(\epsilon^3)
    \end{split}
\end{equation}
\begin{equation}
    \begin{split}
        & K(\Thetab_\infty(\thetab) + \epsilon \Thetab_{\Delta}(\thetab), \Thetab_\infty(\thetab') + \epsilon \Thetab_{\Delta}(\thetab')) \\
        = & K(\Thetab_\infty(\thetab), \Thetab_\infty(\thetab')) + \epsilon \nabla K (\Thetab_\infty(\thetab), \Thetab_\infty(\thetab')) \Thetab_{\Delta}(\thetab) \\
        &+ \epsilon \nabla' K (\Thetab_\infty(\thetab), \Thetab_\infty(\thetab')) \Thetab_{\Delta}(\thetab') + \frac{1}{2} \epsilon^2 \langle \Thetab_{\Delta}(\thetab), \nabla \nabla K (\Thetab_\infty(\thetab), \Thetab_\infty(\thetab')) \Thetab_{\Delta}(\thetab) \rangle \\
        &+ \frac{1}{2} \epsilon^2 \langle \Thetab_{\Delta}(\thetab'), \nabla' \nabla' K (\Thetab_\infty(\thetab), \Thetab_\infty(\thetab')) \Thetab_{\Delta}(\thetab') \rangle \\
        & + \epsilon^2 \langle \Thetab_{\Delta}(\thetab), \nabla' \nabla K (\Thetab_\infty(\thetab), \Thetab_\infty(\thetab')) \Thetab_{\Delta}(\thetab') \rangle + O(\epsilon^3)~.
    \end{split}
\end{equation}
Hence, there is
\begin{equation}
    \begin{split}
        & \mathcal{E}(\Thetab_\infty + \epsilon \Thetab_{\Delta}) - \mathcal{E}(\Thetab_\infty) \\
        =& \epsilon \int_D \left ( - \nabla F(\Thetab_\infty(\thetab)) + \int_D \nabla K(\Thetab_\infty(\thetab), \Thetab_\infty(\thetab')) \mu_0(d \thetab') \right ) \Thetab_{\Delta}(\thetab) \mu_0(d \thetab) \\
        & + \frac{1}{2} \epsilon^2 \Bigg ( \int_D \langle \Thetab_{\Delta}(\thetab), \left ( \nabla \nabla F(\Thetab_\infty(\thetab)) + \int_D \nabla \nabla K(\Thetab_\infty(\thetab) , \Thetab_\infty(\thetab')) \mu_0(d \thetab') \right ) \Thetab_{\Delta}(\thetab) \rangle \mu_0(d \thetab) \\
        & \hspace{30pt} + \int_D \int_D \langle \Thetab_{\Delta}(\thetab), \nabla' \nabla K(\Thetab_\infty(\thetab) , \Thetab_\infty(\thetab')) \Thetab_{\Delta}(\thetab') \rangle \mu_0(d \thetab) \mu_0(d \thetab') \Bigg ) + O(\epsilon^3)~.
    \end{split}
\end{equation}
Since $\Thetab_{\Delta}$ is arbitrary can $\epsilon$ can be taken arbitrarily small, we see that for $\Thetab_\infty$ to be a local minimizer, the first-order condition is, $\forall \thetab \in \supp \mu_0$,
\begin{equation}
    - \nabla F(\Thetab_\infty(\thetab)) + \int_D \nabla K(\Thetab_\infty(\thetab), \Thetab_\infty(\thetab')) \mu_0(d \thetab') = 0~,
\end{equation}
or
\begin{equation}
    \nabla V(\Thetab_\infty(\thetab), \mu_\infty) = 0~,
\end{equation}
and the second-order condition is, $\forall \Thetab_{\Delta}$,
\begin{equation}
\begin{split}
    & \int_D \langle \Thetab_{\Delta}(\thetab), \left ( \nabla \nabla F(\Thetab_\infty(\thetab)) + \int_D \nabla \nabla K(\Thetab_\infty(\thetab) , \Thetab_\infty(\thetab')) \mu_0(d \thetab') \right ) \Thetab_{\Delta}(\thetab) \rangle \mu_0(d \thetab) \\
        & \hspace{30pt} + \int_D \int_D \langle \Thetab_{\Delta}(\thetab), \nabla' \nabla K(\Thetab_\infty(\thetab) , \Thetab_\infty(\thetab')) \Thetab_{\Delta}(\thetab') \rangle \mu_0(d \thetab) \mu_0(d \thetab') \geq 0~,
\end{split}
\end{equation}
or
\begin{equation}
\label{eq:second_order_cond}
\begin{split}
    & \int_D \langle \Thetab_{\Delta}(\thetab), \nabla \nabla V(\Thetab_\infty(\thetab), \mu_\infty) \Thetab_{\Delta}(\thetab) \rangle \mu_0(d \thetab) \\
    & +
    \int_D \int_D \langle \Thetab_{\Delta}(\thetab), \nabla' \nabla K(\Thetab_\infty(\thetab) , \Thetab_\infty(\thetab')) \Thetab_{\Delta}(\thetab') \rangle \mu_0(d \thetab) \mu_0(d \thetab') \geq 0~.
\end{split}
\end{equation}
Suppose for contradiction that $\exists D^- \subseteq D$ with $\mu_0(D^-) > 0$ such that $\nabla \nabla V(\Thetab_\infty(\thetab), \mu_\infty)$ is not positive semidefinite. Define $\Lambda_\infty(\thetab)$ to be the least eigenvalue of $\nabla \nabla V(\Thetab_\infty(\thetab), \mu_\infty)$. Then there is $\Lambda_\infty(\thetab) < 0$ on $D^-$. In addition, $\exists \zeta > 0$, $\exists D^-_0 \subseteq D^-$ with $\mu_0(D^-_0) > 0$ such that $\Lambda_\infty(\thetab) < - \zeta$. For $\thetab \in D^-_0$, let $\Thetab_{\Delta, 0}(\thetab)$ be a normalized eigenvector to $\nabla \nabla V(\Thetab_\infty(\thetab), \mu_\infty)$ associated with its least eigenvalue. Moreover, for $J \in \mathbb{N}^*$ that is large enough, we can select any subset $D^-_J \subset D^-_0$ such that $\mu_0(D^-_J) = \frac{1}{J} < \mu_0(D^-_0)$. Then, define
\begin{equation}
    \Thetab_{\Delta, J}(\thetab) = J^{1/2} \mathds{1}_{\thetab \in D^-_J} \Thetab_{\Delta, 0}(\thetab)~,
\end{equation}
Then, there is
\begin{equation}
\begin{split}
    & \int_D \int_D \langle \Thetab_{\Delta}(\thetab), \nabla' \nabla K(\Thetab_\infty(\thetab) , \Thetab_\infty(\thetab')) \Thetab_{\Delta}(\thetab') \rangle \mu_0(d \thetab) \mu_0(d \thetab')\\
    =& \int_\Omega \left | \int_D \nabla \varphi(\Thetab_\infty(\thetab), \xb) \Thetab_{\Delta, n} \mu_0(d \thetab) \right |^2 \hat{\nu}(d \xb) \\
    =& \int_\Omega \left | J^{1/2} \int_{D^-_J} \nabla \varphi(\Thetab_\infty(\thetab), \xb) \Thetab_{\Delta, 0} \mu_0(d \thetab) \right |^2 \hat{\nu}(d \xb) \\
    \leq & C_{\nabla \varphi}^2 J^{-1}~.
\end{split}
\end{equation}
On the other hand
\begin{equation}
    \begin{split}
        & \int_D \langle \Thetab_{\Delta, J}(\thetab), \nabla \nabla V(\Thetab_\infty(\thetab), \mu_\infty) \Thetab_{\Delta, J}(\thetab) \rangle \mu_0(d \thetab) \\
        =& \int_{D^-_J} J^{-1} \langle \Thetab_{\Delta, 0}(\thetab), \nabla \nabla V(\Thetab_\infty(\thetab), \mu_\infty) \Thetab_{\Delta, 0}(\thetab) \rangle \mu_0(d \thetab) \\
        \leq & - \zeta~.
    \end{split}
\end{equation}
Therefore, for $J$ large enough, we will have 
\begin{equation}
    \begin{split}
        & \int_D \int_D \langle \Thetab_{\Delta}(\thetab), \nabla' \nabla K(\Thetab_\infty(\thetab) , \Thetab_\infty(\thetab')) \Thetab_{\Delta}(\thetab') \rangle \mu_0(d \thetab) \mu_0(d \thetab')\\ 
        & + \int_D \langle \Thetab_{\Delta, n}(\thetab), \nabla \nabla V(\Thetab_\infty(\thetab), \mu_\infty) \Thetab_{\Delta, J}(\thetab) \rangle \mu_0(d \thetab) < 0~,
    \end{split}
\end{equation}
which contradicts \eqref{eq:second_order_cond}. Hence, we can conclude that $\mu_0$-almost surely, $\nabla \nabla V(\Thetab_\infty(\thetab), \mu_\infty)$ is positive semidefinite.

\section{Derivations of the Dynamical Central Limit Theorem}
\label{app_dclt}
\subsection{Proof of Proposition \ref{prop:dclt} (Dynamical CLT - I)}
\label{app-dclt1}
The following derivation is an adaptation of the approach in \cite{braun1977vlasov} for Vlasov interacting particle systems to our scenario.
To start, $\Thetab_t$ and $\Thetab_t^{(m)}$ are governed by the following equations, respectively:
\begin{equation}
\label{eq:Theta_dot}
  \begin{aligned}
  \dot \Thetab_t(\thetab) &= -\nabla V(\Thetab_t(\thetab),
  \mu_t),
  \qquad 
  &\Thetab_0(\thetab) &= \thetab \\
  \dot \Thetab_t^{(m)}(\thetab) &= -\nabla V(\Thetab_t^{(m)}(\thetab),
  \mu_t^{(m)}),
  \qquad 
  &\Thetab_0^{(m)}(\thetab) &= \thetab
  \end{aligned}
\end{equation}
Taking the difference between the two equations in \eqref{eq:Theta_dot} and using the mean value theorem, we get
\begin{equation}
    \begin{split}
        &\dot \Tb_t^{(m)}(\thetab)\\
        =& m^{1/2} \Big ( \dot \Thetab_t^{(m)}(\thetab) - \dot \Thetab_t(\thetab) \Big ) \\
        =& - m^{1/2} \Big ( \nabla V (\Thetab_t^{(m)}(\thetab), \mu_t^{(m)}) - \nabla V(\Thetab_t(\thetab), \mu_t) \Big ) \\
        =& -m^{1/2}\Big ( \nabla V (\Thetab_t^{(m)}(\thetab), \mu_t) - \nabla V(\Thetab_t(\thetab), \mu_t) \Big )
        - m^{1/2} \Big ( \nabla V (\Thetab_t, \mu_t^{(m)}) - \nabla V(\Thetab_t(\thetab), \mu_t) \Big ) \\
        &\quad - m^{1/2} \Big [ \Big ( \nabla V (\Thetab_t^{(m)}(\thetab), \mu_t^{(m)}) - \nabla V(\Thetab_t(\thetab), \mu_t^{(m)}) \Big ) 
        - \Big ( \nabla V (\Thetab_t^{(m)}, \mu_t) - \nabla V(\Thetab_t(\thetab), \mu_t) \Big ) \Big ] \\
        =& - \nabla \nabla V (\tilde{\Thetab}_{t, 1}^{(m)}(\thetab), \mu_t) \Tb_t^{(m)}(\thetab) 
        - \int_D \nabla K(\Thetab_t(\thetab), \thetab') \omega_t^{(m)} (d \thetab')\\
        &\quad- m^{-1/2} \Big ( \int_D \nabla \nabla K(\tilde{\Thetab}_{t, 2}^{(m)}(\thetab), \thetab') \omega_t^{(m)} (d \thetab') \Big )  \Tb_t^{(m)}(\thetab)~,
    \end{split}
\end{equation}
where  $\tilde{\Thetab}_{t. 1}^{(m)} (\thetab)$ and $\tilde{\Thetab}_{t, 2}^{(m)}(\thetab)$ denote points that lie on the line segment between $\Thetab_t(\thetab)$ and $\Thetab_t^{(m)}(\thetab)$.
Using \eqref{eq:omegat_t_n2}, we can substitute $\omega_t^{(m)}$ in the second term at the right hand side, for which we get
\begin{equation}
    \begin{split}
        \int_D \nabla K(\Thetab_t(\thetab), \thetab') \omega_t^{(m)} (d \thetab') 
        =& \int_D \nabla K (\Thetab_t(\thetab), \Thetab_t(\thetab')) \omega_0^{(m)} (d \thetab') \\
        &+ \int_D \nabla' \nabla K(\Thetab_t(\thetab), \tilde{\Thetab}_{t, 3}^{(m)}(\thetab'))  \Tb_t^{(m)}(\thetab') \mu_0(d \thetab') \\
        &+ m^{-1/2} \int_D \nabla' \nabla K (\Thetab_t(\thetab), \tilde{\Thetab}_{t, 3}^{(m)}(\thetab')) \Tb_t^{(m)} (\thetab') \omega_0^{(m)} (d \thetab')~.
    \end{split}
\end{equation}
Therefore, under Assumption \ref{ass:unit_1}, we have
\begin{equation}
\label{eq:dot_T_t_n}
\begin{split}
    \dot \Tb_t^{(m)}(\thetab) &= 
    - \nabla \nabla V (\tilde{\Thetab}_{t, 1}^{(m)}, \mu_t) \Tb_t^{(m)}(\thetab)\\
    &\quad
    - \int_D \nabla' \nabla K(\Thetab_{t}(\thetab), \tilde{\Thetab}_{t, 3}^{(m)}(\thetab'))  \Tb_t^{(m)}(\thetab') \mu_0(d \thetab') \\
    &\quad- \int_D \nabla K (\Thetab_t(\thetab), \Thetab_t(\thetab'))  \omega_0^{(m)} (d \thetab')
    + O(m^{-1/2}) ~.
\end{split}
\end{equation}
Now, we consider the limit as $m \to \infty$. By the standard CLT, we have that $\omega_0^{(m)} (d \thetab)\rightharpoonup \omega_0(d \thetab)$ weakly with respect to $\PP_0$, where $\omega_0(d \thetab)$ is the Gaussian measure with mean zero and covariance defined in~\eqref{eq:covomega0}. On the other hand, by finite-time LLN, we have $\Thetab_t^{(m)}(\thetab) \to \Thetab_t(\thetab)$ pointwise, $\mathbb{P}_0$-almost surely, and as a consequence $\tilde{\Thetab}_{t, 1}^{(m)}(\thetab), \bar{\Thetab}_{t, 3}^{(m)}(\thetab) \to \Thetab_t(\thetab)$ as well. Therefore, $\Tb_t^{(m)}(\thetab) \to \Tb_t(\thetab)$ pointwise, $\mathbb{P}_0$-almost surely, where the limiting $\Tb_t(\thetab)$ solves the equation obtained by taking the limit $m \to \infty$ on both sides of \eqref{eq:dot_T_t_n}, which becomes \eqref{eq:T_t}. \eqref{eq:T_t} should be solved with initial condition $\Tb_0(\thetab) =0$ since $\Tb_0^{(m)}(\thetab) = m^{1/2}( \Thetab_0^{(m)}(\thetab) - \Thetab_0(\thetab) ) =0$.

Finally, taking the limit $m \to \infty$ on both sides of the equation \eqref{eq:omegat_t_n2}, we deduce that $\omega_t^{(m)}(d\thetab)\rightharpoonup \omega_t(d\thetab)$ weakly, in law with respect to $\PP_0$, where the limiting $\omega_t(d\thetab)$ satisfies
\begin{equation}
  \int_D \chi(\thetab) \omega_t(d\thetab) = \int_D 
  \chi(\Thetab_t(\thetab)) \omega_0(d\thetab) +\int_D \nabla
  \chi(\Thetab_t(\thetab)) \cdot \Tb_t(\thetab)
  \mu_0(d\thetab) ~.
\end{equation}
This ends the proof of Proposition~\ref{prop:dclt}. \hfill $\square$

\subsection{Proof of Proposition \ref{prop:dclt2} (Dynamical CLT - II)}
\label{app-dclt2}

Recall from \eqref{eq:T_t} that
\begin{equation}
\begin{split}
    \dot{\Tb_t}(\thetab) &= -\nabla \nabla V(\Thetab_t(\thetab), \mu_t)  \Tb_t(\thetab) - \int_D \nabla' \nabla K(\Thetab_t(\thetab), \Thetab_t(\thetab')) \Tb_t(\thetab') \mu_0(d \thetab') \\
    & \quad- \int_D \nabla K(\Thetab_t(\thetab), \Thetab_t(\thetab')) \omega_0(d \thetab') \\
    &= -\nabla \nabla V(\Thetab_t(\thetab), \mu_t)  \Tb_t(\thetab) - \int_D \nabla K(\Thetab_t(\thetab), \thetab') \omega_t(d \thetab') ~.
\end{split}
\end{equation}
Since $\Tb_0(\thetab)=0$, we can  use Duhamel's principle to deduce that
\begin{equation}
\begin{split}
    \Tb_t(\thetab) &=  - \int_0^t  J_{t, s}(\thetab) \int_D \nabla K (\Thetab_s(\thetab), \thetab') \omega_s(d \thetab')  ds \\
    & = - \int_0^t \int_\Omega J_{t, s}(\thetab) \nabla \varphi(\Thetab_s(\thetab), \xb)  \int_D \varphi(\thetab', \xb) \omega_s(d \thetab') \hat \nu(d \xb) ds \\
    &= - \int_0^t \int_\Omega J_{t, s}(\thetab) \nabla \varphi(\Thetab_s(\thetab), \xb)  g_s(\xb) \hat \nu(d \xb) ds,
\end{split}
\end{equation}
where the tensor $J_{t, s}(\thetab)$ is the Jacobian defined in Proposition \ref{prop:dclt2}.
As a result
\begin{equation}
    \begin{split}
        g_t(\xb) &= \int_D \varphi(\thetab, \xb) \omega_t(d \thetab) \\
        & = \int_D \varphi(\Thetab_t(\thetab),\xb) \omega_0(d \thetab)+ \int_D \nabla \varphi(\Thetab_t(\thetab),\xb) \cdot \Tb_t(\thetab) \mu_0(d \thetab) \\
        &= \int_D \varphi(\Thetab_t(\thetab), \xb) \omega_0(d \thetab)\\
        & \quad  - \int_D \int_0^t \int_\Omega \langle\nabla \varphi(\Thetab_t(\thetab), \xb), J_{t, s}(\thetab) \nabla \varphi(\Thetab_s(\thetab), \xb')\rangle  g_s(\xb') \hat{\nu}(d \xb') ds \mu_0(d \thetab)\\
        &= \bar g_t(\xb)- \int_0^t \int_\Omega \int_D \langle\nabla \varphi(\Thetab_t(\thetab), \xb), J_{t, s}(\thetab) \nabla \varphi(\Thetab_s(\thetab), \xb') \rangle \mu_0(d \thetab) g_s(\xb') \hat{\nu}(d \xb') ds  \\
        &= \bar g_t(\xb) - \int_0^t \int_\Omega \Gamma_{t, s}(\xb, \xb') g_s(\xb') \hat{\nu}(d \xb') ds,
    \end{split}
\end{equation}
with $\bar g_t(\xb)$ and $\Gamma_{t, s}(\xb,\xb')$ defined in \eqref{eq:gtbar0} and \eqref{eq:Gammats}, respectively. This is \eqref{eq:gtvolt}. \hfill$\square$

\section{Long-Time Behavior of the Fluctuations}
\label{app.mcbound}

\subsection{Proof of Theorem~\ref{prop:ginfty} ($\mu_0 = \mu_\infty$ case)}
\label{app:long_time_reg_asymp}
With the argument outlined in Section \ref{sec:g_infty}, what remains to be shown is that $\Gamma_{t-s}^\infty$ is positive-semidefinite as a Volterra kernel, according to the definition in \cite{gripenberg_londen_staffans_1990}. We will utilize the following known result: 
\begin{proposition}[Gripenberg et al. \cite{gripenberg_londen_staffans_1990}]
\label{prop:gripenberg}
Let $k: [0, \infty) \to \mathbb{R}^{n \times n}$ be a convolution-type kernel for a linear Volterra equation in $\mathbb{R}^n$ . If $\forall \eta \in \mathbb{R}^n$, the function $t \mapsto \langle \eta, k(t) \eta \rangle$ is a nonnegative, nonincreasing and convex function on $(0, \infty)$, then $k$ is nonnegative, meaning that
$\forall \phi: [0, \infty) \to \mathbb{R}^{n}$ with compact support, there is
\begin{equation}
    \int_0^\infty \int_0^t \langle\phi(t), k(t-s) \phi(s)\rangle ds dt \geq 0~.
\end{equation}
\end{proposition}
\noindent
Thus, to take advantage of this proposition, we need to verify that $\forall \eta \in \mathbb{R}^n$, $\langle\eta, \Gamma^\infty_t \eta\rangle$ is
    
\textit{(1)} \textit{nonnegative:}
    \begin{equation}
    \begin{split}
        &\langle\eta, \Gamma^{\infty}_t \eta \rangle\\
        &= \int_{\Omega \times\Omega} \int_D \big \langle \nabla \varphi(\Thetab_\infty(\thetab), \xb), e^{- t \nabla \nabla V_\infty(\Thetab_\infty(\thetab))} \nabla \varphi(\Thetab_\infty(\thetab), \xb) \big \rangle \eta(\xb) \eta(\xb') \mu_0(d \thetab) \hat{\nu}(d \xb) \hat{\nu}(d \xb') \\
        &= \int_D \Big\langle\bb(\thetab), e^{- t \nabla \nabla V_\infty(\Thetab_\infty(\thetab))} \bb(\thetab)\Big \rangle \mu_0(d \thetab)
        \geq  0~,
    \end{split}
\end{equation}
where
\begin{equation}
    \bb(\thetab) = \int_\Omega \nabla \varphi(\Thetab_\infty(\thetab), \xb) \eta(\xb) \hat{\nu}(d\xb)
\end{equation}
because by assumption, $\forall \thetab \in D$, $\nabla \nabla V_\infty(\Thetab_\infty(\thetab))$ is positive semidefinite, and hence  
$e^{- t \nabla \nabla V_\infty(\Thetab_\infty(\thetab))}$ is a positive semidefinite operator;

\textit{(2)} \textit{nonincreasing:}
Taking derivative with respect to time,
\begin{equation}
    \begin{split}
        \frac{d}{dt} \langle\eta, \Gamma^{\infty}_t \eta\rangle &= - \int_D \Big \langle \bb(\thetab), \nabla \nabla V(\Thetab_\infty(\thetab)) e^{- t \nabla \nabla V_\infty(\Thetab_\infty(\thetab))} \bb(\thetab) \Big \rangle \mu_0(d \thetab) \leq 0,
    \end{split}
\end{equation}
because again, $\nabla \nabla V_\infty(\Thetab_\infty(\thetab))$ is positive semidefinite;

\textit{(3)} \textit{convex:}
Taking one more derivative with respect to time,
\begin{equation}
    \begin{split}
        \frac{d^2}{dt^2} \langle\eta, \Gamma^{\infty}_t \eta\rangle&=  \int_D \Big \langle \bb(\thetab), (\nabla \nabla V(\Thetab_\infty(\thetab)))^2  e^{- t \nabla \nabla V_\infty(\Thetab_\infty(\thetab))} \bb(\thetab)\Big \rangle \mu_0(d \thetab) \geq  0,
    \end{split}
\end{equation}
Therefore, we can apply Proposition \ref{prop:gripenberg}  to conclude that $\Gamma_{t-s}^\infty$ is PSD as a Volterra kernel, and so $\int_{t_0}^T \int_{t_0}^t \langle g_t, \Gamma^{\infty}_{t-s} g_s \rangle ds dt \geq 0$.

\subsection{Proof of Theorem \ref{thm:long_time_unreg} (Unregularized case)}
\label{app:long_time_unreg}
Recall that
\begin{equation}
\begin{split}
    \lim_{m\to\infty} m \EE_{0} \|f^{(m)}_t - f_t \|_{\hat{\nu}}^2 = \EE_{0} \|g_t\|_{\hat{\nu}}^2 
    &= \EE_{0} \int_{\Omega} \big | \int_D \varphi(\thetab, \xb) \omega_t(d \thetab) \big |^2 \hat{\nu}(d \xb) \\
    &= \EE_{0}  \int_{D \times D} K(\thetab, \thetab') \omega_t(d \thetab) \omega_t(d \thetab')~,
\end{split}
\end{equation}
where, with a slight abuse of notation, in this equation $\EE_0$ also denotes expectation over the randomness of the Gaussian distribution $\omega_0$ defined in Proposition~\ref{prop:dclt}.
From \eqref{eq:omega_t} in Proposition~\ref{prop:dclt},  this can be further expanded into
\begin{equation}
\begin{split}
  &\EE_{0}  \int_{D \times D} K(\thetab, \thetab') \omega_t(d \thetab) \omega_t(d \thetab')\\
  & = \EE_{0}  \int_{D \times D} \langle \Tb_t(\thetab), \nabla \nabla' K(\Thetab_t(\thetab), \Thetab_t(\thetab'))\Tb_t(\thetab') \rangle \mu_0(d\thetab)\mu_0(d\thetab') \\
  & + 2 \EE_{0}  \int_{D \times D} \nabla K(\Thetab_t(\thetab), \Thetab_t(\thetab'))\Tb_t(\thetab)  \mu_0(d\thetab)\omega_0(d\thetab')\\
  & + \EE_{0}  \int_{D \times D} K(\Thetab_t(\thetab), \Thetab_t(\thetab'))  \omega_0(d\thetab)\omega_0(d\thetab')~.
\end{split}
\end{equation}
The last term at the RHS is equal to $\EE_{0} \|\bar{g}_t\|_{\hat{\nu}}^2$ with $\bar{g}_t$ defined in \eqref{eq:gtbar0}. 
Using \eqref{eq:covomega0}, it can be explicitly computed as
\begin{equation}
\begin{split}
   \EE_{0} \|\bar{g}_t\|_{\hat{\nu}}^2 = & \EE_0  \int_{D \times D} K(\Thetab_t(\thetab), \Thetab_t(\thetab'))  \omega_0(d\thetab)\omega_0(d\thetab') \\
  &= \int_{D \times D} K(\Thetab_t(\thetab), \Thetab_t(\thetab'))  \left(\mu_0(d\thetab)\delta_{\thetab}(d\thetab') - \mu_0(d\thetab)\mu_0(d\thetab')\right)\\
  &= \int_{D } K(\thetab, \thetab)\mu_t(d\thetab) - \int_{D \times D} K(\thetab, \thetab)\mu_t(d\thetab)\mu_t(d\thetab') \\
  &= \int_{D } K(\thetab, \thetab)\mu_t(d\thetab) - \|f_t\|_{\hat{\nu}}^2~.
\end{split}
\end{equation}
Thus,
\begin{equation}
    \begin{split}
        \lim_{t\to\infty} \EE_{0} \|\bar{g}_t\|_{\hat{\nu}}^2 =& \lim_{t\to\infty}  \int_{D } K(\thetab, \thetab)\mu_t(d\thetab) - \|f_t\|_{\hat{\nu}}^2 \\
        =&  \int_{D } K(\thetab, \thetab)\mu_\infty(d\thetab) - \|f_\infty\|_{\hat{\nu}}^2 \\
        =& \EE_{0} \|\bar{g}_\infty\|_{\hat{\nu}}^2
    \end{split}
\end{equation}
and so
\begin{equation}
    \lim_{T\to\infty} \fint_0^T \EE_{0} \|\bar{g}_t\|_{\hat{\nu}}^2 dt = \EE_{0} \|\bar{g}_\infty\|_{\hat{\nu}}^2~,
\end{equation}
where here and below we denote $\fint_0^t [\cdot]  \ dt=\frac{1}{t} \int_0^t [\cdot] \ dt$. 
As a result, to prove \eqref{eq:mcbound_lem33} or \eqref{eq:lim_fluc_0} in Theorem~\ref{thm:long_time_unreg}, it suffices to establish that
\begin{equation}
\label{eq:diff_term_bound}
    \lim_{T \to \infty} \fint_0^T \mathfrak{D}_t dt \leq 0~,
\end{equation}
or
\begin{equation}
\label{eq:diff_term_sharp_bound}
    \lim_{T \to \infty} \fint_0^T \mathfrak{D}_t dt \leq - \EE_0 \| \bar{g}_\infty \|_{\hat{\nu}}^2~,
\end{equation}
respectively, where we defined 
\begin{equation}
\label{eq:diff_term}
\begin{split}
  \mathfrak{D}_t &:= \EE_{0}  \int_{D \times D} K(\thetab, \thetab') \omega_t(d \thetab) \omega_t(d \thetab')- \EE_0  \int_{D \times D} K(\Thetab_t(\thetab), \Thetab_t(\thetab'))  \omega_0(d\thetab)\omega_0(d\thetab')\\
  & = \EE_{0}  \int_{D \times D} \langle \Tb_t(\thetab), \nabla \nabla' K(\Thetab_t(\thetab), \Thetab_t(\thetab')) \Tb_t(\thetab') \rangle \mu_0(d\thetab)\mu_0(d\thetab') \\
  & + 2 \EE_{0}  \int_{D \times D} \nabla K(\Thetab_t(\thetab), \Thetab_t(\thetab'))\Tb_t(\thetab)  \mu_0(d\thetab)\omega_0(d\thetab')~.
\end{split}
\end{equation}
To this end, we examine \eqref{eq:T_t} as an infinite-dimensional ODE. 
With the Hilbert space $\mathcal{V}(D)$ defined in Appendix~\ref{app:notations} and $\bb_t$,
$\mathcal{A}_t^{(K)}$ and $\mathcal{A}_t^{(V)}$ defined by \eqref{eq:b_t}, \eqref{eq:A_t_K} and \eqref{eq:A_t_V}, respectively,
we can rewrite \eqref{eq:T_t} as the following ODE on $\mathcal{V}(D)$:
\begin{equation}
\label{eq:T_t_simple}
\begin{split}
    \dot{\Tb}_t = - (\mathcal{A}_t^{(K)} + \mathcal{A}_t^{(V)}) \Tb_t - \bb_t,
\end{split}
\end{equation}
We can also rewrite \eqref{eq:diff_term} as
\begin{equation}
    \mathfrak{D}_t = \langle \Tb_t, \mathcal{A}_t^{(K)} \Tb_t \rangle_0 + 2 \langle \Tb_t, \bb_t \rangle_0~.
\end{equation}
From \eqref{eq:T_t_simple}, we can deduce that
\begin{equation}
\label{eq:T_tmo}
    \begin{aligned}
        \frac{1}{2}\frac{d}{dt} \| \Tb_t \|_0^2 =  - \langle \Tb_t, \mathcal{A}_t^{(V)} \Tb_t \rangle_0 - \langle \Tb_t, \mathcal{A}_t^{(K)} \Tb_t \rangle_0 - \langle \Tb_t, \bb_t \rangle_0,
    \end{aligned}
\end{equation}
or equivalently
\begin{equation}
    \langle \Tb_t, \mathcal{A}_t^{(K)} \Tb_t \rangle_0 + \langle \Tb_t, \bb_t \rangle_0 = -\frac{1}{2} \frac{d}{dt} \| \Tb_t \|_0^2 - \langle \Tb_t, \mathcal{A}_t^{(V)} \Tb_t \rangle_0~.
\end{equation}
Therefore, we can rewrite
\eqref{eq:diff_term} as
\begin{equation}
    \begin{split}
        \mathfrak{D}_t =& 2 \left ( \langle \Tb_t, \mathcal{A}_t^{(K)} \Tb_t \rangle_0 + \langle \Tb_t, \bb_t \rangle_0 \right ) - \langle \Tb_t, \mathcal{A}_t^{(K)} \Tb_t \rangle_0 \\
        = & 2 \left ( -\frac{1}{2} \frac{d}{dt} \| \Tb_t \|_0^2 - \langle \Tb_t, \mathcal{A}_t^{(V)} \Tb_t \rangle_0 \right ) - \langle \Tb_t, \mathcal{A}_t^{(K)} \Tb_t \rangle_0 \\
        = & - \frac{d}{dt} \| \Tb_t \|_0^2 - 2 \langle \Tb_t, \mathcal{A}_t^{(V)} \Tb_t \rangle_0 - \langle \Tb_t, \mathcal{A}_t^{(K)} \Tb_t \rangle_0
    \end{split}
\end{equation}
and as a result, since $\Tb_0 = 0$,
\begin{equation}
\label{eq:diff_term_3_terms}
    \begin{split}
        \fint_0^T \mathfrak{D}_t dt = -\frac{1}{T}\| \Tb_T \|_0^2  - 2 \fint_0^T \langle \Tb_t, \mathcal{A}_t^{(V)} \Tb_t \rangle_0 dt - \fint_0^T \langle \Tb_t, \mathcal{A}_t^{(K)} \Tb_t \rangle_0 dt~.
    \end{split}
\end{equation}
Note that for all $t$, $\mathcal{A}_t^{(K)}$ is a positive semidefinite (PSD) operator on $\mathcal{V}(D)$, as $\forall \xib \in \mathcal{V}(D)$,
\begin{equation}
\begin{split}
    \langle \mathcal{A}_t^{(K)} \xib, \xib \rangle_{0} &= \EE_0 \int_{D\times D} \langle \xib(\thetab), \nabla \nabla' K(\Thetab_t(\thetab), \Thetab_t(\thetab'))  \xib(\thetab')\rangle \mu_0(d \thetab) \mu_0(d \thetab')  \\
    & = \EE_0 \int_\Omega \Big|\int_D \nabla \varphi(\Thetab_t(\thetab))\cdot  \xib(\thetab) \mu_0(d \thetab) \Big|^2 \hat{\nu}(d \xb) \geq  0~.
\end{split}
\end{equation}
This implies that $\fint_0^T \langle \Tb_t, \mathcal{A}_t^{(K)} \Tb_t \rangle_0 dt \geq 0$. Hence, to establish \eqref{eq:diff_term_bound}, it is sufficient to show that
\begin{equation}
\label{eq:T_t_A_t_V}
    \lim_{T \to \infty} \fint_0^T \langle \Tb_t, \mathcal{A}_t^{(V)} \Tb_t \rangle_0 dt = 0~.
\end{equation}
To this end, we need two lemmas that are proved below in Appendices~\ref{app.pf_lem_C1} and \ref{app.pf_lemma_C2}, respectively:
\begin{lemma}
\label{lem:integrable_unreg}
Assuming \eqref{eq:interpolate} and \eqref{eq:assump_unreg_maintext_loss} together with Assumptions~\ref{ass:unit_1}, \ref{ass:iid_init} and \ref{ass:bounded_supp}, we have
\begin{align}
    \int_0^\infty \| \mathcal{A}_t^{(V)} \|_{0} dt <& \infty \label{eq:A_V_integrable} \\
    \int_0^\infty \| \mathcal{A}_\infty^{(K)}- \mathcal{A}_t^{(K)} \|_{0} dt <& \infty \label{eq:A_K_diff_integrable} \\
    \int_0^\infty \| \bb_t - \bb_\infty \|_0 dt <& \infty \label{eq:b_diff_integrable}
\end{align}
\end{lemma}
\begin{lemma}
\label{lem:T_t_bounded_unreg}
Assuming \eqref{eq:interpolate} and \eqref{eq:assump_unreg_maintext_loss} together with Assumptions~\ref{ass:unit_1}, \ref{ass:iid_init} and \ref{ass:bounded_supp}, we have
\begin{equation}
\label{eq:T_t_bounded}
    \sup_{t < \infty} \|\Tb_t\|_0^2  < \infty~.
\end{equation}
\end{lemma}
\noindent
With these two lemmas, we can show that
\begin{equation}
\begin{split}
    \left | \int_0^T \langle \Tb_t, \mathcal{A}_t^{(V)} \Tb_t \rangle_0 dt \right | \leq & \int_0^T  \| \mathcal{A}_t^{(V)} \|_{0} \| \Tb_t \|_0^2 dt \\
    \leq & \left ( \int_0^T \| \mathcal{A}_t^{(V)} \|_0 dt \right ) \sup_{t < \infty} \|\Tb_t\|_0^2 \\
    < & \infty~,
\end{split}
\end{equation}
and therefore \eqref{eq:T_t_A_t_V} is satisfied. This finishes the proof of \eqref{eq:mcbound_lem33} under \eqref{eq:interpolate} and \eqref{eq:assump_unreg_maintext_loss} together with Assumptions~\ref{ass:unit_1}, \ref{ass:iid_init} and \ref{ass:bounded_supp}. 

Next, we show \eqref{eq:lim_fluc_0} under the additional condition of Assumption~\ref{ass:shallownn}. Thanks to \eqref{eq:diff_term_3_terms} and \eqref{eq:T_t_A_t_V}, it is sufficient to establish that
\begin{equation}
\label{eq:T_t_A_t_K}
    \lim_{T \to \infty} \fint_0^T \langle \Tb_t, \mathcal{A}_t^{(K)} \Tb_t \rangle_0 dt = \EE_0 \| \bar{g}_\infty \|_{\hat{\nu}}^2~.
\end{equation}
Heuristically, if $\Tb_\infty := \lim_{t \to \infty} \Tb_t$ exists, then from \eqref{eq:T_t_simple}, it has to satisfy
\begin{equation}
    - \bb_\infty = \left (\mathcal{A}_\infty^{(V)} + \mathcal{A}_\infty^{(K)} \right ) \Tb_\infty = \mathcal{A}_\infty^{(K)} \Tb_\infty~,
\end{equation}
as $\mathcal{A}_\infty^{(V)} = 0$ (because $\nabla \nabla V(\thetab, \mu_\infty) = \int_{\Omega} \varphi(\thetab, \xb) (f_\infty(\xb) - f_*(\xb)) \hat{\nu}(d \xb) = 0$ under the assumption of \eqref{eq:interpolate}). This equation implies that
\begin{equation}
    \left ( \Tb_\infty \right )^{||} = - \left ( \mathcal{A}_\infty^{(K)} \right )^\dag \bb_\infty~,
\end{equation}
where $\left ( \Tb_\infty \right )^{||}$ denotes the component of $\Tb_\infty$ in the range of $\mathcal{A}_\infty^{(K)}$, and $\left ( \mathcal{A}_\infty^{(K)} \right )^\dag$ denotes the Moore-Penrose pseudoinverse of $\mathcal{A}_\infty^{(K)}$. As a result,
\begin{equation}
\begin{split}
    \langle \Tb_\infty, \mathcal{A}_\infty^{(K)} \Tb_\infty \rangle_0 = & \langle \left (\Tb_\infty \right )^{||}, \mathcal{A}_\infty^{(K)} \left ( \Tb_\infty \right )^{||} \rangle_0 \\ = & \langle - \left ( \mathcal{A}_\infty^{(K)} \right )^\dag \bb_\infty, - \mathcal{A}_\infty^{(K)} \left ( \mathcal{A}_\infty^{(K)} \right )^\dag \bb_\infty \rangle_0 \\
    = & \langle \bb_\infty, \left ( \mathcal{A}_\infty^{(K)} \right )^\dag \bb_\infty \rangle_0~.
\end{split}
\end{equation}
Rigorously, without assuming the existence of $\Tb_\infty$, we can establish that
\begin{lemma}
\label{lem:g_infty_strict_bound_unreg}
Assuming \eqref{eq:interpolate} and \eqref{eq:assump_unreg_maintext_loss} together with Assumptions~\ref{ass:unit_1}, \ref{ass:iid_init} and \ref{ass:bounded_supp}, we have
\begin{equation}
\begin{split}
    \lim_{t \to \infty } \fint_0^t \langle \Tb_s, \mathcal{A}_s^{(K)} \Tb_s \rangle_0 ds \geq \langle \bb_\infty, \left ( \mathcal{A}_\infty^{(K)} \right )^{\dag} \bb_\infty \rangle_0~.
\end{split}
\end{equation}
As a consequence,
\begin{equation}
    \lim_{t \to \infty} \fint_0^t \EE_0 \|g_s\|_{\hat{\nu}}^2 dt \leq \EE_0 \|\bar{g}_\infty\|_{\hat{\nu}}^2 - \langle \bb_\infty, \left (  \mathcal{A}_\infty^{(K)} \right )^\dag \bb_\infty \rangle_0~.
\end{equation}
\end{lemma}
\noindent
This lemma is proved in \ref{app.pf_lemma_C3}. It implies that we only need to show that
\begin{equation}
    \langle \bb_\infty, \left (  \mathcal{A}_\infty^{(K)} \right )^\dag \bb_\infty \rangle_0 = \EE_0 \| \bar{g}_\infty \|_{\hat{\nu}}^2~. 
\end{equation}
This requires us to further exploit the relationship among $\mathcal{A}_\infty^{(K)}$, $\bb_\infty$ and $\bar{g}_\infty$. With the Hilbert space $\mathcal{W}_L(\Omega)$ defined in Appendix~\ref{app:notations} and $\mathcal{B}_t$ defined by \eqref{eq:B_t},
we can rewrite \eqref{eq:A_t_K} as
\begin{equation}
    \mathcal{A}_t^{(K)} = \mathcal{B}_t \mathcal{B}_t^\intercal~.
\end{equation}
Further, recall that
\begin{equation}
\label{eq:g_t_omega_t}
    g_t = \int_D \varphi(\thetab, \cdot) \omega_t(d \thetab) = \int_D \varphi(\Thetab_t(\thetab), \cdot) \omega_0 (d \thetab)  + \int_D \nabla \varphi(\Thetab_t(\thetab), \cdot) \cdot \Tb_t(\thetab) \mu_0(d \thetab) 
\end{equation}
\begin{equation}
\label{eq:gtbar}
\bar{g}_t = \int_D \varphi(\thetab, \cdot) \bar{\omega}_t(d \thetab) = \int_D \varphi(\Thetab_t(\thetab),\cdot) \omega_0(d\thetab)~.
\end{equation}
Therefore, we can write
\begin{equation}
\label{eq:g_g_bar}
    g_t = \bar{g}_t + \mathcal{B}_t^\intercal \Tb_t~,
\end{equation}
and
\begin{equation}
\label{eq:b_t_g_t}
    \bb_t = \mathcal{B}_t \bar{g}_t,
\end{equation}
Similar formulas hold when we replace $t$ by $\infty$. With these relations, we see that
\begin{equation}
    \begin{split}
        \left \langle \bb_\infty, \left ( \mathcal{A}_\infty^{(K)} \right )^{\dag} \bb_\infty \right \rangle_0 = & \left \langle \mathcal{B}_\infty \bar{g}_\infty, \left( \mathcal{B}_\infty \mathcal{B}_\infty^\intercal \right )^\dag \mathcal{B}_\infty \bar{g}_\infty \right \rangle_0 \\
        = & \EE_0 \| \left ( \mathcal{B}_\infty \right )^\dag \mathcal{B}_\infty \bar{g}_\infty \|_{\hat{\nu}}^2~,
    \end{split}
\end{equation}
because $\left( \mathcal{B}_\infty \mathcal{B}_\infty^\intercal \right )^\dag = (\mathcal{B}_\infty^\intercal)^\dag (\mathcal{B}_\infty)^\dag$.
Since $ \left ( \mathcal{B}_\infty \right )^\dag \mathcal{B}_\infty$ is the projection operator (matrix) onto the range of $\mathcal{B}_\infty^\intercal$ in $\mathbb{R}^n$, it is then sufficient to prove that
\begin{lemma}
\label{lem:g_bar_in_range}
Under Assumptions~\ref{ass:shallownn}, \ref{ass:unit_1}, \ref{ass:iid_init} and \ref{ass:bounded_supp}, $\mathbb{P}_0$-almost surely, $\bar{g}_\infty \in \mathrm{Ran}(\mathcal{B}_\infty^\intercal)$.
\end{lemma}
\noindent
Lemma \ref{lem:g_bar_in_range} is proven in Appendix~\ref{app.pf_lemma_C4} and it concludes the proof of \eqref{eq:lim_fluc_0} in Theorem \ref{thm:long_time_unreg}. \\

To show that $\| g_t \|_{\hat{\nu}}$ decreases monotonically when  $\mu_0 = \mu_\infty$, note that in this case $\mu_t = \mu_\infty$, $\forall t \geq 0$, and so $\mathcal{A}_t^{(V)} = \mathcal{A}_\infty^{(V)} = 0$, $\mathcal{A}_t^{(K)} = \mathcal{A}_\infty^{(K)}$ and $\bb_t = \bb_\infty$, $\forall t \geq 0$. Thus, \eqref{eq:T_t_simple} becomes
\begin{equation}
\label{eq:T_t_simple_asymp}
\begin{split}
    \dot{\Tb}_t = - \mathcal{A}_\infty^{(K)} \Tb_t - \bb_\infty,
\end{split}
\end{equation}
As will be shown in Lemma \ref{lem:in_range}, $\bb_\infty$ is in the range of $\mathcal{A}_\infty^{(K)}$. Therefore,
defining
\begin{equation}
    \ub_\infty = (\mathcal{A}^{(K)}_\infty)^\dag \bb_\infty~,
\end{equation}
and
\begin{equation}
    \zb_t = \Tb_t + \ub_\infty~,
\end{equation}
there is
\begin{equation}
    \dot{\zb}_t = - \mathcal{A}^{(K)}_\infty \zb_t~,
\end{equation}
whose solution can be written analytically as
\begin{equation}
    \zb_t = e^{- t \mathcal{A}^{(K)}_\infty} \zb_0 = e^{- t \mathcal{A}^{(K)}_\infty} \ub_\infty~.
\end{equation}
Thus,
\begin{equation}
    \Tb_t = \zb_t -\ub_\infty = -(I - e^{- t \mathcal{A}^{(K)}_\infty}) \ub_\infty
\end{equation}
Therefore, as $\bb_\infty = \mathcal{B}_\infty \bar{g}_\infty$, there is
\begin{equation}
\begin{split}
    {g}_t = & \bar{g}_\infty + \mathcal{B}_\infty^\intercal \Tb_t \\
    =& \bar{g}_\infty - \mathcal{B}_\infty^\intercal (I - e^{- t \mathcal{A}^{(K)}_\infty}) \ub_\infty \\
    =& \bar{g}_\infty - \mathcal{B}_\infty^\intercal (I - e^{- t \mathcal{A}^{(K)}_\infty}) (\mathcal{A}^{(K)}_\infty)^\dag \mathcal{B}_\infty \bar{g}_\infty~.
\end{split}
\end{equation}
Hence,
\begin{equation}
    \begin{split}
        |g_\infty|^2 =& |\bar{g}_\infty|^2 - 2 (*) + (**)~,
    \end{split}
\end{equation}
where
\begin{equation}
\begin{split}
    (*) =& \left( \mathcal{B}_\infty \bar{g}_\infty \right )^\intercal (I - e^{- t \mathcal{A}^{(K)}_\infty}) (\mathcal{A}^{(K)}_\infty)^\dag \mathcal{B}_\infty \bar{g}_\infty \\
    = & \bb_\infty^\intercal (I - e^{- t \mathcal{A}^{(K)}_\infty}) (\mathcal{A}^{(K)}_\infty)^\dag \bb_\infty
\end{split}
\end{equation}
and
\begin{equation}
\begin{split}
    (**) = & (\mathcal{B}_\infty \bar{g}_\infty)^\intercal (\mathcal{A}^{(K)}_\infty)^\dag (I - e^{- t \mathcal{A}^{(K)}_\infty})  \mathcal{B}_\infty \mathcal{B}_\infty^\intercal (I - e^{- t \mathcal{A}^{(K)}_\infty}) (\mathcal{A}^{(K)}_\infty)^\dag \mathcal{B}_\infty \bar{g}_\infty \\
    = & \bb_\infty^\intercal  (I - e^{- t \mathcal{A}^{(K)}_\infty})  \mathcal{B}_\infty \mathcal{B}_\infty^\intercal (I - e^{- t \mathcal{A}^{(K)}_\infty}) \bb_\infty~.
\end{split} 
\end{equation}

In the ERM setting, $\mathcal{A}_\infty^{(K)}$ is PSD with a finite number of nonzero eigenspaces. Consider a set of its orthonormal eigenfunctions that span those nonzero eigenspaces, $v_1, ..., v_k$, corresponding to eigenvalues $\lambda_1, ..., \lambda_k > 0$, respectively. 
As $\bb_\infty$ is in the range of $\mathcal{A}_\infty^{(K)}$ by Lemma \ref{lem:in_range}, we can decompose it as
\begin{equation}
    \bb_\infty = \sum_{i = 1}^k c_i v_i
\end{equation}
for some real numbers $c_i$'s. Thus, we can write
\begin{equation}
    \begin{split}
        (*) = & \left ( \sum_{i = 1}^k c_i v_i \right )^\intercal (I - e^{- t \mathcal{A}^{(K)}_\infty}) (\mathcal{A}^{(K)}_\infty)^\dag \left( \sum_{j = 1}^k c_j v_j \right ) \\
        = & \left ( \sum_{i = 1}^k c_i v_i \right )^\intercal \left( \sum_{j = 1}^k c_j \lambda_j^{-1}(1 - e^{-\lambda_j t})v_j \right ) \\
        =& \sum_{i=1}^k \lambda_j^{-1}(1 - e^{-\lambda_j t}) c_i^2~,
    \end{split}
\end{equation}
\begin{equation}
    \begin{split}
        (**) = & \left ( \sum_{i = 1}^k c_i v_i \right )^\intercal (\mathcal{A}^{(K)}_\infty)^\dag (I - e^{- t \mathcal{A}^{(K)}_\infty})  \mathcal{B}_\infty \mathcal{B}_\infty^\intercal (I - e^{- t \mathcal{A}^{(K)}_\infty}) (\mathcal{A}^{(K)}_\infty)^\dag \left ( \sum_{j = 1}^k c_j v_j \right ) \\
        = & \left ( \sum_{i = 1}^k c_i v_i \right )^\intercal (\mathcal{A}^{(K)}_\infty)^\dag (I - e^{- t \mathcal{A}^{(K)}_\infty})  \mathcal{A}^{(K)}_\infty (I - e^{- t \mathcal{A}^{(K)}_\infty}) (\mathcal{A}^{(K)}_\infty)^\dag \left ( \sum_{j = 1}^k c_j v_j \right ) \\
        = &  \left ( \sum_{i = 1}^k c_i v_i \right )^\intercal \left ( \sum_{j = 1}^k \lambda_j^{-1} \left ( 1 - e^{-\lambda_j t} \right )^2 c_j v_j \right ) \\
        =& \sum_{i=1}^k \lambda_j^{-1} \left (1 - e^{-\lambda_j t} \right )^2 c_i^2~.
    \end{split}
\end{equation}
Therefore, 
\begin{equation}
    \begin{split}
        |g_\infty|^2 = & |\bar{g}_\infty|^2 - 2 \sum_{i=1}^k \lambda_j^{-1} \left (1 - e^{-\lambda_j t} \right ) c_i^2 + \sum_{i=1}^k \lambda_j^{-1} \left (1 - e^{-\lambda_j t} \right )^2 c_i^2 \\
        = & |\bar{g}_\infty|^2 + \sum_{i=1}^k \lambda_j^{-1} \left (1 - e^{-\lambda_j t} \right ) \left (-1 - e^{-\lambda_j t} \right ) c_i^2 \\
        = & |\bar{g}_\infty|^2 - \sum_{i=1}^k \lambda_j^{-1} \left (1 - e^{-2 \lambda_j t} \right ) c_i^2~,
    \end{split}
\end{equation}
which is decreasing in time. This completes the proof of Theorem \ref{thm:long_time_unreg}.\hfill $\square$

\subsubsection{Proof of Lemma \ref{lem:integrable_unreg}}
\label{app.pf_lem_C1}
\textit{Proof of \eqref{eq:A_V_integrable}}: $\int_0^\infty \| \mathcal{A}_t^{(V)} \|_{0} dt < \infty$.\\
By the definition of the operator norm induced by $\| \cdot \|_0$ on $\mathcal{V}(D)$, $ \| \mathcal{A}_t^{(V)} \|_0$ is the smallest number $C_t$ such that $\forall \xib$,  there is 
\begin{equation}
    \| \mathcal{A}_t^{(V)} \|_0 = \sup_{\xib \in \mathcal{V}(D), \| \xib \|_0 \neq 0} \frac{\left | \langle \xib, \mathcal{A}_t^{(V)} \xib \rangle_0 \right |}{\| \xib \|^2_0} ~.
\end{equation}
In the unregularized case, a straightforward bound of $\left | \langle \xib, \mathcal{A}_t^{(V)} \xib \rangle_0 \right |$ is
\begin{equation}
\label{eq:asdf}
    \begin{split}
        \left | \langle \xib, \mathcal{A}_t^{(V)} \xib \rangle_0 \right |
        = & \left | \EE_0 \int_D \langle\xib(\thetab),  \nabla \nabla V(\Thetab_t(\thetab), \mu_t)  \xib (\thetab) \rangle\mu_0(d \thetab) \right | \\
        =& \left | \EE_0\int_D  \int_\Omega \langle\xib(\thetab),\nabla \nabla \varphi(\Thetab_t(\thetab), \xb)\xib(\thetab)\rangle (f_t(\xb) - f_*(\xb)) \hat{\nu}(d \xb) \ \mu_0(d \thetab) \right |  \\
        \leq & \EE_0 \int_D \int_\Omega C_{\nabla \nabla \varphi} \left | \xib(\thetab) \right |^2 \left | f_t(\xb) - f_*(\xb) \right | \hat{\nu}(d \xb) \mu_0(d \thetab) \\
        =& C_{\nabla \nabla \varphi} \| \xib \|_0^2 \int_\Omega \left | f_t(\xb) - f_*(\xb) \right | \hat{\nu}(d \xb) \\
        \leq & n^{1/2} C_{\nabla \nabla \varphi} \| \xib \|_0^2 \| f_t - f_* \|_{\hat{\nu}} \\
        =& n^{1/2} C_{\nabla \nabla \varphi} \| \xib \|_0^2 \left ( \mathcal{L}(\mu_t) \right )^{1/2}~.
    \end{split}
\end{equation}

Thus, we have
\begin{equation}
    \| \mathcal{A}_t^{(V)} \|_0 \leq n^{1/2} C_{\nabla \nabla \varphi} \| \xib \|_0^2 \left ( \mathcal{L}(\mu_t) \right )^{1/2}~.
\end{equation}
By the assumption \eqref{eq:assump_unreg_maintext_loss},
we thus have
\begin{equation}
    \int_0^\infty \| \mathcal{A}_t^{(V)} \|_{0} dt \leq 
    n^{1/2} C_{\nabla \nabla \varphi} \int_0^\infty \left ( \mathcal{L}(\mu_t) \right )^{1/2} dt < \infty
\end{equation}
which gives us the desired bound.\hfill $\square$\\\\
\textit{Proof of \eqref{eq:A_K_diff_integrable}}: $\int_0^\infty \| \mathcal{A}_\infty^{(K)}- \mathcal{A}_t^{(K)} \|_{0} dt < \infty$.\\
We have
\begin{equation}
    \begin{split}
        & \langle \xib, (\mathcal{A}_t^{(K)} - \mathcal{A}_\infty^{(K)}) \xib \rangle_0 \\
        =& \EE_0\int_\Omega \Big(\Big ( \int_D \nabla \varphi(\Thetab_t(\thetab), \xb)\cdot  \xib(\thetab) \mu_0(d \thetab) \Big )^2 - \Big ( \int_D \nabla \varphi(\Thetab_\infty(\thetab), \xb)\cdot  \xib(\thetab) \mu_0(d \thetab) \Big )^2 \Big)\hat \nu(d \xb) \\
        = & \EE_0 \int_\Omega \Big ( \int_D \big ( \nabla \varphi(\Thetab_t(\thetab), \xb) + \nabla \varphi(\Thetab_\infty(\thetab), \xb) \big )\cdot  \xib(\thetab) \mu_0(d \thetab) \Big )  \\
        & \times \Big ( \int_D \big ( \nabla \varphi(\Thetab_t(\thetab), \xb) - \nabla \varphi(\Thetab_\infty(\thetab), \xb) \big )\cdot \xib(\thetab) \mu_0(d \thetab) \Big ) \hat \nu(d \xb)~.
    \end{split}
\end{equation}
Hence, the absolute value of the expression above is upper-bounded by
\begin{equation}
    \begin{split}
        & \EE_0 \Big ( \int_D | \nabla \varphi(\Thetab_t(\thetab), \xb) + \nabla \varphi(\Thetab_\infty(\thetab), \xb) |  |\xib(\thetab)| \mu_0(d \thetab) \\
        & \quad\times  \int_D | \nabla \varphi(\Thetab_t(\thetab), \xb) - \nabla \varphi(\Thetab_\infty(\thetab), \xb) |  |\xib(\thetab)| \mu_0(d \thetab) \Big ) \\
        \leq & 2 C_{\nabla \varphi}  C_{\nabla \nabla \varphi}  \| \xib \|_0^2  \Big ( \int_D |\Thetab_t(\thetab) - \Thetab_\infty(\thetab)|^2 \mu_0(d \thetab) \Big )^{1/2}~.
    \end{split}
\end{equation}
Thus, by the assumption \eqref{eq:assump_unreg_maintext_loss}, we have
\begin{equation}
    \begin{split}
        \int_0^\infty \| \mathcal{A}_\infty^{(K)}- \mathcal{A}_t^{(K)} \|_{0} dt \leq & 2 C_{\nabla \varphi}  C_{\nabla \nabla \varphi} \int_0^\infty \left ( \int_D |\Thetab_t(\thetab) - \Thetab_\infty(\thetab)|^2 \mu_0(d \thetab) \right )^{1/2} dt \\
        \leq & 2 C_{\nabla \varphi}  C_{\nabla \nabla \varphi} \int_0^\infty \left ( \mathcal{L}(\mu_t) \right )^{1/2} dt \\
        < & \infty
    \end{split}
\end{equation} \hfill $\square$\\
\textit{Proof of \eqref{eq:b_diff_integrable}}: $\int_0^\infty \| \bb_t - \bb_\infty \|_0 dt < \infty$. \\
There is
\begin{equation}
\begin{split}
    &\bb_t(\thetab) - \bb_\infty(\thetab)\\ =& \int_D \big ( \nabla K(\Thetab_t(\thetab), \Thetab_t(\thetab')) - \nabla K(\Thetab_\infty(\thetab), \Thetab_\infty(\thetab')) \big ) \omega_0(d \thetab') \\
    =& \int_D \int_\Omega \nabla \varphi(\Thetab_t(\thetab), \xb) \cdot \nabla \varphi(\Thetab_t(\thetab'), \xb)^\intercal - \nabla \varphi(\Thetab_\infty(\thetab), \xb) \cdot \nabla \varphi(\Thetab_\infty(\thetab'), \xb)^\intercal \hat{\nu}(d \xb) \omega_0(d \thetab') \\
    =& \int_D \int_\Omega \nabla \varphi(\Thetab_t(\thetab), \xb) \cdot \nabla \varphi(\Thetab_t(\thetab'), \xb)^\intercal - \nabla \varphi(\Thetab_t(\thetab), \xb) \cdot \nabla \varphi(\Thetab_\infty(\thetab'), \xb)^\intercal \hat{\nu}(d \xb) \omega_0(d \thetab') \\
    & + \int_D \int_\Omega \nabla \varphi(\Thetab_t(\thetab), \xb) \cdot \nabla \varphi(\Thetab_\infty(\thetab'), \xb)^\intercal - \nabla \varphi(\Thetab_\infty(\thetab), \xb) \cdot \nabla \varphi(\Thetab_\infty(\thetab'), \xb)^\intercal \omega_0(d \thetab') \\
    =& \int_\Omega \nabla \varphi(\Thetab_t(\thetab), \xb) \cdot \left ( \int_D \left ( \nabla \varphi(\Thetab_t(\thetab'), \xb) - \nabla \varphi(\Thetab_\infty(\thetab'), \xb) \right ) \omega_0(d \thetab') \right )^\intercal \hat{\nu}(d \xb) \\
    & + \int_\Omega \left ( \nabla \varphi(\Thetab_t(\thetab), \xb) - \nabla \varphi(\Thetab_\infty(\thetab), \xb) \right) \left ( \int_D \nabla \varphi(\Thetab_\infty(\thetab'), \xb) \omega_0(d \thetab') \right )^\intercal \hat{\nu}(d \xb)~.
\end{split}
\end{equation}
Thus,
\begin{equation}
\begin{split}
    &\EE_0 \left | \bb_t(\thetab) - \bb_\infty(\thetab) \right |^2 \\
    \leq & \EE_0 \left | \int_\Omega \nabla \varphi(\Thetab_t(\thetab), \xb) \cdot \left ( \int_D \left ( \nabla \varphi(\Thetab_t(\thetab'), \xb) - \nabla \varphi(\Thetab_\infty(\thetab'), \xb) \right ) \omega_0(d \thetab') \right )^\intercal \hat{\nu}(d \xb) \right |^2 \\
    & + \EE_0 \left | \int_\Omega \left ( \nabla \varphi(\Thetab_t(\thetab), \xb) - \nabla \varphi(\Thetab_\infty(\thetab), \xb) \right) \left ( \int_D \nabla \varphi(\Thetab_\infty(\thetab'), \xb) \omega_0(d \thetab') \right )^\intercal \hat{\nu}(d \xb) \right |^2 \\
    \leq & \int_\Omega \left | \nabla \varphi(\Thetab_t(\thetab), \xb) \right |^2 \EE_0 \left | \int_D \left ( \nabla \varphi(\Thetab_t(\thetab'), \xb) - \nabla \varphi(\Thetab_\infty(\thetab'), \xb) \right ) \omega_0(d \thetab') \right |^2 \hat{\nu}(d \xb) \\
    & + \int_\Omega \left | \nabla \varphi(\Thetab_t(\thetab), \xb) - \nabla \varphi(\Thetab_\infty(\thetab), \xb) \right |^2 \EE_0 \left | \int_D \nabla \varphi(\Thetab_\infty(\thetab'), \xb) \omega_0(d \thetab') \right |^2 \hat{\nu}(d \xb) \\
    \leq & C_{\nabla \varphi}^2 \int_\Omega \EE_0 \left | \int_D \left ( \nabla \varphi(\Thetab_t(\thetab'), \xb) - \nabla \varphi(\Thetab_\infty(\thetab'), \xb) \right ) \omega_0(d \thetab') \right |^2 \hat{\nu}(d \xb) \\
    & + C_{\nabla \nabla \varphi}^2 \left | \Thetab_t(\thetab) - \Thetab_\infty(\thetab) \right |^2 \int_\Omega \EE_0 \left | \int_D \nabla \varphi(\Thetab_\infty(\thetab'), \xb) \omega_0(d \thetab') \right |^2 \hat{\nu}(d \xb)~.
\end{split}
\end{equation}
By the property of $\omega_0$, there is
\begin{equation}
    \begin{split}
        \EE_0 \left | \int_D \chi(\thetab) \omega_0(d \thetab) \right |^2 =& \int_D \left | \chi(\thetab) - \int_D \chi(\thetab') \mu_0(d \thetab') \right |^2 \mu_0(d \thetab) \\
        \leq & \int_D \left | \chi(\thetab) \right |^2 \mu_0(d \thetab)
    \end{split}
\end{equation}
for a test function $\chi$ on $D$. Thus,
\begin{equation}
    \begin{split}
        \EE_0 \left | \bb_t(\thetab) - \bb_\infty(\thetab) \right |^2 
    \leq & C_{\nabla \varphi}^2 \int_\Omega \int_D \left | \nabla \varphi(\Thetab_t(\thetab'), \xb) - \nabla \varphi(\Thetab_\infty(\thetab'), \xb) \right |^2 \mu_0(d \thetab') \\
    & + C_{\nabla \nabla \varphi}^2 \left | \Thetab_t(\thetab) - \Thetab_\infty(\thetab) \right |^2 \int_\Omega \int_D  \left | \nabla \varphi(\Thetab_\infty(\thetab'), \xb) \right |^2 \mu_0(d \thetab') \hat{\nu}(d \xb) \\
    \leq & C_{\nabla \varphi}^2 C_{\nabla \nabla \varphi}^2 \int_D \left | \Thetab_t(\thetab') - \Thetab_\infty(\thetab') \right |^2 \mu_0(d \thetab') \\
    & + C_{\nabla \nabla \varphi}^2 C_{\nabla \varphi}^2 \left | \Thetab_t(\thetab) - \Thetab_\infty(\thetab) \right |^2~.
    \end{split}
\end{equation}
Therefore,
\begin{equation}
    \begin{split}
        \| \bb_t - \bb_\infty \|_0^2 =& \EE_0 \int_D \left | \bb_t(\thetab) - \bb_\infty(\thetab) \right |^2 \mu_0(d \thetab) \\
        \leq & 2 C_{\nabla \varphi}^2 C_{\nabla \nabla \varphi}^2 \int_D \left | \Thetab_t(\thetab) - \Thetab_\infty(\thetab) \right |^2 \mu_0(d \thetab)~.
    \end{split}
\end{equation}
Recall that
\begin{equation}
    \begin{split}
        \left | \dot{\Thetab}_t(\thetab) \right | =& \left | \nabla V(\Thetab_t(\thetab), \mu_t) \right | \\
        =& \left | \int_\Omega \left ( f_t(\xb) - f_*(\xb) \right ) \nabla \varphi(\Thetab_t(\thetab), \xb) \hat{\nu}(d \xb) \right | \\
        \leq & C_{\nabla \varphi} \int_\Omega \left | f_t(\xb) - f_*(\xb) \right | \hat{\nu}(d \xb) \\
        \leq & \sqrt{2} C_{\nabla \varphi} \left ( \mathcal{L}(\mu_t) \right )^{1/2}
    \end{split}
\end{equation}
Thus,
\begin{equation}
\begin{split}
    \int_D \left | \Thetab_t(\thetab) - \Thetab_\infty(\thetab) \right |^2 \mu_0(d \thetab)
    = & \int_D \left | \int_t^\infty \dot{\Thetab}_s(\thetab) ds  \right |^2 \mu_0(d \thetab) \\
    \leq & \int_D \left ( \int_t^\infty \left | \dot{\Thetab}_s(\thetab) \right | ds \right )^2  \mu_0(d \thetab) \\
    \leq & 2 C_{\nabla \varphi}^2 \int_D \left ( \int_t^\infty \left ( \mathcal{L}(\mu_s) \right )^{1/2} ds \right )^2 \mu_0(d \thetab) \\
    \leq & 2 C_{\nabla \varphi}^2 \left ( \int_t^\infty \left ( \mathcal{L}(\mu_s) \right )^{1/2} ds \right )^2~.
\end{split}
\end{equation}
Hence, with the assumption of \eqref{eq:assump_unreg_maintext_loss}, we can conclude that
\begin{equation}
\begin{split}
    \int_0^\infty \| \bb_t - \bb_\infty \|_0 dt
    \leq & 2 C_{\nabla \varphi}^2 C_{\nabla \nabla \varphi} \int_0^\infty \int_t^\infty \left ( \mathcal{L}(\mu_s) \right )^{1/2} ds dt \\
    =& 2 C_{\nabla \varphi}^2 C_{\nabla \nabla \varphi} \int_0^\infty t \left ( \mathcal{L}(\mu_t) \right )^{1/2} dt \\
    < & \infty~.
\end{split}
\end{equation} \hfill $\square$

\subsubsection{Proof of Lemma \ref{lem:T_t_bounded_unreg}}
\label{app.pf_lemma_C2}
Our goal is to show that $\| \Tb_t \|_0$ remains bounded for all time. 
First note that, for all $t$, $\mathcal{A}_t^{(K)}$ is a positive semidefinite (PSD) operator on $\mathcal{V}(D)$ since
\begin{equation}
\begin{split}
    \langle \mathcal{A}_t^{(K)} \xib, \xib \rangle_{0} &= \EE_0 \int_{D\times D} \langle \xib(\thetab), \nabla \nabla' K(\Thetab_t(\thetab), \Thetab_t(\thetab'))  \xib(\thetab')\rangle \mu_0(d \thetab) \mu_0(d \thetab')  \\
    & = \EE_0 \int_\Omega \Big|\int_D \nabla \varphi(\Thetab_t(\thetab))\cdot  \xib(\thetab) \mu_0(d \thetab) \Big|^2 \hat{\nu}(d \xb) \geq  0~.
\end{split}
\end{equation}
Second, by Assumption \ref{ass:bounded_supp},
for $\mu_0$-almost-every $\thetab \in D$, $\Thetab_\infty(\thetab) = \lim_{t\to\infty} \Thetab_t(\thetab)$ exists, which allows us to define $\bb_\infty$, $\mathcal{A}_\infty^{(K)}$, and $\mathcal{A}_\infty^{(V)}$ similarly to \eqref{eq:b_t}, \eqref{eq:A_t_K} and \eqref{eq:A_t_V} by replacing $\Thetab_t(\cdot)$ with $\Thetab_\infty(\cdot)$. 
Since we assume that 
\begin{equation}
    \forall \xb_k \in \supp \hat \nu \quad : \quad f_\infty(\xb_k) = \int_D \varphi(\thetab, \xb_k) \mu_\infty(d \thetab) = f_*(\xb_k) 
\end{equation}
we have
\begin{equation}
    \forall \thetab \in D \quad : \quad \nabla \nabla V(\thetab, \mu_\infty) = \int_\Omega \nabla \nabla \varphi(\thetab, \xb) (f_\infty(\xb) - f_*(\xb)) d \xb = 0~.
\end{equation} 
This implies that $\mathcal{A}_\infty^{(V)}$ is the zero operator on $\mathcal{V}(D)$.

Third, we have the following observation:
\begin{lemma}
\label{lem:in_range}
Under Assumptions~\ref{ass:unit_1}, \ref{ass:iid_init} and \ref{ass:bounded_supp}, $\bb_t \in \mathrm{Ran}(\mathcal{A}_t^{(K)})$ for all $t$, and $\bb_\infty \in \mathrm{Ran}(\mathcal{A}_\infty^{(K)})$. Specifically, $\exists \tilde{\ub}_\infty \in \mathcal{V}(D)$ such that $\| \ub_\infty \|_0 < \infty$ and $\mathcal{A}_\infty^{(K)}  \tilde{\ub}_\infty = \bb_\infty$.
\end{lemma}

\noindent \textit{Proof of Lemma \ref{lem:in_range}:}
Recall from \eqref{eq:b_t_g_t} that $\bb_\infty = \mathcal{B}_\infty \bar{g}_\infty$. Define $\tilde{\ub}_\infty = \mathcal{B}_\infty \left (\mathcal{B}_\infty^\intercal \mathcal{B}_\infty \right )^\dag \bar{g}_\infty$. We claim that $\mathcal{A}_\infty^{(K)}  \tilde{\ub}_\infty = \bb_\infty$, because
\begin{equation}
    \begin{split}
        \mathcal{A}_\infty^{(K)} \Tilde{\ub}_\infty = & \left ( \mathcal{B}_\infty \mathcal{B}_\infty^\intercal \right ) \mathcal{B}_\infty  \left (\mathcal{B}_\infty^\intercal \mathcal{B}_\infty \right )^\dag \Bar{g}_\infty \\
        = & \mathcal{B}_\infty \mathcal{B}_\infty^\intercal \left ( \mathcal{B}_\infty \left ( \mathcal{B}_\infty \right )^\dag \right ) \left ( \mathcal{B}_\infty^\intercal \right )^\dag \bar{g}_\infty \\
        = & \mathcal{B}_\infty \left ( \mathcal{B}_\infty^\intercal \left ( \mathcal{B}_\infty^\intercal \right )^\dag \right ) \bar{g}_\infty \\
        =& \mathcal{B}_\infty \bar{g}_\infty \\
        =& \bb_\infty~,
    \end{split}
\end{equation}
where the third equality is because $\mathcal{B}_\infty \left ( \mathcal{B}_\infty \right )^\dag$ is the projection operator onto $\mathrm{Ran}(\mathcal{B}_\infty) = \mathrm{Nul}^\perp(\mathcal{B}_\infty^\intercal)$, and the fourth equality is because $\mathcal{B}_\infty^\intercal \left ( \mathcal{B}_\infty^\intercal \right )^\dag$ is the projection operator onto $\mathrm{Ran}(\mathcal{B}_\infty^\intercal) = \mathrm{Nul}^\perp(\mathcal{B}_\infty)$. 

It remains to establish that $\| \tilde{\ub}_\infty \|_0 < \infty$. To show this, we see that
\begin{equation}
    \begin{split}
    &
        \int_D |\tilde{\ub}_\infty(\thetab)|^2 \mu_0(d \thetab) \\
        &= \int_D \int_{\Omega\times\Omega} \Big (\nabla \varphi(\Thetab_\infty(\thetab), \xb) \big( \mathcal{M}_\infty^\dag \bar{g}_\infty \big)(\xb) \Big )\\
        & \qquad\qquad\qquad \cdot \Big ( \nabla \varphi(\Thetab_\infty(\thetab), \xb') \big( (\mathcal{M}_\infty)^\dag _\infty\bar{g}_\infty \big)(\xb') \Big ) \hat\nu(d\xb) \hat\nu(d \xb') \mu_0(d \thetab')  \\
        &= \int_\Omega \int_\Omega M(\xb, 
        \xb', \mu_\infty) \big( \mathcal{M}_\infty^\dag \bar{g}_\infty \big)(\xb) \big( \mathcal{M}_\infty^\dag \bar{g}_\infty \big)(\xb') \hat\nu(d\xb) \hat\nu(d \xb') \\
        &= \int_\Omega \big( \mathcal{M}_\infty^\dag \bar{g}_\infty \big)(\xb) \cdot \bar{g}_\infty(\xb) \hat\nu(d \xb) \\
        & \leq \lambda_{\min}^{-1} \int_\Omega |\bar{g}_\infty(\xb)|^2 \hat\nu(d \xb)~,
    \end{split}
\end{equation}
where $\lambda_{\min}$ is the least nonzero eigenvalue of the matrix $\mathcal{M}_\infty$ (and hence $\lambda_{\min}^{-1}$ is the largest eigenvalue of $\mathcal{M}_\infty^\dag$).
Since
\begin{equation}
\label{eq:beta_infty}
\begin{split}
    \mathbb{E}_{0} |\bar{g}_\infty(\xb)|^2 &= \mathbb{E}_{0} \Big | \int_D \varphi(\Thetab_\infty(\thetab), \xb) \omega_0(d \thetab) \Big |^2 \\
    &= \int_D \Big ( \varphi(\Thetab_\infty(\thetab), \xb) - \int_D \varphi(\Thetab_\infty(\thetab'), \xb) \mu_0(d \thetab') \Big )^2 \mu_0(d \thetab) \\
    &\leq \int_D \big | \varphi(\Thetab_\infty(\thetab), \xb) \big |^2 \mu_0(d \thetab)~,
\end{split}
\end{equation}
there is 
\begin{equation}
    \begin{split}
        \| \tilde{\ub}_\infty \|_{0}^2 &\leq  \mathbb{E}_{0} \int_D |\tilde{\ub}_\infty(\thetab)|^2 \mu_0(d \thetab) \\
        &\leq  \lambda_{\min}^{-1} \int_\Omega \int_D \big ( \varphi(\Thetab_\infty(\thetab), \xb) \big )^2 \mu_0(d \thetab) \nu(d \xb) \\
        &\leq  \lambda_{\min}^{-1}  C_{\varphi}^2 
        <  \infty~,
    \end{split}
\end{equation}
\hfill (\textit{End of the proof of Lemma \ref{lem:in_range}}) $\square$\\

\noindent Coming back to the prof of Lemma~\ref{lem:T_t_bounded_unreg}, we have shown that, as $t \to \infty$, \eqref{eq:T_t_simple} approaches the asymptotic dynamics
\begin{equation}
    \dot{\Tb}_t = - \mathcal{A}_\infty^{(K)} \Tb_t - \bb_\infty,
\end{equation}
with $\mathcal{A}_\infty^{(K)}$ positive semidefinite and $\bb_\infty$ in the range of $\mathcal{A}_\infty^{(K)}$. This is a stable system.
Hence, the rest of the task is to examine what happens at finite time.
To do so, we perform a change-of-variable with
\begin{equation}
    \zb_t = \Tb_t + \tilde{\ub}_\infty,
\end{equation}
with 
\begin{equation}
\label{eq:u_infty}
    \ub_\infty = \mathcal{B}_\infty (\mathcal{B}_\infty^\intercal \mathcal{B}_\infty)^\dag \bar{g}_\infty
\end{equation} as is defined in the proof of Lemma \ref{lem:in_range}. The dynamics of $\zb_t$ is governed by
\begin{equation}
    \begin{split}
        \dot{\zb}_t = \dot{\Tb}_t 
        =& - (\mathcal{A}_t^{(K)} + \mathcal{A}_t^{(V)}) \Tb_t - \bb_t \\
        =& - \mathcal{A}_t^{(K)} \zb_t - \mathcal{A}_t^{(V)} \zb_t - (\bb_t - (\mathcal{A}_t^{(K)} + \mathcal{A}_t^{(V)}) \tilde{\ub}_\infty) ~.
    \end{split}
\end{equation}
Thus, in integral form,
\begin{equation}
    \zb_t = \Pi(t, 0) \zb_0 + \int_0^t \Pi(t, s) \big ( - \mathcal{A}_s^{(V)} \zb_s - (\bb_s - (\mathcal{A}_s^{(K)} + \mathcal{A}_s^{(V)}) \tilde{\ub}_\infty) \big ) ds,
\end{equation}
where $\Pi(t, s)$ is the fundamental solution (a.k.a. Green's function) associated with the time-variant homogeneous system
\begin{equation}
    \dot{\zb}_t = - \mathcal{A}_t^{(K)} \zb_t~.
\end{equation}
Since $ \mathcal{A}_t^{(K)}$ is positive semidefinite for all $t$, there is $\| \Pi(t, s) \|_{0} \leq 1$ for $t > s$, where with a slight abuse of notation we also use $\| \cdot \|_0$  for the operator norm. Hence,
\begin{equation}
    \begin{split}
        \| \zb_t \|_0 \leq &
        \| \Pi(t, 0)\|_0  \| \zb_0 \|_0 + \int_0^t \| \Pi(t, s) \|_{0}  \Big ( \| \mathcal{A}_s^{(V)} \|_{0}  \| \zb_s \|_0 + \| \bb_s - (\mathcal{A}_s^{(K)} + \mathcal{A}_s^{(V)}) \tilde{\ub}_\infty \|_0 \Big ) ds \\
        \leq & \| \zb_0 \|_0 + \int_0^t \Big ( \| \mathcal{A}_s^{(V)} \|_{0}  \| \zb_s \|_0 + \| \bb_s - (\mathcal{A}_s^{(K)} + \mathcal{A}_s^{(V)}) \tilde{\ub}_\infty \|_0 \Big ) ds~.
    \end{split}
\end{equation}
By Gr\"onwall's inequality, we thus have
\begin{equation}
    \begin{split}
        \| \zb_t \|_0 \leq \Big ( \| \zb_0 \|_0 + \int_0^t \| \bb_s - (\mathcal{A}_s^{(K)} + \mathcal{A}_s^{(V)} ) \tilde{\ub}_\infty \|_0 ds \Big )  e^{\int_0^t \| \mathcal{A}_s^{(V)} \|_{0} ds}~.
    \end{split}
\end{equation}
Therefore, $\| \zb_t \|_0$ remains bounded for all time if we can show that 
\begin{equation}
\label{eq:integrable_168}
    \int_0^\infty \| \bb_t - (\mathcal{A}_t^{(K)} + \mathcal{A}_t^{(V)}) \tilde{\ub}_\infty \|_0 dt < \infty, \qquad \int_0^\infty \| \mathcal{A}_t^{(V)} \|_{0} dt < \infty~.
\end{equation}
Since
\begin{equation}
    \begin{split}
        \| \bb_t - (\mathcal{A}_t^{(K)} + \mathcal{A}_t^{(V)}) \tilde{\ub}_\infty \|_0 \leq \| \bb_t - \bb_\infty \|_0 + \| (\mathcal{A}_t^{(K)} - \mathcal{A}_\infty^{(K)}) \tilde{\ub}_\infty \|_0 + \| \mathcal{A}_\infty^{(V)} \tilde{\ub}_\infty \|_0
    \end{split}
\end{equation}
we see that \eqref{eq:integrable_168} is guaranteed by Lemmas \ref{lem:integrable_unreg} and \ref{lem:in_range}. 

This completes the proof of Lemma \ref{lem:T_t_bounded_unreg}.\hfill $\square$

\subsubsection{Proof of Lemma \ref{lem:g_infty_strict_bound_unreg}}
\label{app.pf_lemma_C3}
From \ref{lem:T_t_bounded_unreg}, we have that
\begin{equation}
    \lim_{t \to \infty} \left \| \fint_0^t \dot{\Tb}_s ds \right \|_0 = \lim_{t \to \infty} \left \| \frac{1}{t} \left ( \Tb_t - \Tb_0 \right ) ds \right \|_0= 0~.
\end{equation}
By \eqref{eq:T_t_simple},
we then obtain that
\begin{equation}
    \lim_{t \to \infty} \left \| \fint_0^t \left (\mathcal{A}_s^{(K)} \Tb_s + \bb_s \right ) ds + \fint_0^t \mathcal{A}_s^{(V)} \Tb_s ds \right \|_0 = 0~.
\end{equation}
By \eqref{eq:A_V_integrable} in Lemma \ref{lem:integrable_unreg} as well as Lemma \ref{lem:T_t_bounded_unreg},
we know that
\begin{equation}
     \lim_{t \to \infty} \left \| \fint_0^t \mathcal{A}_s^{(V)} \Tb_s ds \right \|_0 = 0~.
\end{equation}
Therefore,
\begin{equation}
    \lim_{t \to \infty} \left \| \fint_0^t \left (\mathcal{A}_s^{(K)} \Tb_s + \bb_s \right ) ds \right \|_0 = 0~.
\end{equation}
Next, 
by \eqref{eq:A_K_diff_integrable} and \eqref{eq:b_diff_integrable} in Lemma \ref{lem:integrable_unreg} as well as Lemma \ref{lem:T_t_bounded_unreg},
we know that
\begin{equation}
    \lim_{t \to \infty} \left \| \fint_0^t \left ( \mathcal{A}_s^{(K)} \Tb_s + \bb_s \right ) ds - \fint_0^t \left ( \mathcal{A}_\infty^{(K)} \Tb_s + \bb_\infty \right ) ds \right \|_0 = 0~.
\end{equation}
Therefore,
\begin{equation}
    \lim_{t \to \infty} \left \| \fint_0^t \left ( \mathcal{A}_\infty^{(K)} \Tb_s + \bb_\infty \right ) ds \right \|_0 = 0~.
\end{equation}
With $\tilde{\ub}_\infty$ defined in \eqref{eq:u_infty}, as $\bb_\infty = \mathcal{A}_\infty^{(K)} \ub_\infty$, there is
\begin{equation}
    \lim_{t \to \infty} \left \| \mathcal{A}_\infty^{(K)} \left ( \fint_0^t \Tb_s ds - \ub_\infty \right ) \right \|_0 = 0~.
\end{equation}
Let $ \xib^{||}$ denote the component of a vector field $\xib \in \mathcal{V}(D)$ that is in the range of $\mathcal{A}_\infty^{(K)}$. In the ERM setting, $\mathcal{A}_\infty^{(K)}$ has a least nonzero eigenvalue that is positive, and hence the above implies that
\begin{equation}
     \lim_{t \to \infty} \left \| \left ( \fint_0^t \Tb_s ds - \tilde{\ub}_\infty \right )^{||} \right \|_0 = 0
\end{equation}
or
\begin{equation}
    \lim_{t \to \infty} \left \| \left ( \fint_0^t \Tb_s ds \right )^{||} - \tilde{\ub}_\infty  \right \|_0 = 0
\end{equation}
and therefore, as $\mathrm{Nul}(\mathcal{A}_\infty^{(K)}) = \mathrm{Nul}(\mathcal{B}_\infty \mathcal{B}_\infty^\intercal) = \mathrm{Nul}(\mathcal{B}_\infty^\intercal)$, it follows that 
\begin{equation}
    \lim_{t \to \infty} \left \| \mathcal{B}_\infty^\intercal \left ( \fint_0^t \Tb_s ds \right ) - \mathcal{B}_\infty^\intercal \tilde{\ub}_\infty  \right \|_0 = 0~.
\end{equation}
Similar to \eqref{eq:A_K_diff_integrable}, it can be shown that $\int_0^\infty \| \mathcal{B}_t - \mathcal{B}_\infty \|_0 dt < \infty$. Therefore, we have
\begin{equation}
    \lim_{t \to \infty} \left \| \left ( \fint_0^t \mathcal{B}_s^\intercal \Tb_s ds \right ) - \mathcal{B}_\infty^\intercal \tilde{\ub}_\infty  \right \|_0 = 0~.
\end{equation}

Now,
\begin{equation}
    \begin{split}
        \fint_0^t \langle \Tb_s, \mathcal{A}_s^{(K)} \Tb_s \rangle_0 ds 
        = & \fint_0^t \langle \mathcal{B}_s^\intercal \Tb_s, \mathcal{B}_s^\intercal \Tb_s \rangle_{\hat{\nu}, 0} ds \\
        \geq & \left \langle \left ( \fint_0^t \mathcal{B}_s^\intercal \Tb_s ds \right ), \left ( \fint_0^t \mathcal{B}_s^\intercal \Tb_s ds \right ) \right \rangle_{\hat{\nu}, 0}~.
    \end{split}
\end{equation}
Hence,
\begin{equation}
\begin{split}
    \lim_{t \to \infty } \fint_0^t \langle \Tb_s, \mathcal{A}_s^{(K)} \Tb_s \rangle_0 ds \geq & \lim_{t \to \infty} \left \langle \left ( \fint_0^t \mathcal{B}_s^\intercal \Tb_s ds \right ), \left ( \fint_0^t \mathcal{B}_s^\intercal \Tb_s ds \right ) \right \rangle_{\hat{\nu}, 0} \\
    = & \left \langle \mathcal{B}_\infty^\intercal \tilde{\ub}_\infty, \mathcal{B}_\infty^\intercal \tilde{\ub}_\infty \right \rangle_{\hat{\nu}, 0} \\
    = & \left \langle \mathcal{B}_\infty^\intercal \left ( \mathcal{A}_\infty^{(K)} \right )^\dag \bb_\infty, \mathcal{B}_\infty^\intercal \left ( \mathcal{A}_\infty^{(K)} \right )^\dag \bb_\infty \tilde{\ub}_\infty \right \rangle_{\hat{\nu}, 0} \\
    = & \left \langle \left ( \mathcal{A}_\infty^{(K)} \right )^{\dag} \bb_\infty, \left ( \mathcal{A}_\infty^{(K)} \right ) \left ( \mathcal{A}_\infty^{(K)} \right )^{\dag} \bb_\infty \right \rangle_0 \\
    = & \left \langle \bb_\infty, \left ( \mathcal{A}_\infty^{(K)} \right )^{\dag} \bb_\infty \right \rangle_0~.
\end{split}
\end{equation}
\hfill $\square$
\subsubsection{Proof of Lemma \ref{lem:g_bar_in_range}}
\label{app.pf_lemma_C4}
Since
\begin{equation}
    \bar{g}_\infty(\xb) = \int_D \varphi(\thetab, \xb) \omega_0(d \thetab)~,
\end{equation}
we know that when viewed as a $L$-dimensional random vector, $\bar{g}_\infty$ has the distribution
\begin{equation}
    \bar{g}_\infty \sim \mathcal{N}(0, \bar{C}_\infty)~,
\end{equation}
where
\begin{equation}
\begin{split}
    \left ( \bar{C}_\infty \right )_{ij} :=& \EE_0 \left [ \bar{g}_\infty(\xb_i) \bar{g}_\infty(\xb_j) \right ]  \\
    =& \int_D \varphi(\thetab, \xb_i) \varphi(\thetab, \xb_j) \mu_\infty(d \thetab) - \int_D \varphi(\thetab, \xb_i) \mu_\infty(d \thetab) \int_D \varphi(\thetab', \xb_j) \mu_\infty(d \thetab')~,
\end{split}
\end{equation}
by the covariance of $\omega_0$, \eqref{eq:covomega0}. Thus, we decompose $\bar{C}_\infty$ as $\bar{C}_\infty = \bar{C}_\infty^{(1)} - \bar{C}_\infty^{(2)}$, with
\begin{equation}
    \left ( \bar{C}_\infty^{(1)} \right )_{ij} = \int_D \varphi(\thetab, \xb_i) \varphi(\thetab, \xb_j) \mu_\infty(d \thetab)~,
\end{equation}
\begin{equation}
    \left ( \bar{C}_\infty^{(2)} \right )_{ij} = \int_D \varphi(\thetab, \xb_i) \mu_\infty(d \thetab) \int_D \varphi(\thetab', \xb_j) \mu_\infty(d \thetab')~.
\end{equation}
Since $\bar{C}_\infty$ is PSD, its square root $\left ( \bar{C}_\infty \right )^{1/2}$ is well-defined. By the property of multivariate Gaussian, we can write
\begin{equation}
    \bar{g}_\infty \overset{\mathrm{d}}{=} \left ( \bar{C}_\infty \right )^{1/2} w~,
\end{equation}
where $\overset{\mathrm{d}}{=}$ denotes equality in distribution, and $w \in \mathbb{R}^n$ follows the distribution
\begin{equation}
    w \sim \mathcal{N}(0, \mathrm{Id}_n)~.
\end{equation}
This means that almost surely, $\bar{g}_\infty \in \mathrm{Ran}\left( \left ( \bar{C}_\infty \right )^{1/2} \right)$, and which would imply that $\bar{g}_\infty \in \mathrm{Ran}\left( \bar{C}_\infty  \right)$. This means that almost surely, we can write
\begin{equation}
    \bar{g}_\infty = \bar{C}_\infty^{(1)} w^{(1)} - \bar{C}_\infty^{(2)} w^{(2)}
\end{equation}
for some pair of $w^{(1)}, w^{(2)} \in \mathbb{R}^n$. Our goal is then to show that both $\bar{C}_\infty^{(1)} w^{(1)}$ and $\bar{C}_\infty^{(2)} w^{(2)}$ belong to $\mathrm{Ran}(\mathcal{B}_\infty^\intercal)$. Here, under Assumption \ref{ass:shallownn}, since $\varphi(\thetab, \xb) = c \hat{\varphi}(\zb, \xb)$ when $\thetab = \begin{bmatrix}c & \zb\end{bmatrix}^\intercal$, there is
\begin{equation}
    \nabla \varphi(\thetab, \xb) = \begin{bmatrix} \hat{\varphi}(\zb, \xb) \\ c \nabla_{\zb} \hat{\varphi}(\zb, \xb) \end{bmatrix}~.
\end{equation}
Therefore, first, we have
\begin{equation}
\begin{split}
    \left ( \bar{C}_\infty^{(1)} w^{(1)} \right )_i =& \int_D \varphi(\thetab, \xb_i) \left ( \sum_{j=1}^n \varphi(\thetab, \xb_j) w^{(1)}_j \right ) \mu_\infty(d \thetab) \\
    = & \int_D \nabla \varphi(\thetab, \xb_i)^\intercal \cdot \begin{bmatrix} c(\thetab) \left ( \sum_{j=1}^n \varphi(\thetab, \xb_j) w^{(1)}_j \right ) \\ 0 \end{bmatrix} \mu_\infty(d \thetab) \\
    = & \mathcal{B}_\infty^\intercal \xib^{(1)}~,
\end{split}
\end{equation}
with
\begin{equation}
    \xib(\thetab)^{(1)} = \begin{bmatrix} c(\thetab) \left ( \sum_{j=1}^n \varphi(\thetab, \xb_j) w^{(1)}_j \right ) \\ 0 \end{bmatrix} ~.
\end{equation}
This means that $\left ( \bar{C}_\infty^{(1)} w^{(1)} \right ) \in \mathrm{Ran}(\mathcal{B}_\infty^\intercal)$.

Second, there is
\begin{equation}
\begin{split}
    \left ( \bar{C}_\infty^{(2)} w^{(2)} \right )_i =& \left ( \int_D \varphi(\thetab, \xb_i) \mu_\infty(d \thetab) \right ) \left ( \sum_{j=1}^n w^{(2)}_j \int_D \varphi(\thetab', \xb_j) \mu_\infty(d \thetab') \right ) \\
    = & \int_D \nabla \varphi(\thetab, \xb_i)^\intercal \cdot \begin{bmatrix} c(\thetab) \left ( \sum_{j=1}^n w^{(2)}_j \int_D \varphi(\thetab', \xb_j) \mu_\infty(d \thetab') \right ) \\ 0 \end{bmatrix} \mu_\infty(d \thetab) \\
    = & \mathcal{B}_\infty^\intercal \xib^{(2)}~,
\end{split}
\end{equation}
with
\begin{equation}
    \xib(\thetab)^{(2)} = \begin{bmatrix} c(\thetab) \left ( \sum_{j=1}^n w^{(2)}_j \int_D \varphi(\thetab', \xb_j) \mu_\infty(d \thetab') \right ) \\ 0 \end{bmatrix}
\end{equation}
This means that $\left ( \bar{C}_\infty^{(2)} w^{(2)} \right ) \in \mathrm{Ran}(\mathcal{B}_\infty^\intercal)$. Hence the lemma is proved.\hfill $\square$

\subsection{Proof of Theorem~\ref{thm:min_eig_ddV} (Under assumptions on the curvature in the long-time)}
\label{app:curvature_2}
When the limiting measure $\mu_\infty$ does not necessarily interpolate the training data, such as in the regularized case, we have the following condition on $\Tb_t$ which guarantees that~\eqref{eq:mcbound_lem33} holds:
\begin{lemma}
\label{lem:34}
If 
\begin{equation}
\label{eq:T_t_ddV_cond}
    \lim_{T \to \infty} \EE_{0} \int_0^T \int_D \langle \Tb_t(\thetab),\nabla \nabla V(\Thetab_t(\thetab),
    \mu_t) \Tb_t(\thetab)\rangle\mu_0(d\thetab) dt \geq 0~,
\end{equation}
(including when this limit is $+\infty$) then~\eqref{eq:mcbound_lem33} holds.
\end{lemma}
\noindent
\textit{Proof of Lemma~\ref{lem:34}:}
With $\mathfrak{D}_t$ defined in \eqref{eq:diff_term}, for \eqref{eq:mcbound_lem33} to hold, it is sufficient to show that 
\begin{equation}
    \lim_{T \to \infty} \fint_0^T \mathfrak{D}_t dt \leq 0~.
\end{equation}
Recall from \eqref{eq:diff_term_3_terms} that
\begin{equation}
    \begin{split}
        \fint_0^T \mathfrak{D}_t dt = - \frac{1}{T}\| \Tb_T \|_0^2 - 2 \fint_0^T \langle \Tb_t, \mathcal{A}_t^{(V)} \Tb_t \rangle_0 dt - \fint_0^T \langle \Tb_t, \mathcal{A}_t^{(K)} \Tb_t \rangle_0 dt~.
    \end{split}
\end{equation}
Since $\Tb_0 = 0$ and $\mathcal{A}_t^{(K)}$ is PSD, we see that the assumption \eqref{eq:T_t_ddV_cond} is sufficient.
\hfill $\square$

Note that  condition~\eqref{eq:T_t_ddV_cond} is natural since we know from Proposition~\ref{th:ltflow} that $\lim_{t \to \infty} \nabla \nabla V(\Thetab_t(\thetab), \mu_t) = \nabla \nabla V (\Thetab_\infty(\thetab),\mu_\infty)$ exists and is positive semidefinite  $\mu_0$-almost surely. 
This lemma then allows us to prove Theorem~\ref{thm:min_eig_ddV}:
\textit{Proof of Theorem~\ref{thm:min_eig_ddV}:}
Our goal is to verify \eqref{eq:T_t_ddV_cond} in order to apply Lemma \ref{lem:34}. We first see that
\begin{equation}
\label{eq:T_t_ddV_bound}
    \begin{split}
        & \EE_{0}\int_D \langle \Tb_t(\thetab),\nabla \nabla V(\Thetab_t(\thetab),
    \mu_t) \Tb_t(\thetab)\rangle\mu_0(d\thetab) \\
    \geq & \EE_0 \int_D \lambda_{\min}(\nabla \nabla V(\Thetab_t(\thetab), \mu_t)) |\Tb_t(\thetab)|^2 \mu_0(d \thetab) \\
    \geq & \EE_0 \int_D \min \left \{\lambda_{\min}(\nabla \nabla V(\Thetab_t(\thetab), \mu_t)), 0 \right \} |\Tb_t(\thetab)|^2 \mu_0(d \thetab) \\
    = & \int_D \min \left \{\lambda_{\min}(\nabla \nabla V(\Thetab_t(\thetab), \mu_t)), 0 \right \} \left ( \EE_0 |\Tb_t(\thetab)|^2 \right ) \mu_0(d \thetab) \\
    \geq & \int_D \min \left \{\lambda_{\min}(\nabla \nabla V(\Thetab_t(\thetab), \mu_t)), 0 \right \} \left  ( \sup_{\thetab \in \supp \mu_0} \EE_0 |\Tb_t(\thetab)|^2 \right ) \mu_0(d \thetab) \\
    \geq & \| \Tb_t \|_{\sup}^2 \left ( \int_D \min \left \{\lambda_{\min}(\nabla \nabla V(\Thetab_t(\thetab), \mu_t)), 0 \right \} \mu_0(d \thetab) \right )~,
    \end{split}
\end{equation}
where we define, for $\xib \in \mathcal{V}(D)$,
\begin{equation}
    \| \xib \|_{\sup} := \sup_{\thetab \in \supp \mu_0} \Big ( \EE_0  |\xib(\thetab)|^2 \Big )^{1/2}~,
\end{equation}
which is a norm on $\mathcal{V}(D)$. 

Hence, if we assume that $\left | \int_D \min \left \{\lambda_{\min}(\nabla \nabla V(\thetab, \mu_t)), 0 \right \} \mu_0(d \thetab) \right |$ is small asymptotically, then what remains is to upper-bound $\| \Tb_t \|_{\sup}$. Recall from \eqref{eq:T_t_simple} that the dynamics of $\Tb_t$ is governed by
\begin{equation}
\begin{split}
    \dot{\Tb}_t = - (\mathcal{A}_t^{(K)} + \mathcal{A}_t^{(V)}) \Tb_t - \bb_t,
\end{split}
\end{equation}
Thus, in the $\| \cdot \|_{\sup}$ norm defined above, we have
\begin{equation}
\label{eq:T_t_sup_dot}
\begin{split}
    \frac{d}{dt} \| \Tb_t \|_{\sup} \leq & \| - (\mathcal{A}_t^{(K)} + \mathcal{A}_t^{(V)}) \Tb_t - \bb_t \|_{\sup} \\
    \leq & \| \mathcal{A}_t^{(K)} \Tb_t \|_{\sup} + \| \mathcal{A}_t^{(V)} \Tb_t \|_{\sup} + \| \bb_t \|_{\sup}~.
\end{split}
\end{equation}
We then want to bound the growth of $\| \Tb_t \|_{\sup}$ by upper-bounding the RHS.
Note that for $\xib \in \mathcal{V}(D)$,
\begin{equation}
\label{eq:A_V_sup_norm}
\begin{split}
    \| \mathcal{A}_t^{(V)} \xib \|_{\sup}^2 
    = & \sup_{\thetab \in D} \EE_0 |(\mathcal{A}_t^{(V)} \xib)(\thetab)|^2 \\
    =& \sup_{\thetab \in D} \EE_0 | \nabla \nabla V(\Thetab_t(\thetab), \mu_t) \xib(\thetab) |^2 \\
    \leq & \sup_{\thetab \in D} | \nabla \nabla V(\Thetab_t(\thetab), \mu_t)|^2 \EE_0 |\xib(\thetab) |^2 \\
    \leq & (C_{\nabla \nabla \varphi} C_{\varphi} + \lambda)^2 \sup_{\thetab \in D} \EE_0 |\xib(\thetab) |^2 \\
    = & (C_{\nabla \nabla \varphi} C_{\varphi} + \lambda)^2 \| \xib \|_{\sup}^2~,
\end{split}
\end{equation}
\begin{equation}
\label{eq:A_K_sup_norm}
\begin{split}
    \| \mathcal{A}_t^{(K)} \xib \|_{\sup}^2 
    = & \sup_{\thetab \in D} \EE_0 |(\mathcal{A}_t^{(K)} \xib)(\thetab)|^2 \\
    =& \sup_{\thetab \in D} \EE_0 |\int_D \nabla' \nabla K(\Thetab_t(\thetab), \Thetab_t(\thetab')) \xib(\thetab') \mu_0(d \thetab')|^2 \\
    \leq & \sup_{\thetab \in D} \EE_0 \int_D | \nabla' \nabla K(\Thetab_t(\thetab), \Thetab_t(\thetab')) |^2 |\xib(\thetab')|^2 \mu_0(d \thetab') \\
    \leq & \sup_{\thetab \in D} (C_{\nabla \varphi})^4 \int_D \EE_0 |\xib(\thetab')|^2 \mu_0(d \thetab') \\
    \leq & (C_{\nabla \varphi})^4 \sup_{\thetab' \in D} \EE_0 |\xib(\thetab')|^2 \\
    = & (C_{\nabla \varphi})^4 \| \xib \|_{\sup}^2~.
    \end{split}
\end{equation}
Thus, 
\begin{equation}
    \| \mathcal{A}_t^{(K)} \Tb_t \|_{\sup} + \| \mathcal{A}_t^{(V)} \Tb_t \|_{\sup} \leq ( C_{\nabla \varphi}^2 + C_{\nabla \nabla \varphi} C_{\varphi} + \lambda ) \| \Tb_t \|_{\sup}~.
\end{equation}
To bound $\| \bb_t \|_{\sup}$, we recall that
\begin{equation}
    \begin{split}
        \bb_t(\thetab) =& \int_D \nabla K (\Thetab_t(\thetab), \Thetab_t(\thetab')) \omega_0(d \thetab') \\
        =& \int_\Omega \nabla \varphi(\Thetab_t(\thetab), \xb) \bar{g}_t(\xb) \hat{\nu}(d \xb)~,
    \end{split}
\end{equation}
with
\begin{equation}
    \bar{g}_t(\xb) = \int_D \varphi(\Thetab_t(\thetab), \xb) \omega_0(d \thetab)~.
\end{equation}
This implies that $\forall \thetab \in \supp \mu_0$,
\begin{equation}
    |\bb_t(\thetab)| \leq \frac{1}{n} C_{\nabla \varphi} \sum_{l=1}^n |\bar{g}_t(\xb_l)|
\end{equation}
and so
\begin{equation}
\begin{split}
    \EE_0 |\bb_t(\thetab)|^2 \leq & C_{\nabla \varphi}^2 \EE_0 \left ( \frac{1}{n} \sum_{l=1}^n |\bar{g}_t(\xb_l)| \right )^2 \\
    \leq &  C_{\nabla \varphi}^2 \EE_0 \left ( \frac{1}{n} \sum_{l=1}^n |\bar{g}_t(\xb_l)|^2 \right ) \\
    \leq & C_{\nabla \varphi}^2 \frac{1}{n} \sum_{l=1}^n \EE_0 |\bar{g}_t(\xb_l)|^2~.
\end{split}
\end{equation}
On the other hand, similar to \eqref{eq:beta_infty}, we have
\begin{equation}
\begin{split}
    \mathbb{E}_{0} |\bar{g}_t(\xb)|^2 &= \mathbb{E}_{0} \Big | \int_D \varphi(\Thetab_t(\thetab), \xb) \omega_0(d \thetab) \Big |^2 \\
    &= \int_D \Big ( \varphi(\Thetab_t(\thetab), \xb) - \int_D \varphi(\Thetab_t(\thetab'), \xb) \mu_0(d \thetab') \Big )^2 \mu_0(d \thetab) \\
    &\leq \int_D \big | \varphi(\Thetab_t(\thetab), \xb) \big |^2 \mu_0(d \thetab)\\
    &\leq (C_\varphi)^2~,
\end{split}
\end{equation}
Thus, there is $\forall \thetab \in \supp \mu_0$,
\begin{equation}
    \EE_0 |\bb_t(\thetab)|^2 \leq (C_{\nabla \varphi})^2 (C_{\varphi})^2
\end{equation}
and so
\begin{equation}
    \| \bb_t \|_{\sup} \leq C_{\nabla \varphi} C_{\varphi}~.
\end{equation}
Therefore, based on \eqref{eq:T_t_sup_dot}, we have
\begin{equation}
    \frac{d}{dt} \| \Tb_t \|_{\sup} \leq ( (C_{\nabla \varphi})^2 + C_{\nabla \nabla \varphi} C_{\varphi} + \lambda ) \| \Tb_t \|_{\sup} + C_{\nabla \varphi} C_{\varphi}~.
\end{equation}
Since $\Tb_0 = 0$, we thus have
\begin{equation}
    \begin{split}
        \| \Tb_t \|_{\sup} \leq &  C_{\nabla \varphi} C_{\varphi} \int_0^t e^{( (C_{\nabla \varphi})^2 + C_{\nabla \nabla \varphi} C_{\varphi} + \lambda )(t-s)} ds \\
        = & C_{\nabla \varphi} C_{\varphi} e^{( (C_{\nabla \varphi})^2 + C_{\nabla \nabla \varphi} C_{\varphi} + \lambda )t} \int_0^t e^{-( (C_{\nabla \varphi})^2 + C_{\nabla \nabla \varphi} C_{\varphi} + \lambda)s}ds \\
        \leq & \frac{C_{\nabla \varphi} C_{\varphi}}{(C_{\nabla \varphi})^2 + C_{\nabla \nabla \varphi} C_{\varphi} + \lambda} e^{( (C_{\nabla \varphi})^2 + C_{\nabla \nabla \varphi} C_{\varphi} + \lambda )t}
    \end{split}
\end{equation}
Now, using \eqref{eq:T_t_ddV_bound}, we see that in order for \eqref{eq:T_t_ddV_cond} to hold, it is sufficient to have
\begin{equation}
    \lim_{t \to \infty}  e^{( (C_{\nabla \varphi})^2 + C_{\nabla \nabla \varphi} C_{\varphi} + \lambda ) t} \left ( \int_D \min \left \{\lambda_{\min}(\nabla \nabla V(\thetab, \mu_t)), 0 \right \} \mu_0(d \thetab) \right ) = 0
\end{equation}
and therefore sufficient to have
\begin{equation}
    -\int_D \min \left \{\lambda_{\min}(\nabla \nabla V(\thetab, \mu_t)), 0 \right \} \mu_0(d \thetab) \sim O \left (e^{- ((C_{\nabla \varphi})^2 + C_{\nabla \nabla \varphi} C_{\varphi} + \lambda)t} \right )
\end{equation}
\hfill $\square$

To intuitively  understand~\eqref{eq:curv}, note that we know from \eqref{eq:Vlim} in Proposition~\ref{th:ltflow} that $\Lambda_t(\thetab)\to0$ $\mu_0$-almost surely as $t\to\infty$. Condition~\eqref{eq:curv} can therefore be satisfied  by having $\Lambda_t(\thetab)$ converge to zero sufficiently fast in the regions of $D$ where it is negative, or having the measure of these regions with respect to $\mu_0$ converge to zero sufficiently fast, or both.
\subsection{Proof of Theorem~\ref{thm:long_time_reg_maintext} (Regularized case)}
\label{app:long_time_reg}

Recall from Proposition \ref{prop:dclt2} that the dynamics of $g_t$ is governed by
\begin{equation}
\label{eq:gtvolt_app}
    g_t(\xb) +\int_0^t \int_\Omega \Gamma_{t,s}(\xb,\xb') g_s(\xb') \hat\nu(d\xb') ds= \bar g_t(\xb)~,
\end{equation}
with
\begin{equation}
    \Gamma_{t,s} (\xb,\xb') = \int_D \langle \nabla \varphi(\Thetab_t(\thetab), \xb), J_{t,s} (\thetab) \nabla \varphi(\Thetab_s(\thetab), \xb')\rangle \mu_0(d\thetab)~,
\end{equation}
with $J_{t,s}$ being the Jacobian of the flow $\Thetab_t$.

In the ERM setting, $\supp \hat{\nu}$ is singular, thus we have $\hat{\nu}(d \xb) = n^{-1} \sum_{l=1}^n \delta_{\xb_l}(d \xb)$, where $n$ is the total number of training data points. We define $\mathcal{W}_n(\Omega)$ together with the inner product $\langle \cdot, \cdot \rangle_{\hat{\nu}, 0}$ and the norm $\| \cdot \|_{\hat{\nu}, 0}$ as in Appendix \ref{app:notations}. We will also continue to consider $g_t$ and $\bar{g}_t$ equivalently as $n$-dimensional vectors,
\begin{equation}
    \begin{pmatrix} g_t(\xb_1) & \cdots & g_t(\xb_n) \end{pmatrix}^T, \qquad \begin{pmatrix} \bar{g}_t(\xb_1) & \cdots & \bar{g}_t(\xb_n) \end{pmatrix}^T~,
\end{equation}
respectively. Thus, $\Gamma_{t, s}$ can also be represented by the $n \times n$ matrix
\begin{equation}
    \begin{pmatrix} \Gamma_{t, s}(\xb_1, \xb_1) & \cdots &  \Gamma_{t, s}(\xb_1, \xb_n) \\
    \vdots & & \vdots \\
    \Gamma_{t, s}(\xb_n, \xb_1) & \cdots & \Gamma_{t, s}(\xb_n, \xb_n) \end{pmatrix}~.
\end{equation}
Under such an abuse of notations, we can simplify \eqref{eq:gtvolt_app} into
\begin{equation}
\label{eq:gtvolt_vec}
    g_t + \int_0^t \Gamma_{t, s} g_s ds = \bar{g}_t~.
\end{equation}
Thus, the goal is to prove that
\begin{equation}
\label{eq:limsup_g_t}
    \lim_{t \to \infty} \sup \fint_0^t \EE_0 \|{g}_t\|_{\hat{\nu}}^2 dt \leq \EE_0 \|\bar{g}_\infty \|_{\hat{\nu}}^2~.
\end{equation}
As in \eqref{eq:limGamm}, we also define
\begin{equation}
    \Gamma^\infty_{t - s} (\xb,\xb') = \int_D \langle \nabla \varphi(\thetab,\xb), e^{-(t-s) \nabla \nabla V_\infty(\thetab)} \nabla \varphi(\thetab,\xb')\rangle \mu_\infty(d\thetab)~,
\end{equation}
where for simplicity, we write $V_t(\cdot)$ for $V(\cdot, \mu_t)$ and $V_\infty(\cdot)$ for $V(\cdot, \mu_\infty)$.
Then the heuristic argument outlined in Section \ref{sec:g_infty} before Theorem \ref{prop:ginfty} amounts to rewriting \eqref{eq:gtvolt_vec} as
\begin{equation}
\label{eq:volterra_our}
    g_t + \int_0^t \Gamma^{\infty}_{t - s} g_s ds = \bar{g}_t + \int_0^t (\Gamma^{\infty}_{t - s} - \Gamma_{t, s}) g_s ds
\end{equation}
and then arguing that 1) $\Gamma^\infty$ is a nonnegative convolution-type Volterra kernel, and 2) the second term on the RHS is small. Rigorously, we need to introduce an extra level of complication: for every $t_0 > 0$, we can rewrite \eqref{eq:gtvolt_vec} into
\begin{equation}
    \begin{split}
        g_t &=  \bar{g}_t - \int_{t_0}^t \Gamma_{t, s} g_s ds - \int_0^{t_0} \Gamma_{t, s} g_s ds \\
        &= \bar{g}_t - \int_{t_0}^t \Gamma^{\infty}_{t-s} g_s ds + \int_{t_0}^t (\Gamma^{\infty}_{t-s} - \Gamma_{t, s}) g_s ds - \int_0^{t_0} \Gamma_{t, s} g_s ds~.
    \end{split}
\end{equation}
Then, for any $T > t_0$, by multiplying $g_t$ and integrating from $t_0$ to $T$, we get
\begin{equation}
\begin{split}
\label{eq:g_t_L2_orig}
    & \int_{t_0}^T \| g_t \|_{\hat \nu}^2 dt + \int_{t_0}^T \int_{t_0}^t \langle g_t, \Gamma^{\infty}_{t-s} g_s \rangle_{\hat \nu} ds dt \\
    \leq & \int_{t_0}^T \langle g_t, \bar{g}_t \rangle_{\hat \nu} dt + \int_{t_0}^T \langle g_t, \int_{t_0}^t (\Gamma^{\infty}_{t-s} - \Gamma_{t, s}) g_s ds \rangle_{\hat \nu} dt + \int_{t_0}^T \langle g_t, \int_0^{t_0} \Gamma_{t, s} g_s  ds \rangle_{\hat \nu} dt~.
\end{split}
\end{equation}

Then firstly, the second term on the LHS is nonnegative because of the nonnegativity of $\Gamma_t^\infty$ as a convolution-type Volterra kernel, as proven in Appendix \ref{app:long_time_reg_asymp}.

Hence, we have
\begin{equation}
\begin{split}
\label{eq:g_t_L2_4_terms}
    \int_{t_0}^T \| g_t \|_{\hat \nu}^2 dt
    &\leq  \int_{t_0}^T \langle g_t, \bar{g}_t \rangle_{\hat \nu} dt + \int_{t_0}^T \langle g_t, \int_{t_0}^t (\Gamma^{\infty}_{t-s} - \Gamma_{t, s}) g_s ds \rangle_{\hat \nu} dt \\
    & \quad + \int_{t_0}^T \left\langle g_t, \int_0^{t_0} \Gamma_{t, s} g_s  ds \right\rangle_{\hat \nu} dt~.
\end{split}
\end{equation}

By Cauchy-Schwartz,
\begin{equation}
    \int_{t_0}^T \langle g_t, \bar{g}_t \rangle_{\hat \nu} dt \leq \left( \int_{t_0}^T \| g_t \|_{\hat \nu}^2 dt \right)^\frac{1}{2}  \left( \int_{t_0}^T \| \bar{g}_t \|_{\hat \nu}^2 dt \right)^\frac{1}{2}~,
\end{equation}
\begin{equation}
\begin{split}
    &\int_{t_0}^T \left\langle g_t, \int_{t_0}^t (\Gamma^{\infty}_{t-s} - \Gamma_{t, s}) g_s ds \right\rangle_{\hat \nu} dt \\
    \leq & \left(  \int_{t_0}^T \| g_t \|_{\hat \nu}^2 dt \right)^{\frac{1}{2}}  \left( \int_{t_0}^T \left \| \int_{t_0}^t (\Gamma^{\infty}_{t-s} - \Gamma_{t, s})  g_s ds \right\|_{\hat \nu}^2 dt \right)^\frac{1}{2} \\
    \leq & \left( \int_{t_0}^T \| g_t \|_{\hat \nu}^2 dt \right)  \left( \int_{t_0}^T \int_{t_0}^t \| \Gamma^{\infty}_{t-s} - \Gamma_{t, s} \|_{\hat \nu}^2 ds dt \right)^{\frac{1}{2}} ~, 
\end{split}
\end{equation}
and
\begin{equation}
    \begin{split}
        &\int_{t_0}^T \left\langle g_t, \int_0^{t_0} \Gamma_{t, s} g_s  ds \right\rangle_{\hat \nu} dt\\
        \leq & \left(  \int_{t_0}^T \| g_t \|_{\hat \nu}^2 dt \right)^{\frac{1}{2}}  \left( \int_{t_0}^T \left\| \int_{0}^{t_0} \Gamma_{t, s} g_s ds \right\|_{\hat \nu}^2 dt \right)^\frac{1}{2} \\
        \leq & \left(\int_{t_0}^T \| g_t \|_{\hat \nu}^2 dt \right)^{\frac{1}{2}}  \left( \int_{t_0}^T \left( \int_0^{t_0} \| \Gamma_{t, s} \|_{\hat \nu}^2 ds \right)  \left( \int_0^{t_0} \| g_s \|_{\hat \nu}^2 ds \right) dt \right)^\frac{1}{2} \\
        \leq & \left(  \int_{t_0}^T \| g_t \|_{\hat \nu}^2 dt \right)^{\frac{1}{2}}  \left(  \int_{0}^{t_0} \| g_t \|_{\hat \nu}^2 dt \right)^{\frac{1}{2}}  \left( \int_{t_0}^T \int_0^{t_0} \| \Gamma_{t, s} \|_{\hat \nu}^2 ds dt \right)^\frac{1}{2}~.
    \end{split}
\end{equation}
Therefore, putting everything together, we have
\begin{equation}
    \begin{split}
        \left( \int_{t_0}^T \| g_t \|_{\hat \nu}^2 dt \right)^\frac{1}{2} \leq & \left( \int_{t_0}^T \| \bar{g}_t \|_{\hat \nu}^2 dt \right)^\frac{1}{2}\\
        &+ \left(  \int_{t_0}^T \| g_t \|_{\hat \nu}^2 dt \right)^{\frac{1}{2}} \left( \int_{t_0}^T \int_{t_0}^t \| \Gamma^{\infty}_{t-s} - \Gamma_{t, s} \|_{\hat \nu}^2 ds dt \right)^{\frac{1}{2}} \\
        & + \left( \int_0^{t_0} \| g_t \|_{\hat \nu}^2 dt \right)^{\frac{1}{2}} \left( \int_{t_0}^T \int_0^{t_0} \| \Gamma_{t, s} \|_{\hat{\nu}}^2 ds dt \right)^\frac{1}{2} ~,
    \end{split}
\end{equation}
and hence, using $\fint_a^b \cdot \ dt$ to denote the averaged integral $\frac{1}{b-a} \int_a^b \cdot \ dt$,
\begin{equation}
    \begin{split}
        \left( \fint_{t_0}^T \| g_t \|_{\hat \nu}^2 dt \right)^\frac{1}{2} \leq & \left( \fint_{t_0}^T \| \bar{g}_t \|_{\hat \nu}^2 dt \right)^\frac{1}{2}\\& + \left(  \fint_{t_0}^T \| g_t \|_{\hat \nu}^2 dt \right)^{\frac{1}{2}} \left( \int_{t_0}^T \int_{t_0}^t \| \Gamma^{\infty}_{t-s} - \Gamma_{t, s} \|_{\hat \nu}^2 ds dt \right)^{\frac{1}{2}} \\
        & + \left( \int_0^{t_0} \| g_t \|_{\hat \nu}^2 dt \right)^{\frac{1}{2}} \left( \fint_{t_0}^T \int_0^{t_0} \| \Gamma_{t, s} \|_{\hat{\nu}}^2 ds dt \right)^\frac{1}{2} ~,
    \end{split}
\end{equation}

or
\begin{equation}
\label{eq:g_t_L2_avg}
    \begin{split}
        & \left( 1 - \left[ \int_{t_0}^T \int_{t_0}^t \| \Gamma^{\infty}_{t-s} - \Gamma_{t, s} \|_{\hat \nu}^2 ds dt \right]^{\frac{1}{2}} \right) \left( \fint_{t_0}^T \| g_t \|_{\hat \nu}^2 dt \right)^\frac{1}{2} \\
        \leq & \left( \fint_{t_0}^T \| \bar{g}_t \|_{\hat \nu}^2 dt \right)^\frac{1}{2} 
        + \left( \int_0^{t_0} \| g_t \|_{\hat \nu}^2 dt \right)^{\frac{1}{2}} \left( \fint_{t_0}^T \int_0^{t_0} \| \Gamma_{t, s} \|_{\hat{\nu}}^2 ds dt \right)^\frac{1}{2}~.
    \end{split}
\end{equation}
\begin{lemma}
\label{lem:integrable_reg}
Under all assumptions in Theorem \ref{thm:long_time_reg_maintext} except for \eqref{eq:unif_rate_char} being replaced by
a weaker condition, 
\begin{equation}
\label{eq:assump_reg_maintext}
     \int_0^\infty \int_D\big ( |\Thetab_t(\thetab) - \Thetab_\infty(\thetab)| + |U_t(\thetab)|^2 \big ) e^{C_1 (U_t(\thetab) + \bar{U}_t)}  \mu_0(d \thetab) dt < \infty~,
\end{equation}
we have
\begin{equation}
    \lim_{t_0 \to \infty} \int_{t_0}^\infty \int_{t_0}^t \| \Gamma^{\infty}_{t-s} - \Gamma_{t, s} \|_{\hat \nu}^2 ds dt = 0
\end{equation}
and
$\forall t_0 > 0$,
\begin{equation}
    \lim_{T \to \infty} \fint_{t_0}^T \int_0^{t_0} \| \Gamma_{t, s} \|^2 ds dt = 0~.
\end{equation}
We will prove in Appendix~\ref{app.interpret_reg} that \eqref{eq:unif_rate_char} indeed implies \eqref{eq:assump_reg_maintext}. 
\end{lemma}
The lemma will be proved in Appendix \ref{app.pf_lemma_C8}, and let us first proceed with the proof of the theorem assuming this lemma. Suppose for contradiction that \eqref{eq:limsup_g_t} does not hold, meaning that
\begin{equation}
\label{eq:contra_assump}
    \lim_{T \to \infty} \sup \big ( \fint_0^T \|g_t\|_{\hat \nu}^2 dt \big )^\frac{1}{2} = \|\bar{g}_\infty \|_{\hat \nu} + \epsilon
\end{equation}
for some $\epsilon > 0$. We will select a pair of $t_0$ and $T$ for which the inequality \eqref{eq:g_t_L2_avg} cannot be satisfied.
Firstly, by the convergence of $\bar g_t$ to $\bar g_\infty$, $\exists t_a > 0$ such that $\forall t_1, t_2 > t_a$, 
\begin{equation}
    \left( \fint_{t_1}^{t_2} \|\bar{g}_t\|_{\hat \nu}^2 dt \right)^\frac{1}{2} \leq \|\bar{g}_\infty\|_{\hat \nu} + \tfrac{1}{6}\epsilon~.
\end{equation}
Secondly, by our assumption \eqref{eq:contra_assump} and the first part of Lemma \ref{lem:integrable_reg}, $\exists t_0 > t_a$ such that both
\begin{equation}
\label{eq:131}
    \left( \fint_{0}^{t_0} \|g_t\|_{\hat \nu}^2 dt \right)^\frac{1}{2} \leq \|\bar{g}_\infty \|_{\hat \nu} + 2 \epsilon
\end{equation}
and 
\begin{equation}
    \int_{t_0}^\infty \int_{t_0}^t \| \Gamma^{\infty}_{t-s} - \Gamma_{t, s} \|_{\hat \nu}^2 ds dt < \frac{\epsilon}{6 \|\bar{g}_\infty\|_{\hat \nu} + 3 \epsilon}
\end{equation}
are satisfied.
In particular, \eqref{eq:131} implies
\begin{equation}
    \left( \int_{0}^{t_0} |g_t|^2 dt \right)^\frac{1}{2} \leq t_0^\frac{1}{2} \cdot ( |\bar{g}_\infty| + 2 \epsilon )
\end{equation}
Let 
\begin{equation}
    \delta = \left( \frac{\epsilon}{6 t_0^\frac{1}{2} \cdot (\|\bar{g}_\infty\|_{\hat \nu} + 2
    \epsilon)} \right)^2 > 0~.
\end{equation}
By the second part of Lemma \ref{lem:integrable_reg}, $\exists t_b > t_0$ such that $\forall T > t_b$,
\begin{equation}
    \fint_{t_0}^T \int_0^{t_0} \| \Gamma_{t, s} \|^2 ds dt < \delta 
\end{equation}
so that the last term in \eqref{eq:g_t_L2_avg} satisfies
\begin{equation}
    \left( \int_0^{t_0} \| g_t \|_{\hat \nu}^2 dt \right)^{\frac{1}{2}}  \left( \fint_{t_0}^T \int_0^{t_0} \| \Gamma_{t, s} \|^2 ds dt \right)^\frac{1}{2} < \tfrac{1}{6}\epsilon
\end{equation}
By our assumption \eqref{eq:contra_assump}, we can choose a $T > t_b$ such that
\begin{equation}
    \big ( \fint_{0}^{T} \|g_t\|_{\hat \nu}^2 dt \big )^\frac{1}{2} \geq \|\bar{g}_\infty\|_{\hat \nu} + \tfrac{2}{3}\epsilon~.
\end{equation}
Since 
\begin{equation}
    \left( \fint_{0}^{t_0} \|g_t\|_{\hat \nu}^2 dt \right)^\frac{1}{2} \leq \|\bar{g}_\infty\|_{\hat \nu} + {2}\epsilon,
\end{equation}
we can assume without loss of generality that $\frac{T}{t_0}$ is large enough so that
\begin{equation}
    \big ( \fint_{t_0}^{T} \|g_t\|_{\hat \nu}^2 dt \big )^\frac{1}{2} \geq \|\bar{g}_\infty\|_{\hat \nu} + \tfrac{1}{2} \epsilon~.
\end{equation}
Thus, back to the inequality \eqref{eq:g_t_L2_avg}, 
the LHS is strictly lower-bounded by
\begin{equation}
\begin{split}
    \big ( \|\bar{g}_\infty\|_{\hat \nu} + \tfrac{1}{2} \epsilon \big ) \big( 1 - \frac{\epsilon}{6 \|\bar{g}_\infty\|_{\hat \nu} + 3 \epsilon} \big ) = & \|\bar{g}_\infty\|_{\hat \nu} + \tfrac{1}{3}\epsilon,
\end{split}
\end{equation}
whereas the RHS is strictly upper-bounded by
\begin{equation}
    \|\bar{g}_\infty\|_{\hat \nu} + \tfrac{1}{6} \epsilon + \tfrac{1}{6} \epsilon = \|\bar{g}_\infty\|_{\hat \nu} + \tfrac{1}{3} \epsilon~.
\end{equation}
This gives contradiction and we are done with the proof of Theorem~\ref{thm:long_time_reg_maintext}. \hfill $\square$

\subsubsection{Proof of Lemma~\ref{lem:integrable_reg}}
\label{app.pf_lemma_C8}
It remains to prove Lemma~\ref{lem:integrable_reg}. To do so we will need an auxiliary result, that we state and prove first:
\begin{lemma}
\label{lem.delta_gamma}
Let $\Delta \Gamma_{t, s} := \Gamma_{t, s} - \Gamma_{t-s}^\infty$.
If $\nabla \nabla V$ is uniformly positive definite with eigenvalues lower-bounded by $\lambda$, then there exists constants $C$ and $C'$ whose values depend on $|D'|$, $C_\varphi$, $C_{\nabla \varphi}$, $C_{\nabla \nabla \varphi}$, and $L_{\nabla \nabla \varphi}$ such that
\begin{equation}
    \| \Delta \Gamma_{t, s} \|_{\hat{\nu}} \leq C e^{-\lambda(t-s)} \int_D \left(| \Delta\Thetab_t(\thetab)| + \left (| \Delta\Thetab_s(\thetab) | + U_s(\thetab) \right) e^{C' (U_s(\thetab) + \bar{U}_s)} \right )
    \mu_0(d\thetab)
\end{equation}
where $\Delta\Thetab_t(\thetab)=\Thetab_t(\thetab) - \Thetab_\infty(\thetab)$.
\end{lemma}

\textit{Proof of Lemma \ref{lem.delta_gamma}:}
To bound $\| \Delta \Gamma_{t, s} \|_{\hat \nu}$, we bound $\| \Delta \Gamma_{t, s} \eta \|_{\hat \nu}$ for $\eta \in \mathbb{R}^n$. 
Note that $\Delta \Gamma_{t, s} \eta$ can be obtained in the following way. Consider the two systems
\begin{equation}
    \left\{\begin{aligned}
    \frac{d}{dt} \xib_t(\thetab) &= - \nabla \nabla V_t(\Thetab_t(\thetab))  \xib_t(\thetab) \\
    \xib_s(\thetab) &= \int_\Omega \nabla \varphi(\Thetab_s(\thetab), \xb') \eta (\xb') \hat{\nu}(d \xb')
    \end{aligned}\right.
\end{equation}
\begin{equation}
    \left\{\begin{aligned}
    \frac{d}{dt} {\xib}'_t(\thetab) &= - \nabla \nabla V_\infty(\Thetab_\infty(\thetab))  \xib'_t(\thetab) \\
    \xib'_s(\thetab) &= \int_\Omega \nabla \varphi(\Thetab_\infty(\thetab), \xb') \eta (\xb') \hat{\nu}(d \xb')
    \end{aligned}\right.
\end{equation}
Then there is
\begin{equation}
\begin{split}
    (\Gamma_{t, s} \eta) (\xb) =& \int_D \nabla \varphi(\Thetab_t(\thetab), \xb) \cdot\xib_t(\thetab) \mu_0(d \thetab) \\
    (\Gamma_{t-s}^\infty \eta) (\xb) =& \int_D \nabla \varphi(\Thetab_\infty(\thetab), \xb) \cdot\xib'_t(\thetab) \mu_0(d \thetab)
\end{split}
\end{equation}
and hence
\begin{equation}
\begin{split}
    (\Delta \Gamma_{t, s} \eta)(\xb) =& \int_D \nabla \varphi(\Thetab_t(\thetab), \xb) \xib_t(\thetab) \mu_0(d \thetab) - \int_D \nabla \varphi(\Thetab_\infty(\thetab), \xb) \xib'_t(\thetab) \mu_0(d \thetab) \\
     =& 
    \int_D \nabla \varphi(\Thetab_t(\thetab), \xb)\cdot \big (\xib_t(\thetab) - \xib'_t(\thetab) \big ) \mu_0(d \thetab) \\
      & + \int_D \big( \nabla \varphi(\Thetab_t(\thetab), \xb) - \nabla \varphi(\Thetab_\infty(\thetab), \xb) \big )\cdot \xib'_t(\thetab) \mu_0(d \thetab)~.
\end{split}
\end{equation}
We will first try to bound $\xib_t(\thetab) - \xib'_t(\thetab)$ as a function of $\eta$. Define $\Delta \xib_t(\thetab) = \xib_t(\thetab) - \xib'_t(\thetab)$. Then
\begin{equation}
\begin{split}
    \frac{d}{dr}{\Delta \xib}_r(\thetab) =& - \big (\nabla \nabla V_r(\Thetab_r(\thetab)) - \nabla \nabla V_\infty(\Thetab_\infty(\thetab)) \big ) \xib_r - \nabla \nabla V_\infty(\Thetab_\infty(\thetab)) \Delta \xib_r(\thetab) \\
    =& - \nabla \nabla V_\infty(\Thetab_\infty(\thetab)) \Delta \xib_r(\thetab) \\
    & - \big (\nabla \nabla V_r(\Thetab_r(\thetab)) - \nabla \nabla V_\infty(\Thetab_\infty(\thetab)) \big ) \xib'_r \\
    & - \big (\nabla \nabla V_r(\Thetab_r(\thetab)) - \nabla \nabla V_\infty(\Thetab_\infty(\thetab)) \big ) \Delta \xib_r~.
\end{split}
\end{equation}
Thus,
\begin{equation}
\begin{split}
    \Delta \xib_t(\thetab) &= e^{- (t-s) \nabla \nabla V_\infty(\Thetab_\infty(\thetab))} \Delta \xib_s(\thetab) \\
    & + \int_s^t e^{- (t-r) \nabla \nabla V_\infty(\Thetab_\infty(\thetab))} \big (\nabla \nabla V_r(\Thetab_r(\thetab)) - \nabla \nabla V_\infty(\Thetab_\infty(\thetab)) \big ) \xib'_r(\thetab) dr \\
    & + \int_s^t e^{- (t-r) \nabla \nabla V_\infty(\Thetab_\infty(\thetab))} \big (\nabla \nabla V_r(\Thetab_r(\thetab)) - \nabla \nabla V_\infty(\Thetab_\infty(\thetab)) \big ) \Delta \xib_r(\thetab) dr~.
\end{split}
\end{equation}
Since $\nabla \nabla V_\infty(\Thetab_\infty(\thetab)) - \lambda I_d$ is positive semidefinite, we first have
\begin{equation}
    | \xib'_r(\thetab) | \leq e^{-\lambda(r-s)} | \xib'_s(\thetab) |
\end{equation}
as well as
\begin{equation}
\begin{split}
    | \Delta \xib_t(\thetab) |
    \leq &\, e^{-\lambda(t-s)} | \Delta \xib_s(\thetab) | \\
    & + \int_s^t e^{-\lambda(t-r)}\| \nabla \nabla V_r(\Thetab_r(\thetab)) - \nabla \nabla V_\infty(\Thetab_\infty(\thetab)) \| | \xib'_r(\thetab) | dr \\
    & + \int_s^t e^{-\lambda(t-r)} \| \nabla \nabla V_r(\Thetab_r(\thetab)) - \nabla \nabla V_\infty(\Thetab_\infty(\thetab)) \| | \Delta \xib_r(\thetab) | dr \\
    \leq &\, e^{-\lambda(t-s)} | \Delta \xib_s(\thetab) | \\
    & + \int_s^t e^{-\lambda(t-s)} \| \nabla \nabla V_r(\Thetab_r(\thetab)) - \nabla \nabla V_\infty(\Thetab_\infty(\thetab)) \| | \xib'_s(\thetab) | dr \\
    & + \int_s^t e^{-\lambda(t-r)} \| \nabla \nabla V_r(\Thetab_r(\thetab)) - \nabla \nabla V_\infty(\Thetab_\infty(\thetab)) \| | \Delta \xib_r(\thetab) | dr~.
\end{split}
\end{equation}
To prepare for an application of Gronwall's inequality, we introduce a change-of-variable by defining, for $r \in [s, t]$, 
\begin{equation}
    \overline{\Delta \xib}_r(\thetab) = e^{\lambda(t-s)} \Delta \xib_r(\thetab)~.
\end{equation}
Then we can rewrite the equation above as
\begin{equation}
\begin{split}
    | \overline{\Delta \xib}_t(\thetab) | =&\, e^{\lambda(r-s)} | \Delta \xib_t(\thetab) |  \\
    \leq &\, | \Delta \xib_s(\thetab) | + \int_s^t \| \nabla \nabla V_r(\Thetab_r(\thetab)) - \nabla \nabla V_\infty(\Thetab_\infty(\thetab)) \| | \xib'_s(\thetab) | dr \\
    & + \int_s^t e^{\lambda(r-s)} \| \nabla \nabla V_r(\Thetab_r(\thetab)) - \nabla \nabla V_\infty(\Thetab_\infty(\thetab)) \| | \Delta \xib_r(\thetab) | dr \\
    \leq &\, | \overline{\Delta \xib}_s(\thetab) |  + \int_s^t \| \nabla \nabla V_r(\Thetab_r(\thetab)) - \nabla \nabla V_\infty(\Thetab_\infty(\thetab)) \| | \xib'_s(\thetab) | dr \\
    & + \int_s^t \| \nabla \nabla V_r(\Thetab_r(\thetab)) - \nabla \nabla V_\infty(\Thetab_\infty(\thetab)) \| | \overline{\Delta \xib}_r(\thetab) | dr ~.
\end{split}
\end{equation}
Thus, by Gronwall's inequality,
\begin{equation}
\begin{split}
    | \overline{\Delta \xib}_t(\thetab) | \leq &\,
    \Big ( | \overline{\Delta \xib}_s(\thetab) |
    + \int_s^t \| \nabla \nabla V_r(\Thetab_r(\thetab)) - \nabla \nabla V_\infty(\Thetab_\infty(\thetab)) \| | \xib'_s(\thetab) | dr \Big ) \\
    & \times e^{\int_s^t \| \nabla \nabla V_r(\Thetab_r(\thetab)) - \nabla \nabla V_\infty(\Thetab_\infty(\thetab)) \| dr}~,
\end{split}
\end{equation}
or, back in the original variable that we are interested in,
\begin{equation}
    \begin{split}
    | \Delta \xib_t(\thetab) | \leq &\,
    \Big ( | \Delta \xib_s(\thetab) |
    + \int_s^t \| \nabla \nabla V_r(\Thetab_r(\thetab)) - \nabla \nabla V_\infty(\Thetab_\infty(\thetab)) \| | \xib'_s(\thetab) | dr \Big ) \\
    & \times e^{-\lambda (t-s)+\int_s^t \| \nabla \nabla V_r(\Thetab_r(\thetab)) - \nabla \nabla V_\infty(\Thetab_\infty(\thetab)) \| dr}~.
\end{split}
\end{equation}

Now, we have
\begin{equation}
    \begin{split}
        & |\Delta \Gamma_{t, s} \eta(\xb) |\\
        \leq & \,
        \|\int_D \nabla \varphi(\Thetab_t(\thetab), \xb)^\intercal \cdot \Delta \xib_t(\thetab) \mu_0(d \thetab) \|_{\hat \nu} \\
        & + \| \int_D \big( \nabla \varphi(\Thetab_t(\thetab), \xb) - \nabla \varphi(\Thetab_\infty(\thetab), \xb) \big )^\intercal \xib'_t(\thetab) \mu_0(d \thetab) \|_{\hat \nu} \\
        \leq &\, C_{\nabla \varphi}  \int_D | \Delta \xib_t(\thetab) | \mu_0(d \thetab) + C_{\nabla \nabla \varphi} \int_D |\Delta\Thetab_t(\thetab)| | \xib'_t(\thetab) | \mu_0(d \thetab) \\
        \leq &\, C_{\nabla \varphi}  e^{-\lambda(t-s)}  \int_D \Big ( | \Delta \xib_s(\thetab) |
    + \int_s^t \| \nabla \nabla V_r(\Thetab_r(\thetab)) - \nabla \nabla V_\infty(\Thetab_\infty(\thetab)) \| | \xib'_s(\thetab) | dr \Big ) \\
    & \hspace{65pt}  e^{\int_s^t \| \nabla \nabla V_r(\Thetab_r(\thetab)) - \nabla \nabla V_\infty(\Thetab_\infty(\thetab)) \| dr} \mu_0(d \thetab) \\
    & + C_{\nabla \nabla \varphi}  e^{-\lambda(t-s)} \int_D |\Delta\Thetab_t(\thetab) | | \xib'_s(\thetab) | \mu_0(d \thetab)~.
    \end{split}
\end{equation}
Note that we have,
\begin{equation}
    | \xib'_s(\thetab) | = | \int_\Omega \nabla \varphi(\Thetab_\infty(\thetab), \xb') \eta (\xb') \hat{\nu}(d \xb') | \leq C_{\nabla \varphi} \sup_{1 \leq p \leq P} |\eta(\xb_p)| \leq P^\frac{1}{2} C_{\nabla \varphi} \| \eta \|_{\hat \nu},
\end{equation}
\begin{equation}
\begin{split}
    | \Delta \xib_s(\thetab) | =&\, |\int_\Omega \big ( \nabla \varphi(\Thetab_s(\thetab), \xb) - \nabla \varphi(\Thetab_\infty(\thetab), \xb) \big ) \eta(\xb') \hat{\nu}(d \xb')| \\
    \leq & \, \int_\Omega | \nabla \varphi(\Thetab_s(\thetab), \xb') - \nabla \varphi(\Thetab_\infty(\thetab), \xb') | | \eta(\xb') | \hat{\nu}(d \xb') \\
    \leq &\, n^\frac{1}{2} C_{\nabla \nabla \varphi} |\Delta\Thetab_s(\thetab) | \| \eta \|_{\hat \nu} 
\end{split}
\end{equation}
and, since $\nabla \nabla V_r(\thetab) = \int_\Omega \nabla \nabla \varphi(\thetab, \xb) (f_r(\xb) - f_*(\xb)) \hat{\nu}(d \xb)$ and  $\nabla \nabla V_\infty(\thetab) = \int_\Omega \nabla \nabla \varphi(\thetab, \xb) (f_\infty(\xb) - f_*(\xb)) \hat{\nu}(d \xb)$,
\begin{equation}
\begin{split}
    &\| \nabla \nabla V_r(\Thetab_r(\thetab)) - \nabla \nabla V_\infty(\Thetab_\infty(\thetab)) \|\\
    \leq &\, \| \nabla \nabla V_r(\Thetab_r(\thetab)) - \nabla \nabla V_r(\Thetab_\infty(\thetab)) \| \\
    & + \| \nabla \nabla V_r(\Thetab_\infty(\thetab)) - \nabla \nabla V_\infty(\Thetab_\infty(\thetab)) \| \\
    \leq &\, L_{\nabla \nabla \varphi} C_\varphi| \Delta\Thetab_r(\thetab)| + C_{\nabla \nabla \varphi} \| \Delta f_r \|_{\hat{\nu}, \infty},
\end{split}
\end{equation}
where we use $\| f \|_{\hat{\nu}, \infty}$ to denote $\sup_{\xb \in \supp \hat{\nu}} |f(\xb)|$ and we defined $\Delta f_t = f_t - f_\infty$.

As a result, we have
\begin{equation}
    \begin{split}
        &\| \Delta \Gamma_{t, s} \|_{\hat \nu}\\
        \leq &\, \| \eta \|_{\hat \nu}^{-1} \|\Delta \Gamma_{t, s} \eta \|_{\hat \nu} \\ 
        \leq &\, C_{\nabla \varphi}  e^{-\lambda(t-s)}  \int_D \Big ( C_{\nabla \nabla \varphi} |\Delta\Thetab_s(\thetab)  | +\int_s^t \big ( L_{\nabla \nabla \varphi}   C_\varphi  | \Delta\Thetab_r(\thetab)  | 
     + C_{\nabla \nabla \varphi} \| \Delta f_r  \|_{\hat \nu, \infty} \big )  C_{\nabla \varphi} dr \Big ) \\
    & \hspace{95pt} \times e^{\int_s^t L_{\nabla \nabla \varphi}  C_\varphi  | \Delta\Thetab_r(\thetab)  |
     + C_{\nabla \nabla \varphi} \| \Delta f_r  \|_{\hat \nu, \infty} dr} \mu_0(d \thetab) \\
    & + C_{\nabla \nabla \varphi}  e^{-\lambda(t-s)}  \int_D C_{\nabla \varphi} |\Delta\Thetab_t(\thetab) |   \mu_0(d \thetab)~.
    \end{split}
\end{equation}

Therefore, using $C_0$, $C_1$, etc. to represent constants that depend on $C_\varphi$, $C_{\nabla \varphi}$, $C_{\nabla \nabla \varphi}$, $C_{\nabla \nabla \varphi}$ and $L_{\nabla \nabla \varphi}$, we have
\begin{equation}
    \begin{split}
        & \| \Delta \Gamma_{t, s} \|_{\hat \nu}\\
        \leq & \, C_0  e^{-\lambda(t-s)}  \Big ( \int_D  | \Delta\Thetab_t(\thetab)  | \mu_0(d \thetab) + \int_D | \Delta\Thetab_s(\thetab)  |  e^{C_1 \int_s^t | \Delta\Thetab_r(\thetab)  | + \| \Delta f_r  \|_{\hat \nu, \infty} dr} \mu_0(d \thetab) \\
        & \hspace{45pt} + \int_D \big ( \int_s^t | \Delta \Thetab_r(\thetab)  | + \| f_r - f_\infty \|_{\hat \nu, \infty}  dr \big )  e^{C_1  \int_s^t | \Delta \Thetab_r(\thetab)  | + \| \Delta f_r \|_{\hat \nu, \infty} dr} \mu_0(d \thetab) \Big )~.
    \end{split}
\end{equation}
Note that $\| \Delta f_r  \|_{\hat \nu, \infty}$ can be further upper-bounded by $C_\varphi \int_D | \Delta\Thetab_r(\thetab)  | \mu_0(d \thetab)$.
Furthermore, defining
\begin{equation}
    \overline{\Delta\Thetab_t } = \int_D | \Delta\Thetab_t(\thetab)| \mu_0(d \thetab) 
\end{equation}
we can write the bound above as
\begin{equation}
    \begin{split}
        \| \Delta \Gamma_{t, s} \|_{\hat \nu}
        \leq &\, C_0  e^{-\lambda(t-s)}  \Big ( \int_D  | \Delta\Thetab_t(\thetab)  | \mu_0(d \thetab)  + \int_D | \Delta\Thetab_s(\thetab)  |  e^{C_1  \int_s^t | \Delta\Thetab_r(\thetab) | + \overline{\Delta\Thetab_r } dr} \mu_0(d \thetab) \\
        & \hspace{45pt} + \int_D \big ( \int_s^t | \Delta\Thetab_r(\thetab)  | + \overline{\Delta\Thetab_r }  dr \big )  e^{C_1 \int_s^t | \Delta\Thetab_r(\thetab) | + \overline{\Delta\Thetab_r} dr} \mu_0(d \thetab) \Big )~.
    \end{split}
\end{equation}

Finally, let
\begin{equation}
    U_t(\thetab) = \int_t^\infty | \Delta\Thetab_t(\thetab)  | dt
\end{equation}
and 
\begin{equation}
    \bar{U}_t = \int_D U_t(\thetab) \mu_0(d \thetab) = \int_t^\infty \overline{\Delta\Thetab_t } dt~.
\end{equation}
Then there is
\begin{equation}
\begin{split}
    \| \Delta \Gamma_{t, s} \|_{\hat \nu} \leq &\,  C_0  e^{-\lambda(t-s)}  \int_D \left(| \Delta\Thetab_t(\thetab) | + \left (| \Delta\Thetab_s(\thetab)  | + U_s(\thetab) + \bar{U}_s \right)  e^{C_1  (U_s(\thetab) + \bar{U}_s)} \right )
    \mu_0(d\thetab) \\
    \leq & \, 2 C_0  e^{-\lambda(t-s)}  \int_D \left(| \Delta\Thetab_t(\thetab) | + \left (| \Delta\Thetab_s(\thetab)  | + U_s(\thetab)\right)  e^{C_1  (U_s(\thetab) + \bar{U}_s)} \right )
    \mu_0(d\thetab)~.
\end{split}
\end{equation}
\hfill \textit{(End of the proof of Lemma \ref{lem.delta_gamma}.)} $\square$

\textit{Proof of Lemma \ref{lem:integrable_reg}:}
Lemma \ref{lem.delta_gamma} entails that, $\exists C, C' > 0$ such that
\begin{equation}
\begin{split}
    \| \Delta \Gamma_{t, s} \|_{\hat \nu}^2 \leq & \,C e^{-2\lambda(t-s)} \Bigg ( \int_D \Bigg(| \Delta\Thetab_t(\thetab) | + \Big (| \Delta\Thetab_s(\thetab)  | + U_s(\thetab) \Big) e^{C' (U_s(\thetab) + \bar{U}_s)} \Bigg )
    \mu_0(d\thetab) \Bigg )^2 \\
    \leq &\, 4 C e^{-2\lambda(t-s)} \int_D |\Delta\Thetab_t(\thetab)|^2 + \Big ( |\Delta\Thetab_s(\thetab) |^2 + U_s(\thetab)^2 \Big ) e^{2 C' (U_s(\thetab) + \bar{U}_s)} \mu_0(d \thetab) \\
    \leq &\, 4 C |D'| e^{-2\lambda(t-s)} \int_D |\Delta\Thetab_t(\thetab) | + \Big ( |\Delta\Thetab_s(\thetab) | + U_s(\thetab)^2 \Big ) e^{2 C' (U_s(\thetab) + \bar{U}_s)} \mu_0(d \thetab), 
\end{split}
\end{equation}
where for the last inequality, we assume that $| D' | \geq 1$ (or, to accommodate the more general case, just replace $| D'|$ by $\max \{ |D'|, 1 \}$).

To prove Lemma \ref{lem:integrable_reg}, the first goal is to show
\begin{equation}
    \lim_{t_0 \to \infty} \int_{t_0}^\infty \int_{t_0}^t \| \Delta \Gamma_{t, s} \|_{\hat \nu}^2 ds dt = 0~.
\end{equation}
There is
\begin{equation}
    \begin{split}
        & \int_{t_0}^\infty \int_{t_0}^t \| \Delta \Gamma_{t, s} \|_{\hat \nu}^2 ds dt \\
        \leq &\, 4 C |D'| \int_D \int_{t_0}^\infty \int_{t_0}^t e^{- 2 \lambda (t-s)} \bigg ( | \Delta\Thetab_t(\thetab)  | + \Big (| \Delta\Thetab_s(\thetab)  | + U_s(\thetab)^2 \Big)  e^{2 C'  (U_s(\thetab) + \bar{U}_s)} \bigg ) ds dt \mu_0(d \thetab) \\
        \leq &\, 4 C |D'| \int_D \bigg ( \int_{t_0}^\infty \Big ( \int_{t_0}^t e^{-2\lambda(t-s)} ds \Big ) | \Delta\Thetab_t(\thetab) | dt  \\
        & \hspace{30pt} + \int_{t_0}^\infty \Big ( \int_{s}^\infty e^{-2\lambda(t-s)} dt \Big ) \Big (| \Delta\Thetab_s(\thetab)  | + U_s(\thetab)^2 \Big)  e^{2 C' (U_s(\thetab) + \bar{U}_s)} ds \bigg ) \mu_0(d \thetab) \\
        \leq &\,  2 C |D'| \lambda^{-1} \int_D \Big ( \int_{t_0}^\infty | \Delta\Thetab_t(\thetab)  | dt + \int_{t_0}^\infty \Big (| \Delta\Thetab_s(\thetab) | + U_s(\thetab)^2 \Big)  e^{2 C' (U_s(\thetab) + \bar{U}_s)} ds \Big ) \mu_0(d \thetab) \\
        \leq & \, 4 C |D'| \lambda^{-1} \int_D \int_{t_0}^\infty \Big ( | \Delta\Thetab_s(\thetab)| + U_s(\thetab)^2 \Big ) e^{2 C'(U_s(\thetab) + \bar{U}_s)} ds \mu_0(d \thetab)~.
    \end{split}
\end{equation}
By our assumption, the RHS is finite for $t_0 > 0$. Hence, by taking $t_0$ large enough, the value of $\int_{t_0}^\infty \int_{t_0}^t \| \Delta \Gamma_{t, s} \|_{\hat \nu}^2 ds dt$ can be made arbitrarily close to zero.

The second goal is to show that
$\forall t_0 > 0$,
\begin{equation}
    \lim_{T \to \infty} \fint_{t_0}^T \int_0^{t_0} \| \Gamma_{t, s} \|^2 ds dt = 0~.
\end{equation}
As a first step, we show that
\begin{equation}
    \lim_{T \to \infty} \fint_{t_0}^T \int_0^{t_0} \| \Gamma_{t-s}^\infty \|_{\hat \nu}^2 ds dt = 0
\end{equation}
because $\forall \eta \in \mathcal{W}_L(\Omega)$, there is
\begin{equation}
\begin{split}
    |\langle \eta, \Gamma_{t-s}^\infty \eta \rangle_{\hat \nu}| =& \int_D \big \langle \bb(\thetab), e^{- t \nabla \nabla V_\infty(\Thetab_\infty(\thetab))} \bb(\thetab)\Big ) \big \rangle \mu_0(d \thetab) \\
    \leq & e^{-\lambda(t-s)} \int_D |\bb(\thetab)|^2 \mu_0(d \thetab) \\
    \leq & e^{-\lambda(t-s)} \| \mathcal{M}_\infty \|_{\hat \nu} \| \eta \|_{\hat \nu}^2~,
\end{split}
\end{equation}
where
\begin{equation}
    \bb(\thetab) = \int_\Omega \nabla \varphi(\Thetab_\infty(\thetab), \xb) \eta(\xb) \hat{\nu}(d\xb).
\end{equation}
and $\mathcal{M}_\infty$ is defined as $\mathcal{M}_\infty := \mathcal{B}_\infty^\intercal \mathcal{B}_\infty$, or concretely, for $\eta \in \mathcal{W}_L(\omega)$,
\begin{equation}
\begin{split}
    (\mathcal{M}_\infty \eta)(\xb) :=& \int_\Omega \Big ( \int_D \nabla \varphi(\Thetab_\infty(\thetab'), \xb)^\intercal \nabla \varphi(\Thetab_\infty(\thetab'), \xb') \mu_0(d \thetab')  \Big ) \eta(\xb') \hat{\nu}(d \xb') \\
    =&  \int_\Omega M(\xb, \xb', \mu_\infty) \eta(\xb') \hat{\nu}(d \xb')~,
\end{split}
\end{equation}
where
\begin{equation}
     M(\xb, \xb', \mu_\infty) := \int_D \nabla \varphi(\Thetab_\infty(\thetab'), \xb) \cdot  \nabla \varphi(\Thetab_\infty(\thetab'), \xb') \mu_0(d \thetab')~.
\end{equation}
In the ERM setting, $\mathcal{M}_\infty$ is effectively an $L \times L$ matrix. 
Thus, 
\begin{equation}
\begin{split}
    \fint_{t_0}^T \int_0^{t_0} \| \Gamma_{t-s}^\infty \|_{\hat \nu}^2 ds dt \leq &  \fint_{t_0}^T \int_0^{t_0} e^{-2\lambda(t - s)} \| \mathcal{M}_\infty \|_{\hat \nu}^2 ds dt \\
    \leq & \| \mathcal{M}_\infty \|_{\hat \nu}^2 \fint_{t_0}^T e^{-2\lambda(t-t_0)} dt \to  0 \qquad \text{as \ \ $T\to\infty$}
\end{split}
\end{equation}
Hence, it is sufficient to show that
\begin{equation}
    \lim_{T \to \infty} \fint_{t_0}^T \int_0^{t_0} \| \Delta \Gamma_{t, s} \|^2 ds dt = 0~.
\end{equation}
We have
\begin{equation}
    \begin{split}
        & \int_{t_0}^T \int_0^{t_0} \| \Delta \Gamma_{t, s} \|^2 ds dt \\
        \leq &\, 4 C |D'| \int_D \bigg( \int_{t_0}^T \Big ( \int_0^{t_0} e^{-2 \lambda(t-s)} ds \Big ) |\Delta\Thetab_t(\thetab)| dt \\
        & \hspace{20pt} + \int_0^{t_0} \Big ( \int_{t_0}^T e^{-2\lambda(t-s)} dt \Big ) \Big ( | \Delta\Thetab_s(\thetab) | + U_s(\thetab)^2 \Big ) e^{2 C' (U_s(\thetab) + \bar{U}_s)} ds \bigg ) \mu_0(d \thetab) \\
        \leq &\, 2 C |D'| \lambda^{-1} \int_D \bigg ( \int_{t_0}^T e^{-2\lambda (t - t_0)} |\Delta\Thetab_t(\thetab) | dt \\
        & \hspace{40pt} + \int_0^{t_0} e^{-2\lambda(t_0 - s)} \Big ( | \Delta\Thetab_s(\thetab) | + U_s(\thetab)^2 \Big ) e^{2 C' (U_s(\thetab) + \bar{U}_s)} ds \bigg ) \mu_0(d \thetab) \\
        \leq &\, 4 C |D'| \lambda^{-1} \int_D \int_0^{\infty} \Big ( |\Delta\Thetab_s(\thetab) | + U_s(\thetab)^2 \Big ) e^{2 C' (U_s(\thetab) + \bar{U}_s)} ds \mu_0(d \thetab)  \\
        < & \infty
    \end{split}
\end{equation}
by assumption \eqref{eq:assump_reg_maintext}. Therefore, 
\begin{equation}
    \begin{split}
        \fint_{t_0}^T \int_{0}^{t_0} \| \Delta \Gamma_{t, s} \|_{\hat \nu}^2 ds dt = \frac{1}{T-t_0} \int_{t_0}^T \int_{0}^{t_0} \| \Delta \Gamma_{t, s} \|_{\hat \nu}^2 ds dt \xrightarrow[T \to \infty]{} 0~.
    \end{split}
\end{equation}
This concludes the proof of Lemma \ref{lem:integrable_reg}.\hfill $\square$

\subsubsection{Interpretation of the Assumption \eqref{eq:assump_reg_maintext}}
\label{app.interpret_reg}
Below, we will illustrate the assumption \eqref{eq:assump_reg_maintext}
\begin{equation}
    Q := \int_D \int_0^\infty \big ( |\Delta\Thetab_t(\thetab) | + U_t(\thetab)^2 \big ) e^{C_1 (U_t(\thetab) + \bar{U}_t)} dt \mu_0(d \thetab) < \infty,
\end{equation}
in Theorem \ref{thm:long_time_reg_maintext} by giving examples that satisfy this condition.

First, consider an example where $\exists \kappa > 0, \alpha > 1$ such that $\forall \thetab \in \supp \mu_0$ and $\forall~t>0$, 
\begin{equation}
    |\Delta\Thetab_t(\thetab) | < \kappa (t+1)^{-\alpha},
\end{equation}
that is, all characteristic flows share a uniform asymptotic convergence rate on the order of $t^{-\alpha}$. Then $\forall \thetab \in \supp \mu_0$,
\begin{equation}
    U_t(\thetab) = \int_t^\infty |\Delta\Thetab_s(\thetab)| ds \leq \frac{\kappa}{\alpha - 1} (t+1)^{- ( \alpha - 1)}
\end{equation}
and thus
\begin{equation}
    \bar{U}_t \leq \frac{\kappa}{\alpha - 1} (t+1)^{- ( \alpha - 1)}~.
\end{equation}
Therefore,
\begin{equation}
    \begin{split}
        Q \leq & \int_D \int_0^\infty \big ( |\Delta\Thetab_t(\thetab) | + U_t(\thetab)^2 \big ) e^{C_1 (U_0(\thetab) + \bar{U}_0)} dt \mu_0(d \thetab) \\
        \leq & \int_0^\infty \Big ( \kappa (t + 1)^{- \alpha} + \big ( \frac{\kappa}{\alpha - 1} \big )^2 (t+1)^{-2(\alpha - 1)} \Big ) e^{\frac{2 C_1 \kappa}{\alpha - 1}} dt,
    \end{split}
\end{equation}
which is finite as long as $\alpha > \frac{3}{2}$. Thus,
\begin{proposition}
If $\exists \kappa > 0, \alpha > \frac{3}{2}$ such that $\forall \thetab \in \supp \mu_0$ and $\forall t\ge0$,
\begin{equation}
    |\Delta \Thetab_t(\thetab)|=|\Thetab_t(\thetab) - \Thetab_\infty(\thetab)| < \kappa (t+1)^{-\alpha},
\end{equation}
then the condition \eqref{eq:assump_reg_maintext} is satisfied.
\end{proposition}

Moreover, the assumption allows flexibility in having non-uniform convergence rate for different characteristic flows, $\Thetab_t(\thetab)$. Suppose that $\exists \kappa: \supp {\mu_0} \to \mathbb{R}_+$ and $\alpha > \frac{3}{2}$ such that $\forall \theta \in \supp \mu_0$,
\begin{equation}
    |\Delta\Thetab_t(\thetab) | < \kappa(\thetab) (t+1)^{-\alpha}~.
\end{equation}
Then
\begin{equation}
    U_t(\thetab) = \int_t^\infty |\Delta\Thetab_s(\thetab) | ds \leq \frac{\kappa}{\alpha - 1} (t+1)^{- ( \alpha - 1)}
\end{equation}
and so
\begin{equation}
\begin{split}
    Q \leq & \int_D \int_0^\infty \big ( |\Delta\Thetab_t(\thetab) | + U_t(\thetab)^2 \big ) e^{2 C_1 (U_0(\thetab))} dt \mu_0(d \thetab) \\
    \leq & \int_D \int_0^\infty \Big ( \kappa(\thetab) (t + 1)^{- \alpha} + \big ( \frac{\kappa(\thetab)}{\alpha - 1} \big )^2 (t+1)^{-2(\alpha - 1)} \Big ) e^{\frac{2 C_1 \kappa(\thetab)}{\alpha - 1}} dt \\
    \leq & C_2 \int_D \big ( \kappa(\thetab) + \kappa(\thetab)^2 \big ) e^{\frac{2 C_1 \kappa(\thetab)}{\alpha - 1}} \mu_0(d \thetab)~.
\end{split}
\end{equation}
Therefore,
\begin{proposition}
Suppose $\exists \alpha > \frac{3}{2}$ and a function $\kappa: \supp \mu_0 \to \mathbb{R}_+$, which satisfies 
\begin{equation}
    \int_D \Big ( \kappa(\thetab) + \kappa(\thetab)^2 \Big ) e^{\frac{2 C_1 \kappa(\thetab)}{\alpha - 1}} \mu_0(d \thetab) < \infty,
\end{equation} 
such that $\forall \thetab \in \supp \mu_0$, 
\begin{equation}
    |\Delta\Thetab_t(\thetab)|=|\Thetab_t(\thetab) - \Thetab_\infty(\thetab)| \leq \kappa(\thetab) (t+1)^{-\alpha}~.
\end{equation}
Then the condition \eqref{eq:assump_reg_maintext} is satisfied.
\end{proposition}

\subsubsection{Relationship between Theorem~\ref{thm:long_time_reg_maintext} and \cite{chizat2019sparse}}
\label{app:chizat}
As a comparison to our result, Chizat \cite[Theorem 3.8]{chizat2019sparse} shows that under assumptions including \eqref{eq:pd} as well as the uniqueness and sparseness of the global minimizer, an alternative type of particle gradient descent (with a different homogeneity degree in the loss function and under the conic metric, which give rise to gradient flow in Wasserstein-Fisher-Rao metric instead of Wasserstein metric) converges to the global minimizer for large enough $n$ (depending exponentially on $d$) with a uniform rate. This implies that in that setting, $\lim_{t \to \infty} \lim_{n \to \infty} n \| f_t^{(n)} - f_t \|_{\hat{\nu}}^2 = \lim_{n \to \infty} \lim_{t \to \infty} n \| f_t^{(n)} - f_t \|_{\hat{\nu}}^2 = 0$, $\mathbb{P}_0$-almost surely.

\section{Properties of the Minimizers of the Regularized Loss}
\label{app.prop_min}
First, under Assumption \ref{ass:shallownn}, i.e. in the shallow neural networks setting, define
\begin{equation}
     \hat{F}(\zb) = \int_\Omega f_*(\xb) \hat{\varphi}(\zb, \xb) \hat{\nu}(d \xb), \qquad \hat{K}(\zb, \zb') = \int_\Omega \hat{\varphi}(\zb, \xb) \hat{\varphi}(\zb', \xb) \hat{\nu}(d \xb)~.
\end{equation}
and
\begin{equation}
\label{eq:Vhatdef}
    \hat{V}(\zb,\mu) = -\hat{F}(\zb) + \int_D c' \hat{K}(\zb, \zb') \mu(d c', d \zb')~.
\end{equation}
We prove Proposition~\ref{th:min}, which we extend into:
\begin{proposition}
  \label{th:min2}
  Under Assumptions \ref{ass:shallownn}, \ref{ass:unit_1}, and \ref{ass:D_hat_compact}, the minimizers of the loss $\mathcal{L}(\mu)$ defined in~\eqref{eq:ermmeasure} 
  are all in the form
  \begin{equation}
    \label{eq:20}
    \mu_\lambda (dc,d\zb) = \delta_{c_\lambda}(dc) \hat \mu_+(d\zb)
    +\delta_{-c_\lambda}(dc) \hat \mu_-(d\zb)
  \end{equation}
  where $c_\lambda \ge0 $ and $\hat\mu_\pm\in \mathcal{P}(\hat D)$ satisfy
  \begin{equation}
  \label{eq:54}
  \begin{aligned}
    &\forall \zb \in \supp \hat \mu_- &:& \ \ -\hat F(\zb) + c_\lambda\int_{
      \hat D} \hat K(\zb,\zb')
    \left(\hat \mu_+(d\zb')-\hat \mu_-(d\zb')\right)= \lambda c_\lambda, \\
    &\forall \zb \in \supp \hat \mu_+ &:& \ \ -\hat F(\zb) + c_\lambda\int_{
      \hat D} \hat K(\zb,\zb')
    \left(\hat \mu_+(d\zb')-\hat \mu_-(d\zb')\right) =-
    \lambda c_\lambda,\\
    &\forall \zb \in \hat D &:& \ \
    \Big|-\hat F(\zb)+  c_\lambda \int_{
      \hat D} \hat K(\zb,\zb')
    \left(\hat \mu_+(d\zb')-\hat \mu_-(d\zb')\right)\Big|
    \le \lambda c_\lambda.
  \end{aligned}
\end{equation}
In addition, the constant $c_\lambda$ is unique and positive if
$\hat{F}(\zb)$ is not identically zero on $\hat D$, the closure of the
supports of $\hat \mu_\pm$ are disjoint (i.e.
$\overline{\supp \hat \mu_+ }\cap \overline{\supp \hat \mu_-
}=\emptyset$), and the function
  \begin{equation}
    \label{eq:22}
    f_\lambda = \int_D c\hat \varphi(\zb,\cdot) \mu_\lambda(dc,d\zb) =
    c_\lambda \int_{\hat D} \hat \varphi(\zb,\cdot)
    \left(\hat \mu_+(d\zb) -\hat \mu_-(d\zb) \right)
  \end{equation}
  is the same for all minimizers and satisfies
  \begin{equation}
    \label{eq:43}
    \tfrac14 \lambda^2 |c_\lambda|^{2} \hat K_M^{-1} \le 
    \left \|f_* - f_\lambda \right \|_{\hat{\nu}}^2, \qquad 
    \|f_* - f_\lambda \|_{\hat{\nu}}^2 + \lambda |c_\lambda|^2
    \le \lambda |\gamma_1(f_*)|^2.
  \end{equation}
  where $\hat K_M = \max_{\zb\in\hat D} \left \|\hat \varphi(\zb,\cdot) \right \|_{\hat{\nu}}^2 = \max_{\zb\in\hat D} \hat K(\zb,\zb)$.
\end{proposition}

\begin{remark}
Note that the proposition automatically implies that
$\gamma_1(f_\lambda) \le \gamma_1(f_*) < \infty$. It also implies that 
\begin{equation}
  \label{eq:122}
  \int_D |c|^q \mu_\lambda (dc,d\zb) =|c_\lambda|^q= |\gamma_\lambda|^q_{\text{TV}}\le  |\gamma_1(f_*)|^q
  \qquad \forall q\in\RR_+
\end{equation}
where $\gamma_\lambda = \int_{\RR} c \mu_\lambda(dc,\cdot)$.  
Finally note that the proposition holds if we replace the empirical loss by the population loss.
\end{remark}

\noindent
\textit{Proof:} 
The fact that this loss can only be minimized by minimizers follows from the compactness of the sets $ \{\mu\in \mathcal{P}(D)\ :\ \mathcal{L}(\mu) \le u, u\in \RR\}$.  The minimizers of $\mathcal{L}(\mu)$ must satisfy the following Euler-Lagrange equations~\cite{serfaty2014coulomb}:
\begin{equation}
  \label{eq:ELb}
  \forall (c,\zb) \in D\quad : \quad -c\hat F(\zb) + c \int_{ D}  c' \hat K(\zb,\zb')
  \mu(dc',d\zb') + \tfrac12 \lambda |c|^2 \equiv c\hat V(\zb) +
  \tfrac12\lambda |c|^2\ge \bar V~,
\end{equation}
with equality on the support of $\mu$ and where $\bar V$ is the expectation of the left hand side with respect to~$\mu(dc, d\zb)$. Minimizing the left hand side of~\eqref{eq:ELb} over $c$ at fixed $\zb$, we deduce that
\begin{equation}
  \label{eq:14}
  \forall \zb \in \hat D\quad : \quad
  \min_c\left(c\hat V(\zb) + \tfrac12\lambda |c|^2\right) 
  \ge \bar V~,
\end{equation}
with equality for $\zb$ in the support of
$\hat \mu=\int_{\RR}\mu(dc,\cdot)$.  This means that for any
$\zb\in \supp\hat \mu$, there can only be one $c=c(\zb)$ in
$\supp \mu$, with $c(\zb)$ satisfying the Euler-Lagrange equation
associated with~\eqref{eq:14}
\begin{equation}
  \label{eq:21}
  \hat V(\zb) + \lambda  c(\zb)=0 \quad \Leftrightarrow
  \quad 
  \hat V(\zb) = -\lambda  c(\zb)
\end{equation}
If we insert this equality back in
$c(\zb) \hat V(\zb) + \frac12\lambda |c(\zb)|^2= \bar V$, we deduce that
$|c(\zb)|=c_\lambda$, with the constant $c_\lambda$ related to
$\bar V$ as
\begin{equation}
  \label{eq:31}
  - \frac{1}{2} \lambda  |c_\lambda|^{2} = \bar V ~,
\end{equation}
and furthermore, $\forall \zb \in \supp {\hat{\mu}}$,
\begin{equation}
\label{eq:V_hat_c_lambda}
    \hat{V}(\zb) = \begin{cases} - \lambda c_\lambda & \text{if } c(\zb) = c_\lambda \\
    \lambda c_\lambda & \text{if } c(\zb) = - c_\lambda \end{cases}~.
\end{equation}
These considerations imply that the minimizer must be of the form~\eqref{eq:20}, and if we combine~\eqref{eq:14} and~\eqref{eq:31} and evaluate the minimum on $c$ explicitly we deduce that $\hat \mu_\pm$ and $c_\lambda$ must satisfy the equations in~\eqref{eq:54}. It is also clear from~\eqref{eq:54} that we must have $\overline{\supp \hat \mu_+ }\cap \overline{\supp \hat \mu_-}=\emptyset$: indeed if there was a point $\zb \in \overline{\supp \hat \mu_+ }\cap \overline{\supp \hat \mu_-}$, then at that point $\hat V(\zb)$ would be discontinuous, which is not possible since this function is continuously differentiable for any $\mu$ by our assumptions on $\hat \varphi$. Finally, to show that we must have that $c_\lambda>0$ if $F(\zb) $ is not identically zero on $\hat D$, note that if $c_\lambda =0$, \eqref{eq:ELb} reduces to
\begin{equation}
  \label{eq:53}
  \forall (c,\zb) \in D\quad : \quad -c \hat F(\zb) + \tfrac12 \lambda |c|^2 \ge0
\end{equation}
which can only be satisfied if $\hat F(\zb)=0$.

To show that $c_\lambda$ and the function in~\eqref{eq:22} are unique, let $\mu_\lambda$ and $\mu'_\lambda$ be two different minimizers and consider
\begin{equation}
  \label{eq:93}
  f_\lambda = \int_D c\hat \varphi(\zb,\cdot) \mu_\lambda(dc,d\zb)
  \quad \text{and} \quad
  f'_\lambda = \int_D c\hat \varphi(\zb,\cdot) \mu'_\lambda(dc,d\zb)
\end{equation}
Let us evaluate the loss on $a \mu_\lambda + (1-a) \mu'_\lambda\in\mathcal{P}(D)$ with $a\in[0,1]$. By convexity of $\mathcal{E}_\lambda$ we have
\begin{equation}
  \label{eq:94}
  \mathcal{L}(a \mu_\lambda + (1-a)
  \mu'_\lambda) \le a \mathcal{L}(\mu_\lambda) + (1-a)
  \mathcal{L}(\mu'_\lambda) = \mathcal{L}(\mu_\lambda)
  = \mathcal{L}(\mu'_\lambda)
\end{equation}
Since $a \mu_\lambda + (1-a) \mu'_\lambda$ cannot have a lower loss than this minimum, we must have equality in~\eqref{eq:94}, which reduces to
\begin{equation}
  \label{eq:95}
  \begin{aligned}
    & \left \| f_* -a f_\lambda - (1-a) f'_\lambda
    \right \|_{\hat{\nu}}^2 +  a \lambda |c_\lambda|^2 +
     (1-a) \lambda |c'_\lambda|^2\\
    & = \left \| f_* -f_\lambda \right \|_{\hat{\nu}}^2 +  \lambda
    |c_\lambda|^2
    \\
    & = \left \| f_* - f'_\lambda \right \|_{\hat{\nu}}^2 +  \lambda
    |c'_\lambda|^2~,
  \end{aligned}
\end{equation}
where $c_\lambda $ and $c_\lambda'$ are associated with $\mu_\lambda$ and $\mu_\lambda'$, respectively. Clearly these equations can only be fulfilled for all $a \in [0,1]$ if $c_\lambda=c_\lambda'$ and $f_\lambda = f_\lambda'$ ${\hat \nu}$-a.e. on $\Omega$.

To establish~\eqref{eq:43}, notice that if $\mu_\lambda$ is a minimizer and $f_\lambda $ is given by~\eqref{eq:22}, then we can derive from \eqref{eq:V_hat_c_lambda} that
\begin{equation}
  -\int_\Omega f_\lambda(\xb) f_*(\xb) \hat{\nu}(d \xb) + \left \|f_\lambda
  \right \|_{\hat{\nu}}^2 
  + \lambda  |c_\lambda|^2= 0.
\end{equation}
This gives, using Cauchy-Schwartz,
\begin{equation}
  \label{eq:45}
  \lambda  |c_\lambda|^2 = \int_\Omega f_\lambda(\xb) (f_*(\xb)-f_\lambda(\xb)) \hat{\nu}(d \xb)
  \le \left \|f_\lambda \right \|_{\hat{\nu}}
    \left \|f_*-f_\lambda \right \|_{\hat{\nu}}~.
\end{equation}
Now notice that
\begin{equation}
  \label{eq:57}
  \begin{aligned}
    \left \|f_\lambda \right \|_{\hat{\nu}}^2 &= c_\lambda^2\int _{\hat D\times \hat D}
    \hat K(\zb,\zb') \left(\hat \mu_+(d\zb) -\hat \mu_-(d\zb) \right)
    \left(\hat \mu_+(d\zb') -\hat \mu_-(d\zb') \right)\le 4c_\lambda^2
    \hat K_M~.
  \end{aligned}
\end{equation}
Using~\eqref{eq:57} in~\eqref{eq:45} and reorganizing gives the first inequality in~\eqref{eq:43}. To establish the second, let $\mu_*\in \mathcal{M}_+(D)$ be the measure that minimizes $\int_D |c| \mu(dc,d\zb)$ under the constraint that $f_* = \int_D c \hat \varphi( \zb,\cdot) \mu_*(dc,d\zb)$, so that $\int_D |c| \mu_*(dc,d\zb)=\gamma_1(f_*)$---the measure $\mu_*$ exists since we assumed that $f_* \in \mathcal{F}_1$. Evaluated on $\mu_*$, the loss is
\begin{equation}
  \label{eq:48}
  \mathcal{L} (\mu_*) = \lambda|\gamma_1(f_*)|^2.
\end{equation}
Any minimizer $\mu_\lambda$ of $\mathcal{L}(\mu)$ must do at least as well, i.e we must have
\begin{equation}
  \label{eq:58}
    \left \|f_* - f_\lambda \right \|_{\hat{\nu}}^2 + \lambda \int_D |c|^2
    \mu_\lambda(dc,d\zb)
    = \left \|f_* - f_\lambda \right \|_{\hat{\nu}}^2 + \lambda |c_\lambda|^2
    \le \lambda |\gamma(f_*)|^2.
\end{equation}
This establishes the second inequality in~\eqref{eq:43}.  \hfill $\square$

\section{Analytical Calculations of the Resampling Error}
\label{app.mc_value}
Derivations similar to the one presented here can be found in \cite{leroux07a, cho2009kernel, bach2017breaking}.
In the setting of ReLU without bias on unit sphere, we take $\hat{D} = \Omega = \mathbb{S}^d \subseteq \mathbb{R}^{d+1}$, $\hat{\varphi}(\zb, \xb) = \max(\langle \zb, \xb \rangle, 0)$, and $\nu$ is equal to the uniform measure on $\mathbb{S}^d$. In this case,
\begin{equation}
\begin{split}
\hat{K}(\zb, \zb') 
=  \int_\Omega \hat{\varphi}(\zb, \xb) \hat{\varphi}(\zb', \xb) \nu(d \xb) = \frac{1}{2 (d+1) \pi} (\sin \alpha + (\pi - \alpha) \cos \alpha),
\end{split}
\end{equation}
with $\alpha$ being the angle between $\zb$ and $\zb'$, and
\begin{equation}
    \begin{split}
        \int_\Omega |\hat{\varphi}(\zb, \xb)|^2 \nu(d \xb) = \frac{1}{2} \int_\Omega (\langle \xb,  \zb \rangle)^2 \nu(d \xb) = \frac{1}{2(d+1)}
    \end{split}
\end{equation}
Thus, taking $\mu_*$ to be the measure representing the teacher network, $\mu_* = \frac{1}{m_t} \sum_{i=1}^{m_t} \delta_{\zb_i}(d \zb) \delta_1(d c)$, we have
\begin{equation}
\begin{split}
    \int_D \|\varphi(\thetab,\cdot)\|_{\nu}^2 \mu_*(d\thetab)  =& \int_D \int_\Omega |\varphi(\thetab, \xb)|^2 \nu(d \xb) \mu_*(d \thetab) \\
    =& \int_D  \frac{c^2}{2(d + 1)} \mu_*(d \thetab) \\
    =& \frac{1}{2(d + 1)} 
\end{split}
\end{equation}
On the other hand,
\begin{equation}
\begin{split}
    \|f_*\|_{\nu}^2 = & \int_\Omega \Big | \int_D \varphi(\thetab, \xb) \mu_*(d \thetab) \Big |^2 \nu(d \xb) \\
    =& \int_D \int_D c c' \hat{K}(\zb, \zb') \mu_*(d \thetab) \mu_*(d \thetab') \\
    =& \frac{1}{m_t^2} \sum_{i, j = 1}^{m_t} \hat{K}(\zb_i, \zb_j)
\end{split}
\end{equation}
In the experiments described in the main text, we take $m_t = 2$, and $\zb_1$ and $\zb_2$ are initialized with a fixed random seed such that their angle, $\alpha_{12}$, equal to $1.766$. Thus, 
\begin{equation}
     \|f_*\|_{\nu}^2 =\frac{1}{4 (d+1) \pi} (0 + \pi) +  \frac{1}{4 (d+1) \pi} (\sin \alpha_{12} + (\pi - \alpha_{12}) \cos \alpha_{12}) \approx 0.012
\end{equation}
Together, we get a numerical value of the RHS of \eqref{eq:mcbound_lem33} if we replace $\mu_\infty$, $f_\infty$ and $\hat{\nu}$ by $\mu_*$, $f_*$ and $\nu$, respectively.

\end{document}